\documentclass[a4paper,fleqn]{cas-sc}

\usepackage[numbers]{natbib}

\usepackage{amssymb}
\usepackage{amsthm}
\usepackage{appendix}
\usepackage{array}
\usepackage{bm}
\usepackage{booktabs}
\usepackage{enumitem}
\usepackage{etoolbox}
\usepackage{float}
\usepackage{framed}
\usepackage{graphicx}
\usepackage{hyperref}
\usepackage{lineno}
\usepackage{mathtools}
\usepackage{mdframed}
\usepackage{multirow}
\usepackage{nicematrix}
\usepackage{setspace}
\usepackage{soul}
\usepackage{tabularx}
\usepackage[most]{tcolorbox}
\usepackage[colorinlistoftodos]{todonotes}
\usepackage{tikz}
\usetikzlibrary{decorations.pathreplacing}
\usepackage{xcolor}

\definecolor{revisionbrown}{RGB}{139,69,19}

\definecolor{jancolor}{RGB}{180,0,180}
\definecolor{alexcolor}{RGB}{0,0,255}

\newif\ifrevision
\revisiontrue

\ifrevision

\else

\fi

\DeclareRobustCommand{\corr}[1]{#1}

\newtheorem{theorem}{Theorem}[section]

\definecolor{myblue}{RGB}{0,0,255}
\definecolor{mygreen}{RGB}{28,172,0} 
\definecolor{myred}{RGB}{255,0,0}
\definecolor{mycyan}{RGB}{0,255,255}
\definecolor{mylilas}{RGB}{170,55,241}
\definecolor{myazure}{RGB}{0,175,225}
\definecolor{myorange}{RGB}{175,125,0}
\definecolor{mygray}{RGB}{150,150,150}
\definecolor{chatbackcolor}{rgb}{0.95,0.95,0.92}

\lstnewenvironment{mm}[1][]{%
  \lstset{
  language=Matlab,
  backgroundcolor=\color{chatbackcolor},
  showstringspaces=false,
  xleftmargin = 0.2cm,
  xrightmargin = 0.1cm,
  framexleftmargin=2pt,
  extendedchars=true,
  captionpos=b,
  basicstyle=\fontsize{7.5}{9.5}\selectfont\rmfamily,
  tabsize=2,
  keywordstyle=\textbf,
  stringstyle=\textit,
  numberstyle=\small,
  breakautoindent  = true,
  breakindent      = 2em,
  breaklines       = true,
  postbreak        = ,
  frame=single,
  keywordstyle=\color{blue},
  numbers=none,%
    numberstyle={\tiny \color{black}},
    numbersep=9pt, 
    emph=[1]{function,for,end,break},
    emphstyle=[1]\color{red}, 
  literate={ö}{{\"o}}1
           {ä}{{\"a}}1
           {ü}{{\"u}}1, #1
}}{}

\newcommand{\Sym}{\operatorname{Sym}}
\newcommand{\Symzero}{\operatorname{Sym}_0}
\newcommand{\Var}{\operatorname{Var}}
\newcommand{\id}{\operatorname{id}}
\newcommand{\Argmin}{\operatorname*{Argmin}}

\newcommand{\rt}[1]{\textcolor{black}{#1}}

\newcommand{\md}{\mathrm{d}}

\newcommand{\calI}{\mathcal{I}}

\newcommand{\calN}{\mathcal{N}}
\newcommand{\calV}{\mathcal{V}}
\newcommand{\calVh}{\mathcal{V}_h}
\newcommand{\calTh}{\mathcal{T}_h}
\newcommand{\bbP}{\mathcal{P}}
\newcommand{\bbV}{\mathbb{V}}

\renewcommand{\j}{j}
\renewcommand{\k}{k}
\renewcommand{\l}{\ell}
\renewcommand{\P}{\boldsymbol{P}}
\renewcommand{\S}{\boldsymbol{S}}
\renewcommand{\u}{\boldsymbol{u}}

\newcommand{\A}{A}
\newcommand{\Akq}{\A_\k^\q}
\newcommand{\Aq}{\hat{\A}_{\q}}
\newcommand{\B}{\mathbf{B}}
\newcommand{\C}{\boldsymbol{C}}
\newcommand{\De}{\mathbb{D}_{\Sigb}}
\newcommand{\DF}{\mathbb{D}_{\boldsymbol{P}}}
\newcommand{\ds}{\md s}
\newcommand{\dx}{\md \x}
\newcommand{\dxi}{\md \boldsymbol{\xi}}
\newcommand{\E}{\boldsymbol{E}}
\newcommand{\errN}{\epsilon_{\ts}^{\l}}
\newcommand{\F}{\boldsymbol{F}}
\newcommand{\Iphi}{I_{\varphi}}
\newcommand{\ii}{i}
\newcommand{\iset}{\ii \in \{1, \dots, \nn\}}
\newcommand{\jset}{\j \in \{1, 2, 3\}}
\newcommand{\Kel}{\mathbf{K}^\textnormal{el}}
\newcommand{\Kt}{\mathbf{K}_{\ts}}
\newcommand{\Ktti}{\widetilde{\mathbf{K}}_\ts}
\newcommand{\kset}{\k \in \{1, \dots, \nt\}}
\newcommand{\LL}{\mathbb{L}}
\newcommand{\m}{m}
\newcommand{\matlab}{MATLAB}
\newcommand{\mset}{\m \in \{1, 2, 3\}}
\newcommand{\N}{\mathbb N}
\newcommand{\Nkl}{N_\k^\p}

\newcommand{\nf}{n_f}
\newcommand{\nint}{n_{int}}
\newcommand{\nmesh}{37}
\newcommand{\nn}{n_n}
\newcommand{\np}{n_p}
\newcommand{\nq}{n_q}
\newcommand{\nt}{n_e}
\newcommand{\Omegah}{\Omega_h}
\newcommand{\omegaq}{\omega_{\q}}
\newcommand{\PP}{\mathbb{P}}
\newcommand{\Phii}{\Phi_{\ii}}
\newcommand{\Phiix}{\Phii(\x)}
\newcommand{\p}{p}
\newcommand{\phikl}{\varphi_{\k,\p}}
\newcommand{\phiklx}{\phikl(\x)}
\newcommand{\phiref}{\hat{\varphi}}
\newcommand{\phirefp}{\phiref_{\p}}
\newcommand{\pset}{\p \in \{1, \dots, \np\}}
\newcommand{\Qiso}{Q}
\newcommand{\q}{q}
\newcommand{\qset}{\q \in \{1, \dots, \nq\}}
\newcommand{\R}{\mathbb R}
\newcommand{\Sigb}{\boldsymbol{\Sigma}}
\newcommand{\sigb}{\boldsymbol{\sigma}}
\newcommand{\Tk}{T_{\k}}
\newcommand{\Tref}{\hat{T}}
\newcommand{\tmax}{t_\textnormal{max}}
\newcommand{\ts}{\tau}

\newcommand{\vh}{\boldsymbol{v}_h}
\newcommand{\vhjx}{\vj_h(\x)}
\newcommand{\vij}{\vj_{\ii}}
\newcommand{\vj}{v^{(\j)}}
\newcommand{\vv}{\boldsymbol{v}}
\newcommand{\vx}{\vv(\x)}

\newcommand{\x}{\boldsymbol{x}}
\newcommand{\xib}{\boldsymbol{\xi}}

\newcommand{\ttBd}{\texttt{B}}
\newcommand{\ttFd}{\texttt{F}}
\newcommand{\ttFDD}{\texttt{F\_{3D}}}
\newcommand{\ttU}{\texttt{U}}

\newlist{stepenum}{enumerate}{4}
\setlist[stepenum,1]{label=\textbf{\arabic*.}, ref=\arabic*}
\setlist[stepenum,2]{label=\textbf{(\roman*)}, ref=\theenumi.\roman*}
\setlist[stepenum,3]{label=\textbf{(\Alph*)},  ref=\theenumii.\Alph*}
\setlist[stepenum,4]{label=\textbf{(\alph*)},  ref=\theenumiii.\alph*}
\setlist[stepenum]{nosep, leftmargin=*}

\def\tsc#1{\csdef{#1}{\textsc{\lowercase{#1}}\xspace}}
\tsc{WGM}
\tsc{QE}
\tsc{EP}
\tsc{PMS}
\tsc{BEC}
\tsc{DE}

\newcommand{\fivevdots}{%
    \mathrel{\vcenter{\offinterlineskip
        \hbox{$.$}\kern2pt
        \hbox{$.$}\kern2pt
        \hbox{$.$}\kern2pt
        \hbox{$.$}\kern2pt
        \hbox{$.$}
    }}%
}
\newcommand{\tenvdots}{%
    \mathrel{\vcenter{\offinterlineskip
        \hbox{$.$}\kern2pt
        \hbox{$.$}\kern2pt
        \hbox{$.$}\kern2pt
        \hbox{$.$}\kern2pt
        \hbox{$.$}\kern2pt
        \hbox{$.$}\kern2pt
        \hbox{$.$}\kern2pt
        \hbox{$.$}\kern2pt
        \hbox{$.$}\kern2pt
        \hbox{$.$}
    }}%
}

\begin{document}
\let\WriteBookmarks\relax
\def\floatpagepagefraction{1}
\def\textpagefraction{.001}

\shorttitle{Mathematical analysis and vectorized FEM implementation of logarithmic-strain-based additive finite elastoplasticity}

\shortauthors{M. Frost et~al.}

\title [mode = title]{
Mathematical analysis and vectorized FEM implementation of logarithmic strain space framework for finite strains: a case study on elastoplasticity} 

\author[1]{M Frost}[
                        orcid=0000-0001-7511-2910]

\ead{mfrost@it.cas.cz}

\credit{Conceptualization, Methodology, Validation, Writing – Original Draft, Supervision}

\affiliation[1]{organization={Institute of Thermomechanics, Czech Academy of Sciences},
    addressline={Dolej\v{s}kova 5}, 
    city={Prague},
    postcode={182~00}, 
    country={Czechia}}

\author[2]{M Kru\v{z}\'{i}k}[orcid=0000-0003-1558-5809]
\ead{kruzik@utia.cas.cz}

\credit{Conceptualization, Methodology, Validation, Writing – Original Draft} 

\affiliation[2]{organization={Institute of Information Theory and Automation, Czech Academy of Sciences},
    addressline={Pod Vod\'{a}renskou v\v{e}\v{z}\'{i} 4}, 
    city={Prague},
    postcode={182~00}, 
    country={Czechia}}

\author[1,3]{A Moskovka}[orcid=0000-0003-0091-151X]
\cormark[1]
\ead{moskovka@it.cas.cz}
\credit{Software, Validation, Data curation, Writing - Original draft preparation, Visualization}

\affiliation[3]{organization={Faculty of Applied Sciences, University of West Bohemia in Pilsen},
    addressline={Technick\'{a} 8},
    city={Plze\v{n}},
    postcode={301~00}, 
    country={Czechia}}

\author[2]{J Valdman}[orcid=0000-0002-6081-5362]
\ead{valdman@utia.cz}
\credit{Methodology, Software, Writing – Review \& Editing}

\cortext[cor1]{Corresponding author}

\begin{abstract}
We investigate the mathematical and computational aspects of the logarithmic strain space framework for finite-strain models in nonlinear continuum mechanics, focusing on the classical von Mises elastoplastic model with linear kinematic hardening. To overcome the lack of coercivity and lower-semicontinuity in the quadratic Hencky energy, we introduce higher-gradient regularization and prove the existence of energetic solutions for rate-independent quasistatic evolution under time-dependent Dirichlet loading. Complementarily, we develop a vectorized MATLAB finite-element implementation of the framework leveraging page-wise matrix functions and operations, and our custom vectorization library. Computational benchmarks show that local tensor operations are significantly accelerated, the logarithmic-strain approach does not substantially increase computational costs, and the implementation robustly handles large rotations and deformations.
\end{abstract}



\maketitle

\section{Introduction}
\label{sec:intro}

Finite strain models of materials are essential in engineering since they allow us to analyze the mechanical response of structures or components experiencing large displacements, rotations, or geometric changes. The relationship between stress and strain then becomes highly nonlinear, and the classical small strain theory no longer applies. The logarithmic strain space approach developed by Miehe, Apel and Lambrecht in \cite{MiApLa2002} is a well-established route for accomplishing the transition from a linear constitutive model to a nonlinear one. In the original paper, the theory was applied to the classical small-strain plasticity, but it was later successfully employed for the development of finite-strain models of thermoviscoplasticity \cite{MiMeDiGoSch2011}, gradient plasticity \cite{anguige2018,kruzik2021,MiWeAl2014}, gradient-enhanced damage--plasticity \cite{HoJi2013, HoBrReSi2022, FrMeSt2023},  or, recently, martensitic phase transformations \cite{ChaWoKa2024, MHVKJSF2025}. The main appeal lies in its modularity, allowing constitutive structures of geometrically linear models to be incorporated into the logarithmic strain space, accompanied by pre- and postprocessing phases that are purely geometric and independent of the constitutive structures themselves. This brings significant benefits for the implementation of such constitutive models into the finite element method (FEM), which is the key instrument for computational structural analysis in engineering.

Although this flexible approach is frequently adopted and implemented in finite element codes, rigorous mathematical existence theories for the resulting additive logarithmic-strain models are much less developed. This is largely because the elastic energy inherited from their linear constitutive siblings is generally neither polyconvex nor quasiconvex and lacks the standard coercivity properties (used in nonlinear elasticity). From the algorithmic standpoint, working in logarithmic strain space requires the computation of generalized logarithmic strain tensors and the corresponding transformation tensors. This entails repeated evaluations of matrix logarithms (see, e.g., \cite{AlHi2012, AlHiRe2013}) and spectral decompositions, and the use of nested loops for repeated computations on objects arising from finite element discretization, which brings an additional computational burden. Over the past two decades, vectorization techniques have been developed to make the operations on large vectors, matrices, or tensors more efficient. However, to the best of the authors' knowledge, a comprehensive algorithmic framework dealing with the numerical aspects of the FEM implementation of the additive logarithmic-strain approach in a concise and unified manner has not been presented yet.

To fill these gaps, the present work provides mathematical analysis and a vectorized FEM implementation of a finite strain version of von Mises elastoplasticity with linear kinematic hardening. The particular choice of the constitutive model was motivated by our aims stemming from applied mathematics and computational practice: it has a sufficiently simple and well-known constitutive structure that allows us to focus on the mathematical and numerical treatment of the model, yet it features the material nonlinearity typical of many real engineering materials; hence, the main challenges for mathematical analysis are preserved, and the constitutive core can still serve as a useful template for applying the approach to more complex material models.

Hence, as the first contribution of the present work, we provide a mathematically rigorous analysis of a suitably regularized version of the model and prove the existence of energetic solutions under time-dependent Dirichlet loading. To overcome the mathematical issues raised above, we employ powerful tools from the calculus of variations, which are now able to treat quasistatic evolutionary problems as well (e.g., \cite{ortiz-repetto, mielke1}) and resort to gradient regularization, where the stored energy density also depends, e.g., on the gradients (of various orders) of the deformation and/or internal variables. \rt{Mathematically, this ensures model compactness, which is essential for proving the existence of solutions via the direct method \cite{ball81}. Such a constitutive material law, referred to as non-simple, was first introduced by Toupin \cite{toupin} and has since been extensively developed, e.g., \cite{greenRivlin, forest, mierouNum, forest2}.}

As the second main contribution of this work, we develop, describe, and implement a fully vectorized version of the finite strain model in the \matlab{} computational environment. In scripting languages, explicit loops over all elements, nodes, or quadrature points of a FEM geometric discretization are typically computationally expensive. By exploiting efficient implementations of linear and multilinear algebra, vectorization techniques enable replacing such loops with operations on large vectors or matrix-like objects (called arrays in \matlab{}). The gain in computational efficiency (related to, e.g., elimination of interpreter overhead, exploiting optimized libraries and efficient processor and cache utilization) generally comes at the expense of increased memory consumption. Vectorization techniques have been successfully applied to second-order linear partial differential equations; see, for example, \cite{Koko2007, Ko2016, WeSo2018}, as well as to selected nonlinear problems \cite{CeSyVa2019, Ro2017}. They have proved particularly effective in finite-element assembly. Vectorized assembly procedures have been developed for standard P1 finite elements \cite{FuPrWi2011, RaVa2013}, edge elements \cite{AnVa2015}, primal hybrid finite elements \cite{MaPoVaAc2025, mallesham2026}, and higher-order adaptive finite elements \cite{InPr2023}, among others. A general vectorization framework applicable to several programming languages was presented in \cite{CuJaSc2016}. Here we combine vectorization strategies at both the global and local levels of the assembly process. We also benefit from vectorizing numerical spectral decomposition recently introduced in \matlab{}. This eliminates the need for analytical eigenvalue calculations based on Cardano's formulas, which are known to suffer from numerical instability when eigenvalues coincide or are close to one another, which may happen in continuum mechanics.

The paper is organized as follows. Section~\ref{sec:gen-frame} introduces the weak formulation of the quasistatic elastoplastic boundary value problem and summarizes the methodology of the logarithmic strain space approach presented in \cite{MiApLa2002}. This is employed in Section~\ref{sec:maths} to formulate the boundary value problem for von Mises elastoplasticity with kinematic hardening in the logarithmic strain space. We also focus on its mathematical analysis there -- we define the energetic solution concept in the spirit of \cite{RIS}, and we discuss the main pitfalls of the additive logarithmic strain formulation with respect to this concept. Section~\ref{sec:discretization} covers the problem's discretizations, introduces the general computational scheme, and briefly summarizes the backbone of the FEM code developed in \cite{CeSyVa2019}. The particular, fully vectorized FEM implementation obtained by extending that code is described in Section~\ref{sec:numerics}. In this respect, we employ the vectorized linear algebra framework described in \cite{MoRaVaVa2026}, whose development has benefited from implicit expansion and page-wise matrix operations in \matlab{}. We provide details on the algorithmization of all new features, including the MATLAB code, which is also available for free download from a code repository \cite{code-MatlabFileExchange}. Finally, we assess the performance of the vectorized code and demonstrate its computational capabilities on two real-world inspired examples in Section~\ref{sec:simulations}.

\section{Additive Lagrangian approach to finite plasticity developed by Miehe and coworkers} \label{sec:gen-frame}

In this section, we briefly introduce the elastoplastic quasistatic boundary value problem and logarithmic strain space approach for finite plasticity developed in \cite{MiApLa2002}. We summarize the main facts and refer the reader to the dedicated literature for details.

\subsection{Quasistatic boundary value problem of elastoplasticity}

\rt{Let us consider a quasistatic, i.e., inertia effects not considered, mechanical problem of a homogeneous solid body occupying a domain $\Omega \subset \mathbb{R}^3$. Its Lipschitz boundary is decomposed as $\partial\Omega = \bar{\Gamma}_\textnormal{D} \cup \bar{\Gamma}_\textnormal{N}$, where $\Gamma_\textnormal{D}$ and $\Gamma_\textnormal{N}$ are open and disjoint. A prescribed displacement function ${\boldsymbol{u}}_{\rm D}$ specifies the Dirichlet boundary conditions on $\Gamma_\textnormal{D}$. We denote the unknown displacement field by ${\boldsymbol{u}}$ and the space of suitable test functions vanishing on $\Gamma_\textnormal{D}$ by $\calV$. The general loading consists of prescribed volume forces $\boldsymbol{f}_{\rm V}$ in $\Omega$ and surface tractions $\boldsymbol{f}_{\rm T}$ on $\Gamma_\textnormal{N}$.}

For implementation in FEM, the principle of virtual work (the weak form of the momentum balance) is essential. Denoting by $\hat{\cdot}$ the current counterparts of the objects defined above in the reference configuration, we can formulate the virtual work principle in the current configuration as 
\begin{equation} \label{weak0SS}
\text{find } \hat{\boldsymbol{u}} \in \{\hat{\boldsymbol{u}}_{\rm D}\} + \hat{\calV} : \int_{\hat{\Omega}} \hat{\sigb}\big(\hat{\boldsymbol{\varepsilon}}(\hat{\nabla}\hat{\boldsymbol{u}}), \hat{\boldsymbol{\varepsilon}}^\textnormal{p}\big)
: \hat{\boldsymbol{\varepsilon}}(\hat{\nabla}\hat{\boldsymbol{v}}) \, \mathrm{d}\hat{x} 
= \int_{\hat{\Omega}}  \hat{\boldsymbol{f}}_{\rm V} \cdot \hat{\boldsymbol{v}} \, \mathrm{d}\hat{x} 
+ \int_{\hat{\Gamma}_{{\rm N}}} \hat{\boldsymbol{f}}_{\rm T} \cdot \hat{\boldsymbol{v}} \, \mathrm{d}\hat{s}, \qquad \forall \, \hat{\boldsymbol{v}} \in \hat{\calV}.
\end{equation}
Above, $\hat{\sigb}$ denotes the Cauchy stress tensor, and 
\begin{equation}
\hat{\boldsymbol{\varepsilon}}({\nabla}{\hat{\boldsymbol{u}}}) = \frac{1}{2}\big( {\nabla}{\hat{\boldsymbol{u}}} + ({\nabla}{\hat{\boldsymbol{u}}})^\top\big)
\end{equation}
is the (infinitesimal) strain tensor. Both tensors are functions of the displacement gradient. The additional dependence of $\hat{\sigb}$ on the internal field $\hat{\boldsymbol{\varepsilon}}^\textnormal{p}$ (or even a set of internal fields) allows dealing with nonlinear path-dependent material response via the so-called internal state variables theory, e.g. \cite{COL-GUR, MAU-MUS-94a, HOR-BAM-10}. Indeed, in many engineering situations, the material response is not only nonlinear (as, e.g., in hyperelasticity), but also path-dependent; hence, the constitutive relation also depends on internal variables (and even on their gradients) that capture the material's microstructural evolution and influence the (total) macroscopic strain. Omitting the dependence on $\hat{\boldsymbol{\varepsilon}}^\textnormal{p}$ in \eqref{weak0SS} would lead to the conventional elasticity instead.

Hence, the model must be completed by specifying the constitutive response, i.e., by recovering stress from the (work-conjugate) current strain, $\hat{\boldsymbol{\varepsilon}}$, and its history encompassed in the internal variable field denoted by $\hat{\boldsymbol{\varepsilon}}^\textnormal{p}$. In the small-strain approximation, the additive decomposition of strains,
\begin{equation}\label{eq:ss-additivity}
    \hat{\boldsymbol{\varepsilon}} = \hat{\boldsymbol{\varepsilon}}^\textnormal{el} + \hat{\boldsymbol{\varepsilon}}^\textnormal{p},
\end{equation}
is a highly useful assumption that has enabled engineers to construct a wide range of constitutive relations. A typical example is von Mises elastoplasticity, where the internal variable, $\hat{\boldsymbol{\varepsilon}}^\textnormal{p}$, captures the evolution of plastic flow in the material at higher stress; see Appendix~\ref{app:von-mises}. More details on the conceptual numerical framework for the small-strain setting, which served as a starting point for the developments below, can be found in \cite{CeSyVa2019}.

When the \emph{geometric} nonlinearities also come into play, we can reformulate the momentum balance in the \emph{reference} configuration in terms of the first Piola--Kirchhoff (nominal) stress tensor, $\boldsymbol{P}$, which is work-conjugate to the deformation gradient,
\begin{equation}\label{eq:def-grad}
    \F = \nabla \boldsymbol{u} + \boldsymbol{I},
\end{equation}
and $\P = \det (\F) \hat{\sigb}\F^{-\top}$. The total Lagrangian formulation of the weak problem \eqref{weak0SS} in the reference configuration, $\Omega$, then takes the form
\begin{equation} \label{weak0FS}
(\mathrm{P}) \quad
 \text{find } \boldsymbol{u} \in \{\boldsymbol{u}_{\rm D}\} + \calV : \quad \int_{\Omega} \boldsymbol{P}(\F(\nabla\boldsymbol{u}), \E^\textnormal{p}) : \nabla\boldsymbol{v} \, \dx = 
\int_{\Omega} \boldsymbol{f}_{\rm V} \cdot \boldsymbol{v} \, \dx + 
\int_{\Gamma_{{\rm N}}} \boldsymbol{f}_{\rm T} \cdot \boldsymbol{v} \, \ds, \qquad \forall \, \boldsymbol{v} \in \calV.
\end{equation}

Considering again a suitable strain tensor, $\boldsymbol{E}(\F(\nabla \boldsymbol{u}))$, 
and some internal variable, $\boldsymbol{E}^\textnormal{p}$, capturing the \emph{material} nonlinearities (see \eqref{eqsEC} and \eqref{eq:log-plast-strain} below), the tensor $\P$ can be recognized as a stress-displacement operator emphasizing the implicit dependence of the first Piola stress on the internal variable(s) in case of internal state variable theory. Such a situation will be treated in the rest of the paper.

Now, the situation is more difficult, since the additivity assumption similar to \eqref{eq:ss-additivity} is generally invalid and constructing models with some (inelastic) contribution to the total deformation is not straightforward. Indeed, even for the family of plasticity models with one internal variable representing the plastic strain, many approaches can be found in the literature \cite{Lee1969, simo, Eterovic1990, Peric1992, Miehe1998, SimHug98CI, Lubarda2004, XiBrMe2006}. One of the popular ones is the logarithmic strain approach formulated in a FEM-convenient modular way in \citep{MiApLa2002}.

\subsection {Construction of finite strain models via the logarithmic strain approach} \label{subsec:miehe-construct}

\begin{figure}
\centering
\vspace{5mm}
\normalsize

\begin{align*}
& \hspace{15mm} \boldsymbol{u} \rightarrow \underbrace{ \big[\boldsymbol{\varepsilon} (\nabla \boldsymbol{u}), \{\boldsymbol{\varepsilon}^\textnormal{p}\}\big] \rightarrow  \fbox{CONSTITUTIVE MODEL} \rightarrow \big[\sigb, \mathbb{D}_{\boldsymbol{\varepsilon}}, \{\boldsymbol{\varepsilon}^\textnormal{p}_{\rm updated}\}\big]}_\text{small-strain realm} \\[5mm]
& \boldsymbol{u} \rightarrow \boldsymbol{F} \xrightarrow{\textbf{(A)}} \underbrace{\big[ \boldsymbol{E}\big(\boldsymbol{F}(\nabla \boldsymbol{u})\big),\{\boldsymbol{E}^\textnormal{p}\}\big] \rightarrow \fbox{CONSTITUTIVE MODEL}  \xrightarrow{\textbf{(B)}}
\big[\Sigb, \De, \{\boldsymbol{E}^\textnormal{p}_{\rm updated}\}\big]}_\text{logarithmic strain space} \xrightarrow{\textbf{(C)}}
\big[\boldsymbol{P}, \DF \big]
\end{align*}

\vspace{-2mm}
\caption{Schematic comparison of the data flow for constitutive models within the conventional small-strain (top) and the logarithmic-strain (bottom) computational frameworks. Inputs to the constitutive model are on the left of the central box, and outputs on the right. Symbols \textbf{(A)}, \textbf{(B)}, and \textbf{(C)} over the arrows denote the preprocessing, constitutive update, and postprocessing phases referred to in Fig.~\ref{fig:alg-setup}, respectively.}
\label{fig:const-mod-scheme}
\end{figure}

\vspace{5mm}

The approach for a simple, algorithmically smooth transition from a small-strain model towards a finite strain formulation was originally developed by Miehe et al.~\cite{MiApLa2002}. The main advantage is that it does not require substantial modifications to the constitutive core; instead, it extends the original computational workflow with preprocessing and postprocessing phases. The key idea is that the original constitutive algorithm operates on objects defined in the logarithmic strain space, and the transformations in the pre- and postprocessing phases link them to the weak formulation \eqref{weak0FS}; see \citep{MiApLa2002} for the theoretical underpinnings.

Fig.~\ref{fig:const-mod-scheme} shows the brief scheme of finite element data flow for an initial boundary value problem. The upper process diagram shows the situation for a constitutive model formulated in the small-strain realm, with (possibly a set of) internal variables enclosed by curly brackets, $\{\cdot\}$. In such a situation, the current and the reference configurations coincide, hence we omit the $\hat{\cdot}$ notation for brevity. The lower diagram shows the modifications made to adopt the logarithmic strain-space approach fitted to the weak formulation \eqref{weak0FS}. In the preprocessing phase denoted \textbf{(A)} (cf. Fig.~\ref{fig:alg-setup}), the logarithmic strain measure, $\E$, is computed from the displacement field, $\u$, via the deformation gradient, $\F$. Then, this strain measure and internal variables enter the constitutive law, which may retain the identical structure as in the small-strain formulation. It provides the generalized stress tensor, $\Sigb$, and the algorithmic stiffness tensor (also called consistent tangent modulus), $\De$, marked as \textbf{(B)}. In the postprocessing phase, \textbf{(C)}, these are transformed into the first Piola--Kirchhoff stress measure, $\P$, and the corresponding algorithmic stiffness tensor, $\DF$. \rt{The postprocessing transformation is purely geometric and uses specific projection tensors. Other computational procedures from the original small-strain model -- such as handling internal variables or resolving the constitutive response (time discretization and state updates) -- remain unchanged.}

In line with the original paper \citep{MiApLa2002} and with respect to clear and direct communication of the mathematical and numerical advances, we focus our study on the von Mises elastoplastic model with linear kinematic hardening of Prager type, so that the only internal variable considered in the small-strain model is identified as the plastic strain tensor, $\boldsymbol{\varepsilon}^\textnormal{p}$. We refer readers to the Appendix~\ref{app:von-mises} for a brief summary of this model. However, as documented by a rich literature \cite{MiWeAl2014, HoBrReSi2022, MiMeDiGoSch2011, MHVKJSF2025}, the framework \emph{a priori} does not constrain the type and number of internal variables. Hence, in principle, the curly brackets in Fig.~\ref {fig:const-mod-scheme} could contain more internal state variables of various types (scalars, tensors, etc.) if additional features are to be captured in the model (e.g., isotropic hardening, damage, phase transition), provided their transition to the logarithmic space is properly treated.

Below, we provide an overview of the key tensorial quantities used to extend the model to finite strains, following the approach illustrated in Fig.~\ref{fig:const-mod-scheme}. The complete model is presented in the next section. The (Lagrangian) right Cauchy--Green tensor, $\C$, is resolved from the deformation gradient introduced in \eqref{eq:def-grad}, as
\begin{equation} \label{eqsFC}    
    \C = \F^{\top} \F.
\end{equation}
The logarithmic (Hencky) strain, $\E$, takes the form
\begin{equation} \label{eqsEC}
    \E = \tfrac{1}{2} \ln \C.
\end{equation}
Defining the transformation tensors $\PP$ (4th order) and $\LL$ (6th order) as
\begin{equation} \label{eqsPPTL}
    \PP = \frac{\partial \E}{\partial \F}, \qquad \LL = \frac{\partial^2 \E}{\partial \F^2},
\end{equation}
allows us to recover the first Piola--Kirchhoff stress, $\P$, as
\begin{equation} \label{eqsP}
    \P = \Sigb:\PP,
\end{equation}
and the consistent tangent modulus $\DF$, as
 \begin{equation} \label{eqsDF}
    \DF = \PP^{\top} : \De : \PP + \Sigb : \LL,
\end{equation}
which are cornerstones for the postprocessing phase. Let us also recall the relation between $\P$ and the second Piola--Kirchhoff stress $\S$
\begin{equation} \label{eqsSsi}
    \S = \F^{-1} \P,
    \end{equation}
and the Cauchy stress    
\begin{equation} \label{eqsCauchy}
    \sigb = \frac{1}{\det \F}\, \P \F^{\top}.
\end{equation}

The original formulation for metal plasticity models is motivated by the equivalence between multiplicative decomposition of the total deformation gradient and additive decomposition of logarithmic Lagrangian strain in the case of coaxial total and plastic deformations \cite{MiApLa2002}. In the paper, the authors also showed reliability and effectiveness for anisotropic plastic response. Let us note that for metals and alloys, we can usually assume that the elastic strain is much lower than the total strain when plastic flow occurs. In \cite{FrMeSt2022}, the authors compared computational results from the multiplicative decomposition-based formulation and the logarithmic strain-based formulation of elastoplasticity with nonlinear isotropic hardening. It was found that ``the combination of large elongations and large shearing results in a significant difference in results''. However, the majority of the problems outlined in the paper ``are not dominant or decisive in practical simulations.'' The deviations were attributed to the mathematical properties of the model, namely the loss of ellipticity in the logarithmic strain space approach, cf. also \cite{neff2016}. \corr{The related mathematical issues are discussed in the next section.}

\section{Formulation and mathematical analysis of the von Mises elastoplastic problem with kinematic hardening in the logarithmic strain space} \label{sec:maths}

\subsection{Additive von Mises elastoplasticity with linear kinematic hardening in the logarithmic strain space}

Motivated by the conventional (rate-independent) von Mises elastoplasticity model summarized in Appendix~\ref{app:von-mises} and by the framework of \cite{MiApLa2002} presented in the previous section, we formulate the logarithmic strain-space adaptation of the model to the finite strain regime. Loosely speaking, this can be accomplished by the following purely formal substitutions in relations in Appendix~\ref{app:von-mises}:
\begin{equation} \label{fosu}
  \boldsymbol{\varepsilon} \rightarrow \E, \qquad \boldsymbol{\varepsilon}^\textnormal{p} \rightarrow \E^{\rm p}, \qquad {\psi} \rightarrow {\Psi}, \qquad {\phi} \rightarrow {\Phi}, \qquad \sigb \rightarrow \Sigb,
\end{equation}
with material parameters retained, and complemented by the relations from Section~\ref{subsec:miehe-construct}.

In a more formal way, we first recall that the Lagrangian logarithmic (Hencky) strain is
\begin{equation} \label{eq:Hencky}
    \E(\boldsymbol F)
    =
    \frac12 \log(\boldsymbol F^\top \boldsymbol F),
\end{equation}
and we define the plastic variable
\begin{equation}\label{eq:log-plast-strain}
    \boldsymbol E^\textnormal{p}:\Omega\to\Symzero(3),
\end{equation}
where $\Sym(3)$ denotes the space of symmetric $3\times3$ matrices and
\begin{equation*}
    \Symzero(3):=\big\{\boldsymbol A\in\Sym(3):\ {\rm tr}\,\boldsymbol A=0\big\}
\end{equation*}
is its subspace of trace-free (deviatoric) matrices; the trace-free constraint expresses plastic incompressibility, see below. The variable $\boldsymbol E^\textnormal{p}$ is interpreted as the logarithmic plastic strain. Let us note that the plastic metric can be reconstructed simply as $\C^\textnormal{p} = \exp(2\E^{\rm p})$; see \cite{MiApLa2002} for details.

Consider for a moment the multiplicative decomposition of the deformation gradient into an elastic part, $\boldsymbol{F}^{\rm el}$, and a plastic part, $\F^{\rm p}$, as $\F = \F^{\rm el}\F^{\rm p}$. Under the assumption of co-axiality of $\F^{\rm p}$ and $\F$, one gets from \eqref{eq:Hencky}
\begin{equation} \label{eq:add-ansatz}
     \E = \frac{1}{2} \ln ({\F^{\rm p}}^{\top}{\F^{\rm el}}^{\top}\F^{\rm el}\F^{\rm p}) = \E^{\rm el} + \E^{\rm p}.
\end{equation}
To arrive at a computationally attractive framework, the logarithmic strain approach in \cite{MiApLa2002} assumes the approximate validity of this relation beyond the coaxiality assumption.

Hence, we define the elastic logarithmic
strain as
\begin{equation}
    \E^\textnormal{el}(\F, \boldsymbol{E}^\textnormal{p}) = \E(\boldsymbol F)-\boldsymbol{E}^\textnormal{p}.
\end{equation}
Such a formulation is particularly attractive because the plastic variable
\(\boldsymbol E^\textnormal{p}\) belongs to the linear space \(\Symzero(3)\), and plastic
incompressibility 
\begin{equation}
    \det \F^\textnormal{p} = 1
\end{equation}
is represented by the (computationally much simpler) linear constraint
\begin{equation}
{\rm tr }\boldsymbol E^\textnormal{p}=0.
\end{equation}
The (volumetric) free-energy density of the material is then given by
\begin{equation}
    \Psi(\boldsymbol{E}, \boldsymbol{E}^\textnormal{p}) = {\Psi^{\rm el}(\boldsymbol{E} - \boldsymbol{E}^\textnormal{p})} + {\Psi^{\rm kin}(\boldsymbol{E}^\textnormal{p})} = \frac{1}{2} (\boldsymbol{E} - \boldsymbol{E}^\textnormal{p}) : \mathbb{C} : (\boldsymbol{E} - \boldsymbol{E}^\textnormal{p}) + \frac{1}{2} H^\mathrm{kin} \boldsymbol{E}^\textnormal{p} : \boldsymbol{E}^\textnormal{p}
\end{equation}
with $\mathbb{C}$ and $H^\mathrm{kin}$ being the conventional elastic stiffness tensor and the kinematic hardening coefficient, respectively.
The generalized thermodynamic forces (in the logarithmic space) take the standard form
\begin{align}
    \Sigb &= \frac{\partial \Psi}{\partial \boldsymbol{E}} = \mathbb{C} : (\boldsymbol{E} - \boldsymbol{E}^\textnormal{p})  \label{eq:log-stress}, \\
    \boldsymbol{X} &= \frac{\partial \Psi^{\rm kin}}{\partial \boldsymbol{E}^\textnormal{p}} = H^{\rm kin} \boldsymbol{E}^\textnormal{p}, \label{eq:back-stress}\\
    \boldsymbol{\Xi} &= -\frac{\partial \Psi}{\partial \boldsymbol{E}^\textnormal{p}} = \Sigb - \boldsymbol{X}.
\end{align}
Here, $\Sigb$ is the logarithmic stress, $\boldsymbol{X}$ is the backstress in the logarithmic strain space, and $\boldsymbol{\Xi}$ is the relative stress (the force driving the dissipation).

The classical yield condition and flow rule (i.e. the relation between the generalized forces and fluxes/velocities) can be equivalently reformulated \cite{han} in terms of the dissipation (pseudo)potential
\begin{equation} \label{PhiEp}
\Phi(\dot{\boldsymbol{E}}^{\rm p})  = 
\begin{cases}
\sqrt{2/3}\,Y\|\dot{\boldsymbol{E}}^{\rm p}\| & \text{if }\;{\rm tr}(   \dot{\boldsymbol{E}}^{\rm p}) = 0, \\
+\infty & \text{else},
\end{cases}
\end{equation}
where $Y$ denotes the (constant) yield stress, together with a thermodynamic selection criterion, which can be loosely stated as a minimization problem
\begin{equation} \label{eq:MPDP} 
   \operatorname*{min}_{\dot{\boldsymbol{E}}^{\rm p}}
\left\{ \dot{\Psi} + \Phi \right\}.
\end{equation}
Formulation \eqref{eq:MPDP} is called the minimum principle of dissipation potential in \cite{HaFi2008}. Its discretized version, the incremental energy minimization principle, is common in mechanics of materials \cite{ortiz-repetto, Petryk2003, MiLaGu2004, StuPe2013}, and for the model at hand provides the same governing equations as the principle of maximum dissipation; see details in \cite{HaFi2008}. It also perfectly fits into the mathematical theory of energetic solutions \cite{RIS}, which is applied below.

\subsection{Regularized boundary value problem with time-dependent Dirichlet data}

Let us first recall why a rigorous treatment of the boundary value problem requires some care. In nonlinear elasticity, the existence of equilibria is usually proved by the direct method of the calculus of variations: one minimizes the stored energy over admissible deformations, which requires \emph{coercivity} of the energy, so that minimizing sequences are bounded in a suitable Sobolev space, and \emph{weak lower semicontinuity} of the energy along such sequences. For integral functionals of the deformation gradient, weak lower semicontinuity is essentially equivalent to the quasiconvexity of the stored energy density in the sense of Morrey \cite{dacorogna}. Convexity in $\F$ is too restrictive, since it is incompatible with frame indifference and with the blow-up of the energy as $\det\F\to0^+$, while quasiconvexity is very difficult to verify. Ball \cite{ball77} therefore introduced \emph{polyconvexity}, i.e., convexity of the energy as a function of $\F$, its cofactor matrix and its determinant. Polyconvexity is easy to check; it is compatible with the blow-up as $\det\F\to0^+$, and together with a simple \corr{coercivity condition on the stored energy density (cf.~\eqref{eq:growth} below)} it guarantees the existence of minimizers.

For rate-independent processes such as elastoplasticity, a natural notion of weak solution is that of an \emph{energetic solution} \cite{mielke2, RIS}, which is defined by a global stability condition and an energy balance, see below, and does not require any differentiability of the energy with respect to the state. Energetic solutions are constructed from the time-incremental minimization problems, i.e., from the incremental energy minimization principle mentioned above. Each incremental problem is a minimization problem to which the direct method applies under the same coercivity and lower semicontinuity requirements as in elasticity, and the limit of vanishing time steps is obtained from uniform energy estimates and compactness arguments \cite{mami1, mielke2}. For finite-strain plasticity based on the multiplicative decomposition, this program has been carried out in \cite{mielke1, mami2, mierouNum}, where gradient terms of the plastic variables provide the necessary compactness.

In the logarithmic strain framework, the plastic variable $\boldsymbol E^\textnormal{p}$ takes values in the linear space $\Symzero(3)$, and the dissipation, the kinematic hardening energy and the gradient term of $\boldsymbol E^\textnormal{p}$ introduced below are convex. Hence, the plastic part of the problem does not create any difficulty. The main difficulty lies instead in the deformation variable. Even for fixed \(\boldsymbol{E}^\textnormal{p}\), the map
\begin{equation}\label{eq:stress-pot_growth}
    \boldsymbol F\mapsto
    \Psi^\textnormal{el}\!\left(
        \frac12\log(\boldsymbol F^\top \boldsymbol F)-\boldsymbol E^\textnormal{p}
    \right)
\end{equation}
is not, in general, polyconvex or quasiconvex \cite{ball77,dacorogna}. Consequently, the direct
method of the calculus of variations does not yield existence of a
minimizer without regularization.  This issue is
already emphasized by Miehe, Apel, and Lambrecht in \cite{MiApLa2002}, who note that their
logarithmic-strain framework is not polyconvex in the elastic range and
that existence of global incremental minimizers is therefore not
guaranteed for arbitrary elastic strains.

There is also a coercivity issue.  A purely quadratic Hencky energy,
for example
\begin{equation*}
    \overline{\Psi}^\textnormal{el}(\E)\sim |\E|^2,
\end{equation*}
grows only logarithmically in the singular values of \(\boldsymbol F\), see the definition \eqref{eq:Hencky}.
Thus, one cannot expect a standard bound of the form
\begin{equation}\label{eq:growth}
    \Psi^\textnormal{el}\!\left(\frac12\log(\boldsymbol F^\top \boldsymbol F)-\boldsymbol E^\textnormal{p}\right)
    \ge c_\textnormal{el}|\boldsymbol F|^{\corr{\varrho}}-C_\textnormal{el}
\end{equation}
\corr{with constants $c_\textnormal{el},C_\textnormal{el}>0$ and some exponent $\varrho>1$.}
This is why we introduce regularizing gradient terms multiplied by a ``small'' parameter $\varepsilon>0$ in \eqref{eq:regul-energy} below. This allows us to obtain both coercivity and compactness and, finally, the existence of an energetic solution for the corresponding rate-independent evolution. Another approach to these issues was considered in \cite{neff2015}.

Let \(\Omega\subset\R^3\) be a bounded domain and let
\(\Gamma_\textnormal{D}\subset\partial\Omega\) have positive surface measure. For simplicity and with respect to simulations in Section~\ref{sec:simulations}, let us consider the only
external loading given through a time-dependent Dirichlet boundary condition. For \(t\in[0,T]\), let
\begin{equation}
    \boldsymbol g_t:=\boldsymbol g(t,\cdot):\R^3\to\R^3
\end{equation}
be an orientation-preserving diffeomorphism\footnote{Let us note that the matrix logarithm of $\C = \F^\top\F$ alone cannot distinguish orientation-preserving from orientation-reversing gradients; hence, this must be checked independently.}, and prescribe
\begin{equation}
    \boldsymbol y(t,\boldsymbol x)
    =
    \boldsymbol g_t(\boldsymbol x)
    \quad\text{on }\Gamma_\textnormal{D}.
\end{equation}
Following the composition approach for time-dependent boundary conditions,
we write
\begin{equation}
    \boldsymbol y(t)
    =
    \boldsymbol g_t\circ\boldsymbol\upsilon(t),
\end{equation}
and impose the fixed boundary condition
\begin{equation}
    \boldsymbol\upsilon(t,\boldsymbol x)
    =
    \boldsymbol x
    \quad\text{on }\Gamma_\textnormal{D}.
\end{equation}
Define
\begin{equation}
    \boldsymbol A_t(\boldsymbol\xi)
    :=
    \nabla\boldsymbol g_t(\boldsymbol\xi), \quad     \boldsymbol F_t(\boldsymbol\upsilon)
    :=
    \boldsymbol A_t(\boldsymbol\upsilon)\nabla\boldsymbol\upsilon.
\end{equation}
The corresponding total logarithmic strain is
\begin{equation}
    \E_t(\boldsymbol\upsilon)
    :=
    \frac12
    \log\!\left(
        \boldsymbol F_t(\boldsymbol\upsilon)^\top
        \boldsymbol F_t(\boldsymbol\upsilon)
    \right).
\end{equation}
Choose constants \(r>3\), \(q\ge2\), \(\varepsilon>0\), \(\kappa>0\), and \(s>3r/(r-3)\). We define the regularized energy by
\begin{multline} \label{eq:regul-energy}
    \mathcal E_\varepsilon
    (t,\boldsymbol\upsilon,\boldsymbol E^\textnormal{p})
    :=
    \int_\Omega
    \Psi^\textnormal{el}\!\left(
        \E_t(\boldsymbol\upsilon)-\boldsymbol E^\textnormal{p}
    \right)\dx
    +
    \int_\Omega
    \left[
        \Psi^\textnormal{kin}(\boldsymbol E^\textnormal{p})
        +
        \frac{\kappa}{q}|\nabla\boldsymbol E^\textnormal{p}|^q
    \right]\dx
    \\
    +
    \varepsilon
    \int_\Omega
    \left[
        \frac1r|\nabla\boldsymbol\upsilon|^r
        +\frac1r|\nabla^2\boldsymbol\upsilon|^r
        +(\det\nabla\boldsymbol\upsilon)^{-s}
    \right]\dx.
\end{multline}
To avoid $\det \F \leq 0$, we set
\begin{equation}
    \mathcal E_\varepsilon
    (t,\boldsymbol\upsilon,\boldsymbol E^\textnormal{p})
    =
    +\infty \quad \textnormal{whenever} \quad     \det\nabla\boldsymbol\upsilon\le0 \quad \textnormal{on a set of positive measure}.
\end{equation}
The admissible state space is
\begin{equation}
    \mathcal Q
    :=
    \mathcal Y\times\mathcal P, \quad
\textnormal{where} \quad \mathcal Y
    := \left\{
        \boldsymbol\upsilon\in W^{2,r}(\Omega;\R^3):
        \boldsymbol\upsilon=\id
        \text{ on }\Gamma_\textnormal{D}
    \right\}, \quad \mathcal P
    :=     W^{1,q}(\Omega;\Symzero(3)).
\end{equation} Let also the dissipation potential $\Phi : \Symzero(3)\to[0,\infty)$ 
be convex, positively one-homogeneous, and satisfy
\begin{equation}
    c_\Phi|\boldsymbol Z|
    \le
    \Phi(\boldsymbol Z)
    \le
    C_\Phi|\boldsymbol Z|.
\end{equation}
We define the dissipation distance between two states $\ts$ and $\ts+1$ as
\begin{equation}
    \mathcal D\left(\boldsymbol E^\textnormal{p}_\ts,\boldsymbol E^\textnormal{p}_{\ts+1}\right)
    :=
    \int_\Omega
    \Phi\left(\boldsymbol E^\textnormal{p}_{\ts+1}-\boldsymbol E^\textnormal{p}_\ts\right)\dx.
\end{equation}

\noindent \textbf{Energetic formulation}

A pair $ t\mapsto (\boldsymbol\upsilon(t),\boldsymbol E^\textnormal{p}(t)) $ is called an \emph{energetic solution} \cite{RIS} if, for every \(t\in[0,T]\), it satisfies
the \emph{global stability} condition
\begin{equation}
    \mathcal E_\varepsilon
    (t,\boldsymbol\upsilon(t),\boldsymbol E^\textnormal{p}(t))
    \le
    \mathcal E_\varepsilon
    (t,\widehat{\boldsymbol\upsilon},\widehat{\boldsymbol E^\textnormal{p}})
    +
    \mathcal D(
        \boldsymbol E^\textnormal{p}(t),\widehat{\boldsymbol E^\textnormal{p}}
    )
\end{equation}
for every admissible competitor
\( (\widehat{\boldsymbol\upsilon},\widehat{\boldsymbol E^\textnormal{p}})\in\mathcal Q\),
and the \emph{energy balance}
\begin{equation}
\begin{aligned}
    &
    \mathcal E_\varepsilon
    (t,\boldsymbol\upsilon(t),\boldsymbol E^\textnormal{p}(t))
    +
    \Var_{\mathcal D}
    (\boldsymbol E^\textnormal{p};[0,t])
    =
    \mathcal E_\varepsilon
    (0,\boldsymbol\upsilon_0,\boldsymbol E^\textnormal{p}_0)
    +
    \int_0^t
    \partial_\tau
    \mathcal E_\varepsilon
    (\tau,\boldsymbol\upsilon(\tau),\boldsymbol E^\textnormal{p}(\tau))
    \,\md \ts.
\end{aligned}
\end{equation}
\corr{Here, $\Var_{\mathcal D}(\boldsymbol E^\textnormal{p};[0,t]):=\sup\sum_{j=1}^{N}\mathcal D\big(\boldsymbol E^\textnormal{p}(t_{j-1}),\boldsymbol E^\textnormal{p}(t_j)\big)$, where the supremum is taken over all partitions $0=t_0<t_1<\dots<t_N=t$, is the total dissipation on $[0,t]$.}
The integral on the right-hand side represents the work induced solely by
the time-dependent Dirichlet boundary condition. Namely, at fixed $(\boldsymbol\upsilon,\boldsymbol E^\textnormal{p})$, define $\boldsymbol L_t :=    \dot{\boldsymbol A}_t \boldsymbol A_t^{-1}$. Since $
    \partial_t
    \boldsymbol F_t(\boldsymbol\upsilon)
    =
    \boldsymbol L_t(\boldsymbol\upsilon)
    \boldsymbol F_t(\boldsymbol\upsilon),
$ we have
\begin{equation}
\begin{aligned}
    \partial_t
    \mathcal E_\varepsilon
    (t,\boldsymbol\upsilon,\boldsymbol E^\textnormal{p})
    =
    \int_\Omega
    D\Psi^\textnormal{el}\!\left(
        \E_t(\boldsymbol\upsilon)-\boldsymbol E^\textnormal{p}
    \right)
    :
    D\E
    \left(
        \boldsymbol F_t(\boldsymbol\upsilon)
    \right)
    \left[
        \boldsymbol L_t(\boldsymbol\upsilon)
        \boldsymbol F_t(\boldsymbol\upsilon)
    \right]
    \dx.
\end{aligned}
\end{equation}

The stability condition can be interpreted as: there is no admissible change of state whose energetic gain is larger than the dissipation needed to perform that change. Let us note that the current state is compared with \emph{every admissible} competitor, not only with infinitesimally close states.\footnote{That is an important feature of the energetic formulation and also one of the reasons why it can describe discontinuous jumps.} The energy equality is easy to interpret: all work done by external loads is spent either on dissipation or stored in the material.

We formulate and prove the theorem which establishes the existence of an energetic solution to the regularized model in Appendix~\ref{app:en-solution}. In Section~\ref{sec:simulations}, we carry out numerical simulations to illustrate the advantage of the modularity of the logarithmic strain framework. We directly adopt the standard formulation of the elastoplastic model from \cite{CeSyVa2019}, which corresponds to setting $\varepsilon = \kappa = 0$ in \eqref{eq:regul-energy}, i.e., we neglect the higher-order regularization. This significantly simplifies the code; e.g., it permits avoiding higher-order approximation of the plastic strain variable in the FEM implementation. However, strictly speaking, the existence or convergence of such an unregularized model is not covered by the theorem. A detailed discussion is left to Appendix~\ref{app:en-solution}.

\section{Discretizations and algorithmic setup}
\label{sec:discretization}

The standard tool for the numerical solution of boundary value problems \eqref{weak0FS} is the finite element method (FEM) \cite{Ciarlet2002, Strang2008, Bathe2006}. The weak (variational) formulation of the problem is restricted to a finite-dimensional space of piecewise polynomial functions defined on a mesh, i.e., on a partition of the body into simple elements (in this paper, hexahedra with trilinear shape functions), and the unknowns are the nodal values of the displacement field. In elastoplasticity, the FEM is combined with a discretization in time, so that at every time step a nonlinear algebraic system has to be solved, while the internal variables are updated locally at the quadrature points of the elements by a return-mapping algorithm \cite{simo2, han}.

\subsection{Temporal and spatial discretizations of the initial boundary value problem in the weak form}

The time-discretization leads to a generic time-incremental (rate-independent) problem
\begin{equation} \label{weak2}
(\mathrm{P}_\ts) \quad
 \text{find } \boldsymbol{u}_\ts \in \{\boldsymbol{u}_{{\rm D},\ts}\} + \calV : \quad \int_{\Omega} \P_\ts(\boldsymbol{u}_\ts)
: \nabla\boldsymbol{v} \, \dx
=
\int_{\Omega} \boldsymbol{f}_{{\rm V},\ts}\cdot \boldsymbol{v} \, \dx
+
\int_{\Gamma_{{\rm N}}} \boldsymbol{f}_{{\rm T},\ts}\cdot \boldsymbol{v} \, \ds,
\quad \forall \, \boldsymbol{v} \in \calV.
\end{equation}
The operator $\P_\ts$ encompasses the implicit dependence on the internal variable $\E^\textnormal{p}_\ts$, which is, however, resolved via the constitutive core (cf. \eqref{eq:MPDP}) based on the value of $\u_\ts$, too. Finally, we perform the spatial discretization of the time-discretized problem \eqref{weak2}. Together with the time incrementation of the time interval $[0,\tmax]$, we obtain the problem in the form:
\begin{equation} \label{weak3}
\begin{split}
(\mathrm{P}_\ts)_h:\quad
&\text{find } \boldsymbol{u}_{\ts,h} \in \{\boldsymbol{u}_{\textnormal{D},\ts,h}\} + \calV_h \text{ such that} \\
&\int_{\Omegah} \P_\ts \, (\boldsymbol{u}_{\ts,h})
: \nabla \vh \, \dx
=
\int_{\Omegah} \boldsymbol{f}_{{\rm V},\ts,h}\cdot \vh \, \dx
+
\int_{\Gamma_{\textnormal{N},h}} \boldsymbol{f}_{{\rm T},\ts,h}\cdot \vh \, \ds,
\quad \forall \, \vh \in \calVh,
\end{split}
\end{equation}
where $\calVh$ from \eqref{Vh} is a finite-dimensional approximation of $\calV$, see \eqref{vxV} below, and
$\boldsymbol{f}_{{\rm V},\ts,h}$, $\boldsymbol{f}_{{\rm T},\ts,h}$, $\boldsymbol{u}_{{\rm D},\ts,h}$ and $\Omegah$  are suitable approximations of $\boldsymbol{f}_{{\rm V},\ts}$, $\boldsymbol{f}_{{\rm T},\ts}$, $\boldsymbol{u}_{{\rm D},\ts}$, and $\Omega$, respectively.

\subsection{Finite element discretization and formulation of the algebraic problem} \label{subsec:FEMbasics}

\begin{table}[htbp]
    \centering
    \renewcommand{\arraystretch}{1.15}
    \begin{NiceTabular}{cl}[left-margin=12mm]
        \toprule
        Symbol & Description \\
        \midrule

        $\ii$ & Index of a global node or basis function ($\iset$) \\
        $\j,\m$ & Spatial component indices ($\j,\m \in \{1,2,3\}$) \\
        $\k$ & Element index ($\kset$) \\
        $\p$ & Index of a local node or basis function ($\pset$) \\
        $\q$ & Local quadrature-point index ($\qset$) \\
        $\ts$ & Discrete time-step index \\
        $\l$ & Newton's method iteration index \\
        \addlinespace[1.5ex]

        $\calTh$ & Set of elements (computational mesh) \\
        $\nt$ & Total number of elements \\
        $\calN$ & Set of mesh nodes \\
        $\nn$ & Total number of nodes (global scalar nodal basis functions) \\
        $\nf$ & Total number of boundary faces \\
        \addlinespace[1.5ex]

        $\Tref$ & Reference element \\
        $\np$ & Number of nodes (local scalar basis functions) per element \\
        $\Nkl$ & $\p$-th local node of the $\k$-th element \\
        \addlinespace[1.5ex]

        $\nq$ & Number of quadrature points per element \\
        $\nint$ & Total number of quadrature points ($\nint = \nt\nq$) \\
        $\Aq$ & $\q$-th quadrature point of the reference element $\Tref$ \\
        $\Akq$ & $\q$-th quadrature point of the $\k$-th element \\
        \addlinespace[1.5ex]

        $\phirefp$ & $\p$-th basis function of the reference element $\Tref$ \\
        $\phikl$ & $\p$-th local basis function of the $\k$-th element \\
        $\Phii$ & $\ii$-th global basis function \\
        \addlinespace

        \bottomrule

        \CodeAfter
        \begin{tikzpicture}[
            category brace/.style={
                decorate,
                decoration={brace, amplitude=4pt},
                thick
            },
            category label/.style={
                rotate=90,
                font=\footnotesize,
                align=center
            }
        ]
            \draw[category brace]
                ([xshift=2mm,yshift=1mm]9-|1)
                -- ([xshift=2mm,yshift=-2mm]2-|1)
                node[midway,xshift=-6mm,category label]
                {Indices};

            \draw[category brace]
                ([xshift=2mm,yshift=0.5mm]14-|1)
                -- ([xshift=2mm,yshift=-2.5mm]9-|1)
                node[midway,xshift=-6mm,category label]
                {Mesh objects\\and counts};

            \draw[category brace]
                ([xshift=2mm,yshift=0.5mm]17-|1)
                -- ([xshift=2mm,yshift=-1.5mm]14-|1)
                node[midway,xshift=-6mm,category label]
                {Element\\and nodes};

            \draw[category brace]
                ([xshift=2mm,yshift=0.5mm]21-|1)
                -- ([xshift=2mm,yshift=-2.5mm]17-|1)
                node[midway,xshift=-6mm,category label]
                {Quadrature};

            \draw[category brace]
                ([xshift=2mm,yshift=2.5mm]24-|1)
                -- ([xshift=2mm,yshift=-2mm]21-|1)
                node[midway,xshift=-6mm,category label]
                {Basis\\functions};
        \end{tikzpicture}
    \end{NiceTabular}
    \caption{Indices and notation used in the finite element discretization.}
    \label{tab:fem-notation}
\end{table}

In order to solve \eqref{weak3}, we employ FEM for the spatial discretization; an overview of related indexing system and the notation is in Table~\ref{tab:fem-notation}. Let $\calTh$ denote a partition of the computational domain $\Omegah$ into closed hexahedral elements. On this mesh, we use trilinear Q1 finite elements in the sense of Ciarlet \cite{Ciarlet2002}.\footnote{Extension to other common element types is straightforward; cf. \cite{CeSyVa2019}.} Each element is geometrically defined by its eight vertices, which belong to the global set of mesh nodes $\calN$.

We \corr{consider a test function and a test space} (in three dimensions) of the form:
\begin{equation} \label{vxV}
\vx = \big(v^{(1)}(\x), v^{(2)}(\x), v^{(3)}(\x)\big), \qquad \calV = \mathcal V^{(1)} \times \mathcal V^{(2)} \times \mathcal V^{(3)}.
\end{equation}
Denote by $\mathcal W_h^s$ the full scalar trilinear finite-element space, including basis functions associated with the Dirichlet boundary, and set $\mathcal V_h^{(j)}
:= \left\{ w_h \in \mathcal W_h^s : w_h = 0 \text{ on } \Gamma_{\mathrm{D},h}^{(j)} \right\}, \; j=1,2,3$, where $\Gamma_{\mathrm{D},h}^{(j)}$ denotes the discrete Dirichlet boundary for the $j$-th displacement component. We then approximate $\vx \in \calV$ in the finite-dimensional subspace
\begin{equation} \label{Vh}
\calV_h = \mathcal V_h^{(1)} \times \mathcal V_h^{(2)} \times \mathcal V_h^{(3)} \subset \calV.
\end{equation}
The full scalar finite-element space $\mathcal W_h^s$ is spanned by the global Lagrange basis functions
\begin{equation} \label{globalPhi}
    \Phiix \in \mathcal W_h^s, \quad \iset,
\end{equation}
such that
\begin{equation}
    \Phii(N_{\tilde{\ii}}) =
    \begin{cases}
        1, \;\;\; \tilde{\ii} = \ii, \\
        0, \;\;\; \tilde{\ii} \neq \ii,
    \end{cases}
\end{equation}
where $\nn$ denotes their number (equivalent to the number of mesh nodes). Homogeneous Dirichlet conditions are imposed by setting the coefficients associated with constrained displacement components to zero. Hence, any component of $\vh(\x)$ is given by a linear combination
\begin{equation} \label{linearcombsingle}
    \vhjx = \sum_{\ii=1}^{\nn} \vij \, \Phiix \, , \qquad \x \in \Omegah, \quad \jset .
\end{equation}
By $\phiklx$ we denote the $\p$-th local basis function on the $\k$-th element, which is the restriction of the corresponding $\Phiix$ to $\Tk$
\begin{equation} \label{localPhi}
    \phiklx = \Phiix\big|_{\Tk}, \quad \ii = \Iphi(\k,\p),
\end{equation}
where $\Iphi: \N \times \N \rightarrow \N$ is the standard local-to-global index map.

We consider a reference element $\Tref$ with vertices
\begin{equation} \label{Tref}
    \begin{split}
        &\hat{N}_1 = [-1, -1, -1], \qquad \hat{N}_2 = [\;\; 1, -1, -1], \qquad \hat{N}_3 = [-1, \;\; 1, -1], \qquad \hat{N}_4 = [\;\; 1, \;\; 1, -1], \\
        &\hat{N}_5 = [-1, -1, \;\; 1], \qquad \hat{N}_6 = [\;\; 1, -1, \;\; 1], \qquad \hat{N}_7 = [-1, \;\; 1, \;\; 1], \qquad \hat{N}_8 = [\;\; 1, \;\; 1, \;\; 1].
    \end{split}
\end{equation}
and denote its local coordinates by $\xib = (\xi_1,\xi_2,\xi_3)$. Then, any element $\Tk \in \calTh, \, \kset,$ is given by the isoparametric mapping, $\Qiso$, between $\Tref$ and $\Tk$ defined by
\begin{equation} \label{IsoTra}
    \Qiso_\k: \xib \in [-1,1]^3 = \Tref \rightarrow \x \in \Tk: \, \Qiso_{\k}(\xib) = \sum\limits_{\p=1}^{\np} \Nkl \, \phirefp(\xib),
\end{equation}
where $\Nkl$ denotes the $\p$-th node of the $\k$-th element of $\calTh$ and $\phiref_{\p}$ is the $\p$-th local basis function on $\Tref$. Using \eqref{IsoTra}, for any $\Tk$ we obtain the following relation:
\begin{equation} \label{philocref}
    \phiklx = \phirefp(\xib), \quad \x = \Qiso_{\k}(\xib).
\end{equation}
\rt{To compute the volume integrals in \eqref{weak3} numerically, we employ the following transformation formulas \cite{Bathe2006}}
\begin{equation} \label{transform}
    \dx = |\det J_{\Tk}(\xib)| \dxi, \quad
    \begin{pmatrix}
        \frac{\partial}{\partial x_1} \\
        \frac{\partial}{\partial x_2} \\
        \frac{\partial}{\partial x_3}
    \end{pmatrix} = 
    \big(J_{\Tk}(\xib)\big)^{-\top}
    \begin{pmatrix}
        \frac{\partial}{\partial \xi_1} \\
        \frac{\partial}{\partial \xi_2} \\
        \frac{\partial}{\partial \xi_3}
    \end{pmatrix}, \qquad
    \kset, 
\end{equation}
and consider a numerical quadrature \cite{Bathe2006} in $\Tref$ in the general form:
\begin{equation} \label{GQTref}
    \int_{\Tref} G(\xib) \; \dxi \approx \sum\limits_{\q=1}^{\nq} \omegaq \, G(\Aq),
\end{equation}
where $G$ denotes the integrand, $\hat{A}_{\q}$ are the quadrature points inside $\Tref$, and $\omegaq$ are the corresponding weights. Then, one can derive the following formula:
\begin{equation} \label{GQT}
    \int_{\Tk} G(\x) \; \dx \approx \sum\limits_{\q=1}^{\nq} \omegaq \, \big|\det J_{\Tk}(\hat{A}_{\q})\big| \; G(\Akq), \qquad \kset,
\end{equation}
where $\Akq$ are the quadrature points of $\Tk \in \calTh$ given by the transformation \eqref{IsoTra} and $J_{\Tk}(\xi_1,\xi_2,\xi_3)$ is the corresponding $3 \times 3$ Jacobian matrix on $\Tk$. Using \eqref{IsoTra}, it can be expressed as
\begin{equation} \label{Jacobian}
    [J_{\Tk}(\xib)]_{\j,\m} = \frac{\partial x_{\j}}{\partial \xi_{\m}} = \sum\limits_{\p=1}^{\np} \frac{\partial \phirefp(\xib)}{\partial \xi_{\m}} [\Nkl]_{\j}, \quad \j, \mset, \;\; \kset,
\end{equation}
where $\j$ and $\m$ index the physical and reference coordinates, respectively, and $[\Nkl]_{\j}$ is the $\j$-th coordinate of node $\Nkl$.

Using \eqref{IsoTra} and \eqref{transform}, we obtain the following formula for the partial derivatives of local basis functions on the whole $\calTh$:
\begin{equation} \label{philocder}
    \frac{\partial \phikl}{\partial x_{\j}}(\Akq) = \sum\limits_{\m=1}^{3}\big[\big(J_{\Tk}(\Aq)\big)^{-1}\big]_{\m,\j} \frac{\partial \phirefp}{\partial \xi_{\m}}(\Aq), \qquad \jset, \; \kset, \; \pset, \; \qset. 
\end{equation}

\noindent
Similarly, the Gauss quadrature rules can also be applied to surface integrals, which can appear, for example, when Neumann boundary conditions are prescribed. Further technical details can be found in \cite{CeSyVa2019}.

Passing from one element to the whole computational mesh, we form an ordered list of all mesh nodes as a vector $(N_1, ..., N_{\nn})$. Now the function $\vh$ may be represented by the vector
\begin{equation}
\mathbf{v}
=
\big(
v_{h}^1(N_1),\,
v_{h}^2(N_1),\,
v_{h}^3(N_1),\,
\ldots,\,
v_{h}^1(N_{\nn}),\,
v_{h}^2(N_{\nn}),\,
v_{h}^3(N_{\nn})
\big)^{\top}
\in \mathbb{R}^{3\nn},
\end{equation}
which belongs to the space
\begin{equation}
\bbV := \left\{
\mathbf{v} \in \mathbb{R}^{3\nn}
\;\middle|\;
v_d = 0,\quad d \in \calI_\textnormal{D}
\right\},
\end{equation}
\rt{where $\calI_\textnormal{D} \subset \{1,\ldots,3\nn\}$ is the index set of prescribed displacement degrees of freedom, and the constraints defining $\bbV$ impose the corresponding
homogeneous Dirichlet boundary conditions.}

\rt{To express $(\mathrm{P}_\ts)_h$ from \eqref{weak3} in algebraic form, we introduce the nodal vectors
$\mathbf{u}_{\textnormal{D},\ts} \in \mathbb{R}^{3\nn}$ and
$\mathbf{u}_\ts \in \{\mathbf{u}_{\textnormal{D},\ts}\} + \bbV$
associated with $\boldsymbol{u}_{\textnormal{D},\ts,h}$ and
$\boldsymbol{u}_{\ts,h}$, respectively. The external and internal force
vectors are specified through their action on the test vector $\mathbf{v}$:}

\begin{equation}
\mathbf{f}_\ts \in \mathbb{R}^{3\nn}, \qquad
\mathbf{v}^{\top}\mathbf{f}_\ts
:=
\int_{\Omegah} \boldsymbol{f}_{\textnormal{V},\ts,h}\cdot \vh \, \mathrm{d}x
+
\int_{\Gamma_{\textnormal{N},h}} \boldsymbol{f}_{\textnormal{T},\ts,h}\cdot \vh \, \mathrm{d}s,
\quad \forall \, \vh \in \calVh,
\end{equation}
\vspace{-1mm}
\begin{equation}
\mathbf{g}_\ts : \mathbb{R}^{3\nn} \to \mathbb{R}^{3\nn}, \qquad
\mathbf{v}^{\top}\mathbf{g}_\ts(\mathbf{u}_\ts)
:=
\int_{\Omegah}
\boldsymbol{P}_\ts \: \!(\boldsymbol{u}_{\ts,h})
:\nabla \vh\,\mathrm{d}x,
\quad \forall \, \vh \in \calVh.
\end{equation}

\noindent
\rt{Consequently, the problem $(\mathrm{P}_\ts)_h$ has the following algebraic form:}
\begin{equation}
(\bbP_\ts)\quad
\text{find } \mathbf{u}_\ts \in \{\mathbf{u}_{\textnormal{D},\ts}\} + \bbV
\;:\quad
\mathbf{v}^{\top}\big(\mathbf{g}_\ts(\mathbf{u}_\ts)-\mathbf{f}_\ts\big)=0,
\quad \forall\, \mathbf{v} \in \bbV .
\end{equation}
\rt{Eliminating the rows associated with prescribed displacement
degrees of freedom reduces $(\bbP_\ts)$ to the following system of nonlinear equations:}
\begin{equation}
(\bbP_{\textnormal{D},\ts}) \quad
\text{find } {\mathbf{u}}_\ts \in \{\mathbf{u}_{\textnormal{D},\ts}\} + \bbV
\;:\;
\widetilde{\mathbf{g}}_\ts(\mathbf{u}_\ts)=\widetilde{\mathbf{f}}_\ts
\quad \text{in } \mathbb{R}^{n},
\qquad n=\dim \bbV,
\end{equation}
\rt{where $\widetilde{\mathbf{g}}_\ts$ and $\widetilde{\mathbf{f}}_\ts$
contain the components of $\mathbf{g}_\ts$ and $\mathbf{f}_\ts$
corresponding to the unconstrained degrees of freedom.}

\rt{We solve this discrete problem using the semismooth Newton method.
For this purpose, we define the tangent matrix associated with a
generalized derivative of $\mathbf{g}_\ts$ at $\mathbf{u}_\ts$:}

\begin{equation}
\Kt : \mathbb{R}^{3\nn} \to \mathbb{R}^{3\nn \times 3\nn}, \qquad
\mathbf{v}^{\top}\Kt(\mathbf{u}_\ts) \mathbf{w}
:=
\int_{\Omegah}
\P^{\circ}_\ts(\boldsymbol{u}_{\ts,h})
\,\nabla \boldsymbol{w}_h
:\nabla \vh \,\mathrm{d}x,
\quad \forall \, \vh,\boldsymbol{w}_h \in \calVh.
\end{equation}
where \corr{$\P^\circ_\ts$ represents a generalized (Newton) derivative of $\P_\ts$ in the sense of semismooth Newton methods, given at the integration points by the consistent tangent modulus $\DF$ from \eqref{eqsDF}. Its preimage in the logarithmic space, $\De$, coincides with the elastic stiffness $\mathbb{C}$ for elastic response and is replaced by the elastoplastic tangent whenever plastic yielding is active; its closed form can be found in \cite{CeSyVa2019}.} \rt{At each Newton iteration, we solve a linear system with the reduced tangent stiffness matrix $\Ktti$, obtained by retaining only the rows and columns associated with unconstrained degrees of freedom:}
\begin{equation}\label{eq:Newton-red}
\Ktti(\mathbf{u}^\l_\ts)\,\delta \widetilde{\mathbf{u}}^\l_\ts = \widetilde{\mathbf{f}}_\ts - \widetilde{\mathbf{g}}_\ts(\mathbf{u}^\l_\ts)
\end{equation}
with $\ts$ denoting the current time step and $\l$ denoting the current iteration. \rt{Here, $\mathbf{u}^\l_\ts$ contains all displacement components, including the prescribed values, whereas $\delta\widetilde{\mathbf{u}}^\l_\ts$ contains only the increments of the unconstrained degrees of freedom.}

\subsection{Algorithmic setup} \label{subsec:framework}

\begin{figure}
\begin{mdframed}[linewidth=0.8pt,innerleftmargin=10pt,innerrightmargin=10pt,
                 innertopmargin=6pt,innerbottommargin=6pt,
    font=\normalfont\normalsize]
\noindent
\begin{stepenum}
  \item \textbf{Generating FEM assemblies (mostly adopted from \cite{CeSyVa2019})}:
    \begin{stepenum}
      \item Discretization of the computational domain $\Omegah$ into elements,
      \item Computing the values of the local basis functions and their derivatives at the Gauss (integration) points, computing element volumes, etc.
      \item Assembling the displacement-gradient matrix $\B$ using \eqref{Bd},
      \item Assembling the global elastic stiffness matrix $\Kel$.
    \end{stepenum}
  \item \textbf{Main computational procedure -- time stepping} (\texttt{For}-loop):
    \begin{stepenum}
      \item Setting the boundary conditions (including vector $\mathbf{f}_\ts$) and adjusting the initial displacements $\boldsymbol{u}_{\ts}^0$,
      \item Iterating via Newton's method for solving a nonlinear system (\texttt{While}-loop):
        \begin{stepenum}
          \item {Preprocessing phase}:    \label{phase:preproc}
            \begin{stepenum}
              \item Determining $\F$ from $\mathbf{u}$ and computing $\C$ according to \eqref{eqsFC}, \label{phase:preproc-C}
              \item Computing the eigenvalues and eigenvectors of $\C$ and then recovering $\E$ from \eqref{eqsEC}, \label{phase:preproc-E}
            \end{stepenum}
          \item Constitutive update phase (resolution of the {constitutive problem}) providing $\Sigb$ and $\De$, \label{phase:constit}
          \item {Postprocessing phase}: \label{phase:postproc}
            \begin{stepenum}
              \item Computing the transformation tensor $\PP$ from \eqref{eqsPPTL} and the product $\Sigb : \LL$,
              \item Recovering stress $\P$ and tangent modulus $\DF$ from \eqref{eqsP} and \eqref{eqsDF}, respectively, 
            \end{stepenum}
            \item Assembly of the internal force vector $\mathbf{g}_\ts$ and global tangent matrix $\Kt$ utilizing $\P$ and $\DF$,
             \item Solving \eqref{eq:Newton-red} for the displacement increment and updating the nodal displacement vector $\boldsymbol{u}_\ts$,
            \item Checking the convergence of Newton's method via \eqref{eq:converg-crit}.
        \end{stepenum}
      \item Updating the data and continuing to the next time step.
          \end{stepenum}
    \item \textbf{Data storage and visualization.}
\end{stepenum}
\end{mdframed}
\caption{Computational scheme for the solution of the logarithmic strain version of the von Mises elastoplastic boundary value problem with the finite element method.}
\label{fig:alg-setup}
\end{figure}

The complete numerical implementation can be described by the algorithmic setup summarized in Fig.~\ref{fig:alg-setup}. The computational core (denoted \textbf{2.}) consists of an outer time-stepping loop in which, at each time step, new displacements are computed together with an inner while-loop arising from Newton’s method. A single Newton iteration requires the assembly of the global tangent matrix $\Kt$, which involves the evaluation of all quantities introduced in this text. In practice, this is typically performed within an inner for-loop iterating over the quadrature points of all elements of $\calTh$.

\section{MATLAB implementation} 
\label{sec:numerics}

Our MATLAB implementation follows the scheme in Fig.~\ref{fig:alg-setup}. We build on the FEM tool introduced in \cite{CeSyVa2019}, whose capabilities have been continuously enhanced, expanded, and streamlined \cite{Svetlik2024, Karatson2025, Sysala2025}. Below, we will detail some crucial parts of our modified code. In the first stage (\textbf{1. Generating FEM assemblies}), we follow the approach in \cite{CeSyVa2019} except for a slightly modified construction of the displacement-gradient matrix, described in Section~\ref{subsec:preproc}.

Details on the second stage (\textbf{2. Main computational procedure -- time stepping}), which brings the main novelties of this work -- i.e., the logarithmic space approach, its vectorization, and details on integration into the workflow -- are presented in subsequent parts of this section. But first we specify the data storage and vectorization libraries we use.

\subsection{Preliminaries} \label{Subsec:preliminaries}

To fully exploit the advantages of vectorization techniques in MATLAB, it is often beneficial to reduce the dimensions of some objects or to produce a multidimensional object. As for the first situation, second-order tensors are in some cases represented by vectors containing their independent components. In such a case, our implementation operates with the Voigt notation to represent a symmetric second-order tensor ($\in \mathbb{R}^{3\times3}_{\rm sym}$) via the corresponding vector (or vice versa)
\begin{equation}
    \begin{pmatrix}
        C_{11} & {\color{black} C_{12}} & {\color{black} C_{13}} \\
        {\color{black} C_{21}} & C_{22} & {\color{black} C_{23}} \\
        {\color{black} C_{31}} & {\color{black} C_{32}} & C_{33}
    \end{pmatrix} \quad \longrightarrow \quad
    (C_{11}, C_{22}, C_{33}, C_{23}, C_{13}, C_{12})^{\top}.
\end{equation}
This also motivates the ordering of the vectors representing non-symmetric second-order tensors ($\in \mathbb{R}^{3\times3}$) as follows
\begin{equation}
    \begin{pmatrix}
        F_{11} & {\color{black} F_{12}} & {\color{black} F_{13}} \\
        {\color{black} F_{21}} & F_{22} & {\color{black} F_{23}} \\
        {\color{black} F_{31}} & {\color{black} F_{32}} & F_{33}
    \end{pmatrix} \quad \longrightarrow \quad (F_{11}, F_{22}, F_{33}, F_{23}, F_{13}, F_{12}, F_{32}, F_{31}, F_{21})^{\top}.
\end{equation}

Moreover, modern vectorization tools allow performing matrix operations on a set of matrices simultaneously through multidimensional arrays, i.e., ``stacks'' of matrices. Hence, for instance, the set of (symmetric) logarithmic stress tensors at all integration points is stored both as a $6 \times \nint$ ``stack'' of vectors (i.e., MATLAB 2-D array denoted \texttt{Si}) for computing the constitutive response and as a $3 \times 3 \times \nint$ ``stack'' of matrices (i.e., MATLAB 3-D array denoted \texttt{Si\_3D}) for the tensor operations in the postprocessing phase. Table~\ref{tab:matlab-notation} provides an overview of the main quantities used in the model together with their counterparts in MATLAB.

\begin{table}[htbp]
    \centering
    \renewcommand{\arraystretch}{1.15}
    \begin{NiceTabular}{lccc}[left-margin=0mm]
        \toprule
        Quantity & Symbol & \matlab{} variable & Array size \\
        \toprule

        Nodal displacements
            & $\mathbf{u}$ & \texttt{U} & $3 \times \nn$ \\
        Displacement-gradient matrix
            & $\B$ & $\ttBd$ & $(9\nint) \times (3\nn)$ \\
        Elastic stiffness matrix
            & $\Kel$ & \texttt{K\_el}
            & $(3\nn) \times (3\nn)$ \\
        \midrule[0.2pt]

        \Block{2-1}{Deformation gradient \eqref{eq:def-grad}}
            & \Block{2-1}{$\F$}
            & \texttt{F} & $9 \times \nint$ \\
            & & \texttt{F\_3D} & $3 \times 3 \times \nint$ \\
        \addlinespace

        Right Cauchy--Green tensor \eqref{eqsFC}
            & $\C$ & \texttt{C\_3D}
            & $3 \times 3 \times \nint$ \\
        Eigenvectors of $\C$
            & $\boldsymbol{N}$ & \texttt{eVec\_3D}
            & $3 \times 3 \times \nint$ \\
        Eigenvalues of $\C$ (diagonal matrices)
            & $\boldsymbol{\lambda}$ & \texttt{eVal\_3D}
            & $3 \times 3 \times \nint$ \\
        \addlinespace

        \Block{2-1}{Logarithmic (Hencky) strain \eqref{eqsEC}}
            & \Block{2-1}{$\E$}
            & \texttt{E} & $6 \times \nint$ \\
            & & \texttt{E\_3D} & $3 \times 3 \times \nint$ \\
        \midrule[0.2pt]

        \Block{2-1}{Logarithmic stress tensor
            (Fig.~\ref{fig:const-mod-scheme} and \eqref{eq:log-stress})}
            & \Block{2-1}{$\Sigb$}
            & \texttt{Si} & $6 \times \nint$ \\
            & & \texttt{Si\_3D} & $3 \times 3 \times \nint$ \\
        Consistent tangent modulus (Fig.~\ref{fig:const-mod-scheme})
            & $\De$ & \texttt{DSi} & $36 \times \nint$ \\
        Plastic strain (Fig.~\ref{fig:const-mod-scheme})
            & $\boldsymbol{E}^\textnormal{p}$ & \texttt{Ep}
            & $6 \times \nint$ \\
        \midrule[0.2pt]

        Transformation tensor \eqref{eqsPPTL}
            & $\PP$ & \texttt{PP\_3D}
            & $6 \times 9 \times \nint$ \\
        Tensor contraction in \eqref{eqsPPTL}, \eqref{eqsDF}
            & $\Sigb : \LL$ & \texttt{SiLL\_3D}
            & $9 \times 9 \times \nint$ \\
        Consistent tangent modulus \eqref{eqsDF}
            & $\DF$ & \texttt{DP\_3D}
            & $9 \times 9 \times \nint$ \\
        \addlinespace

        \Block{2-1}{First Piola--Kirchhoff stress \eqref{eqsP}}
            & \Block{2-1}{$\P$}
            & \texttt{P} & $9 \times \nint$ \\
            & & \texttt{P\_3D} & $3 \times 3 \times \nint$ \\
        \midrule[0.2pt]

        Nodal internal forces
            & $\mathbf{g}$ & \texttt{g} & $3 \times \nn$ \\
        Nodal external forces
            & $\mathbf{f}$ & \texttt{f} & $3 \times \nn$ \\
        Global tangent stiffness matrix
            & $\Kt$ & \texttt{K\_tan}
            & $(3\nn) \times (3\nn)$ \\
        \midrule[0.2pt]

        Second Piola--Kirchhoff stress \eqref{eqsSsi}
            & $\S$ & \texttt{S\_3D}
            & $3 \times 3 \times \nint$ \\
        Cauchy stress \eqref{eqsCauchy}
            & $\sigb$ & \texttt{si\_3D}
            & $3 \times 3 \times \nint$ \\

        \bottomrule

        \CodeAfter
        \begin{tikzpicture}[
            phase brace/.style={
                decorate,
                decoration={brace, amplitude=4pt},
                thick
            },
            phase label/.style={
                rotate=90,
                font=\footnotesize,
                align=center
            }
        ]
            \draw[phase brace]
                ([xshift=0mm,yshift=1mm]12-|1)
                -- ([xshift=0mm,yshift=-1mm]5-|1)
                node[midway,xshift=-5mm,phase label]
                {Preprocessing};

            \draw[phase brace]
                ([xshift=0mm,yshift=1mm]16-|1)
                -- ([xshift=0mm,yshift=-1mm]12-|1)
                node[midway,xshift=-6mm,phase label]
                {Constitutive\\problem};

            \draw[phase brace]
                ([xshift=0mm,yshift=1mm]21-|1)
                -- ([xshift=0mm,yshift=-1mm]16-|1)
                node[midway,xshift=-5mm,phase label]
                {Postprocessing};

            \draw[phase brace]
                ([xshift=0mm,yshift=1mm]24-|1)
                -- ([xshift=0mm,yshift=-1mm]21-|1)
                node[midway,xshift=-6mm,phase label]
                {Newton's\\method};
        \end{tikzpicture}
    \end{NiceTabular}
    \caption{The overview of model quantities, their symbols, corresponding MATLAB array names and sizes. Both eigenvalues and eigenvectors are stored as 3-D arrays only.}
    \label{tab:matlab-notation}
\end{table}

The assembly procedures for most quantities associated with the mesh
$\calTh$ and the finite element discretization are adopted from
\cite{CeSyVa2019}. In particular:
\begin{itemize}
    \item[] \texttt{COORD} -- a $3 \times \nn$ matrix of nodal coordinates;
    \item[] \texttt{ELEM} -- an $8 \times \nt$ matrix containing
    the global indices of the nodes of each element, representing
    the local-to-global mapping $\Iphi$ from \eqref{localPhi};
    \item[] \texttt{SURF} -- a $4 \times \nf$ matrix containing
    the global indices of the nodes of each boundary face;
    \item[] \texttt{HatP} -- a $\np \times \nq$ matrix of values of basis functions, $\phirefp$, at the quadrature points of $\Tref$, cf. \eqref{IsoTra} and \eqref{philocref};
    \item[] \texttt{DHatP1}, \texttt{DHatP2}, \texttt{DHatP3} -- $\np \times \nq$ matrices of the partial derivatives of basis functions, $\phirefp$, at the quadrature points of the reference element, $\Tref$, cf. \eqref{IsoTra} and \eqref{philocref};
    \item[] \texttt{DPhi1}, \texttt{DPhi2}, \texttt{DPhi3} -- $\np \times \nint$ matrices of the partial derivatives of local basis functions $\phikl$ at all quadrature points of $\calTh$, cf. \eqref{philocder};
    \item[] \texttt{DET} -- determinants of Jacobian $3 \times 3$ matrices, $J_{\Tk}(\xi_1,\xi_2,\xi_3)$, at all quadrature points of $\calTh$, cf. \eqref{Jacobian}.
\end{itemize}

\subsection{Vectorization libraries} \label{subsec:vectorization}

Our \matlab{} implementation benefits from the vectorization techniques from \cite{MoRaVaVa2026} for most algebraic operations between 3-D arrays of matrices, which allows us to avoid time-consuming for-loops mostly over quadrature points; such a form of data storage is denoted by suffix $\texttt{\_3D}$ in the implementation. The origins of the implementation date back to the FEM assemblies for nodal \cite{RaVa2013} and edge \cite{AnVa2015} elements. Some original library functions were optimized due to the implicit (also called arithmetic or broadcasting) expansion feature that \matlab{} came with in R2016b. Later, an additional speedup was achieved due to new page-wise functions introduced in \matlab{} releases from 2020 to 2023, which paved the way for new vectorized FEM implementations (e.g. \cite{TcKrPoKrJi2024}):
\begin{description}
    \item \texttt{pagetranspose} -- page-wise matrix transposition (since R2020b)
    \item \texttt{pagemtimes} -- page-wise matrix multiplication (since R2020b)
    \item \texttt{pagemldivide} -- page-wise left matrix divide for solving linear systems (since R2022a)
    \item \texttt{pageeig} -- page-wise spectral decomposition of square matrices (since R2023a)
\end{description}
In particular, the last one provides a significant speedup for the preprocessing part (Section~\ref{subsec:preproc}), as will be documented on the computational-time comparison plots in Section~\ref{subsec:vector-assess}. In our \matlab{} implementation, we heavily employ the following functions based on the libraries developed in \cite{RaVa2013} and \cite{MoRaVaVa2026}:
\begin{description}
    \item \texttt{amt} -- page-wise transposition of an array of matrices,
    \item \texttt{amdet} -- page-wise determinant calculation of an array of square matrices,
    \item \texttt{amtam} -- page-wise multiplication of the transposed array of matrices with another array of matrices,
    \item \texttt{amamt} -- page-wise multiplication of an array of matrices with an array of transposed matrices,
    \item \texttt{asam} -- page-wise multiplication of an array of scalars and an array of matrices, 
\end{description}
In the current work, we extended the library with the following function:
\begin{description}
    \item \texttt{amam} -- page-wise multiplication of two arrays of matrices.
\end{description}

\subsection{Boundary conditions and time stepping} \label{subsec:time-step}

After establishing the objects needed for the finite element method, the boundary conditions are imposed. As in the original code, the (possibly non-uniform) time discretization is defined by prescribing the vector of time steps.\footnote{Let us note that for our constitutive model (as well as for other rate-independent ones), rescaling the input time by any positive, monotone reparametrization leads to an identical rescaling of the model’s output.} The prescribed non-homogeneous Dirichlet and Neumann boundary conditions are discretized as well, hence, the vector of external forces, $\mathbf{f}_\ts$, for each step can be assembled. For effective speed-up, the initial guess of the displacement vector, $\textbf{u}$, in the non-prescribed degrees of freedom is based on a combination of extrapolation from previously converged solutions and purely elastic response for the first few steps.

\subsection{Preprocessing} \label{subsec:preproc}

Given a matrix of all displacements, $\ttU$, from the preceding parts of the algorithm, we sequentially compute:
\begin{itemize}
    \item[$\bullet$] a deformation gradient matrix, $\ttFd$,
    of size $9 \times \nint$ in the Voigt notation at all quadrature (integration) points in all elements, and its 
    3-D version $\ttFDD$ of size $3 \times 3 \times \nint$;
    \item[$\bullet$] a 3-D array of right Cauchy--Green tensor, $\texttt{C\_3D}$, of size $3 \times 3 \times \nint$ using \texttt{amtam} function;
    \item[$\bullet$] 3-D arrays \texttt{eVal\_3D}, \texttt{eVec\_3D} of size $3 \times 3 \times \nint$ of the corresponding eigenvalues (in the form of a diagonal matrix) and eigenvectors of $\texttt{C\_3D}$ using the \texttt{pageeig} function;
    \item[$\bullet$] a 3-D array of Hencky strain, $\texttt{E\_3D}$, of size $3 \times 3 \times \nint$ containing logarithmic strains using \eqref{eqsEC} and \texttt{amamt}, \texttt{amam} functions.
\end{itemize}

Concerning $\F$, given nodal displacements, we need to compute the deformation gradients at all integration (quadrature) points of each element. The crucial part is the construction of displacement-gradient matrix, $\B$. We employ the vector representation of tensors, and for one integration point at one particular element we get the relation
\begin{equation} \label{eqB}
    \begin{pmatrix}
        F_{11}(\Akq) \\
        F_{22}(\Akq) \\
        F_{33}(\Akq) \\
        F_{23}(\Akq) \\
        F_{13}(\Akq) \\
        F_{12}(\Akq) \\
        F_{32}(\Akq) \\
        F_{31}(\Akq) \\
        F_{21}(\Akq)
    \end{pmatrix} \;
    = \;
    \begin{pmatrix}
        D_1 \varphi_{\k,1}^\q & 0 & 0 & \hdots & D_1 \varphi_{\k,\np}^\q & 0 & 0 \\
        0 & D_2 \varphi_{\k,1}^\q & 0 & \hdots & 0 & D_2 \varphi_{\k,\np}^\q & 0 \\
        0 & 0 & D_3 \varphi_{\k,1}^\q & \hdots & 0 & 0 & D_3 \varphi_{\k,\np}^\q \\
        0 & D_3 \varphi_{\k,1}^\q & 0 & \hdots & 0 & D_3 \varphi_{\k,\np}^\q & 0 \\
        D_3 \varphi_{\k,1}^\q & 0 & 0 & \hdots & D_3 \varphi_{\k,\np}^\q & 0 & 0 \\
        D_2 \varphi_{\k,1}^\q & 0 & 0 & \hdots & D_2 \varphi_{\k,\np}^\q & 0 & 0 \\
        0 & 0 & D_2 \varphi_{\k,1}^\q & \hdots & 0 & 0 & D_2 \varphi_{\k,\np}^\q \\
        0 & 0 & D_1 \varphi_{\k,1}^\q & \hdots & 0 & 0 & D_1 \varphi_{\k,\np}^\q \\
        0 & D_1 \varphi_{\k,1}^\q & 0 & \hdots & 0 & D_1 \varphi_{\k,\np}^\q & 0
    \end{pmatrix} \;
    \begin{pmatrix}
        u_1(N_\k^1) \\
        u_2(N_\k^1) \\
        u_3(N_\k^1) \\
        \tenvdots \\
        u_1(N_\k^{\np}) \\
        u_2(N_\k^{\np}) \\
        u_3(N_\k^{\np})
    \end{pmatrix} \;
    + \;
    \begin{pmatrix}
        1 \\
        1 \\
        1 \\
        0 \\
        0 \\
        0 \\
        0 \\
        0 \\
        0
    \end{pmatrix}
\end{equation}
with $D_\j \varphi_{\k,\p}^\q$ being the value of the $\j$-th partial derivative of the $\p$-th local basis function at the $\q$-th integration point of the $\k$-th element.

For the entire computational mesh, a (sparse) global displacement-gradient matrix, $\B$, of size $(9*\nint) \times (3*\nn)$ has to be assembled. It can be expressed as a $\nint \times \nn$ block matrix with blocks $\B_{\k,\q}^{\ii}$, $\ii \in \{1, \hdots, \nn\}$, $\j \in \{1, \hdots, 3\}$, $\k \in \{1, \hdots, \nt\}$, $\q \in \{1, \hdots, \nq\}$ as follows:
\begin{equation}
\label{Bd}
\begin{array}{c@{\qquad}c@{\qquad}c}
\B =
\begin{pmatrix}
    \B_{1,1}^1 & \B_{1,1}^2 & \dots & \B_{1,1}^{\nn} \\
    \B_{1,2}^1 & \B_{1,2}^2 & \dots & \B_{1,2}^{\nn} \\
    \vdots & \vdots & \ddots & \vdots \\
    \B_{1,\nq}^1 & \B_{1,\nq}^2 & \dots & \B_{1,\nq}^{\nn} \\
    \B_{2,1}^1 & \B_{2,1}^2 & \dots & \B_{2,1}^{\nn} \\
    \vdots & \vdots & \ddots & \vdots \\
    \vdots & \vdots & \ddots & \vdots \\
    \vdots & \vdots & \ddots & \vdots \\
    \B_{\nt,1}^1 & \B_{\nt,1}^2 & \dots & \B_{\nt,1}^{\nn} \\
    \vdots & \vdots & \ddots & \vdots \\
    \B_{\nt,\nq}^1 & \B_{\nt,\nq}^2 & \dots & \B_{\nt,\nq}^{\nn}
\end{pmatrix},
\qquad \qquad \qquad
\B_{\k,\q}^{\ii} =
\begin{pmatrix}
    D_1 \Phi_{\k,\q}^{\ii} & 0 & 0 \\
    0 & D_2 \Phi_{\k,\q}^{\ii} & 0 \\
    0 & 0 & D_3 \Phi_{\k,\q}^{\ii} \\
    0 & D_3 \Phi_{\k,\q}^{\ii} & 0 \\
    D_3 \Phi_{\k,\q}^{\ii} & 0 & 0 \\
    D_2 \Phi_{\k,\q}^{\ii} & 0 & 0 \\
    0 & 0 & D_2 \Phi_{\k,\q}^{\ii} \\
    0 & 0 & D_1 \Phi_{\k,\q}^{\ii} \\
    0 & D_1 \Phi_{\k,\q}^{\ii} & 0
\end{pmatrix}
\end{array}
\end{equation}
with $D_{\j} \Phi_{\k,\q}^{\ii} = \frac{\partial \Phii}{\partial x_{\j}}(\Akq)$ being the value of the $\j$-th partial derivative of the $\ii$-th global basis function at the $\q$-th quadrature node of the $\k$-th element. If the $\ii$-th node is not adjacent to the $\k$-th element, then $\Phi_{\k,\q}^{\ii}$ and its derivatives vanish in all quadrature points of $\Tk$, and therefore $\B_{\k,\q}^{\ii}$ is a zero matrix. Code block~\ref{code1} shows the algorithmization of the preprocessing part (denoted \textbf{(A)} in Fig.~\ref{fig:alg-setup}).

\begin{mm}[caption={Vectorized assemblies in the preprocessing part.},label={code1}]
% deformation gradients (matrix)
F(:) = B*U_it(:);  % 9 x nint
F = F + [1 1 1 0 0 0 0 0 0]';
% deformation gradients (3-D array)
F_3D = reshape(F(inds_Voigt,:),3,3,n_int);  % 3 x 3 x nint
% right Cauchy-Green tensors (3-D array)
C_3D = amtam(F_3D,F_3D);  % 3 x 3 x nint
% spectral decompositions of 'C' (3-D arrays)
[eVec_3D,eVal_3D] = pageeig(C_3D);  % (for one slice, C = eVec*eVal*eVec')
eVal_ln_3D = eVal_3D;  % 3 x 3 x nint
eVal_ln_3D(inds_diag) = log(eVal_3D(inds_diag));  % applying logarithm elementwise to the diagonals
% logarithmic (Hencky) strains (3-D array)
E_3D = amamt(amam(eVec_3D,eVal_ln_3D),eVec_3D)/2;  % 3 x 3 x nint
% logarithmic (Hencky) strains (matrix)
E = reshape(E_3D,9,[]);
E = E([1 5 9 8 7 4],:).*[1;1;1;2;2;2];  % 6 x nint
\end{mm}

\subsection{Constitutive update}
\label{subsec:const-mod}

The main benefit of the logarithmic strain approach is its modular nature, which allows for a complete adoption of the constitutive update algorithm used in the small-strain formulation. Hence, the elastoplastic update algorithm from \cite{CeSyVa2019} (corresponding to the non-regularized constitutive model from Section~\ref{sec:maths}) can be directly incorporated. Just the internal variable arrays are now interpreted as the logarithmic plastic strain, $\E^\textnormal{p}$, and as an auxiliary variable related to the logarithmic back stress, $\boldsymbol{X}$, imposed by kinematic hardening; cf. \eqref{eq:back-stress}. Let us note that isotropic elasticity is considered; hence, the elasticity tensor is recovered from only two independent elastic parameters.

The employed constitutive model allows a closed-form, fully vectorized implementation; hence, the constitutive update is resolved for all integration (material) points simultaneously. Note, however, that this might not work for some other constitutive models, where the iterative process is employed, e.g., our recent work \cite{MHVKJSF2025}. Code block~\ref{code6} shows the constitutive function call with updated names of input and output variables.

\vspace{2mm}
\begin{mm}[caption={Vectorized constitutive update function call.},label={code6}]
[Si,DSi,IND_p,Ep,Hard] = constitutive_problem(E,Ep_old,Hard_old,shear,bulk,a,Y);
\end{mm}

\subsection{Postprocessing}

\begin{figure}
\begin{mdframed}[linewidth=0.9pt,innerleftmargin=10pt,innerrightmargin=10pt,
                 innertopmargin=8pt,innerbottommargin=8pt]
\setstretch{0.20}
\begin{enumerate}
\item \textbf{Eigenvalues, eigenvectors, diagonal functions and projections}

Get eigenvalues $\lambda_i$ and eigenvectors $\boldsymbol{N}^i$ from eigenvalue problem 
$\C \boldsymbol{N}^i = \lambda_i \boldsymbol{N}^i$.

Determine $\boldsymbol{n}^i = \F \boldsymbol{N}^i$, and let $\boldsymbol{N}_i = \boldsymbol{N}^i$, $\boldsymbol{n}_i = \boldsymbol{n}^i$ (since constant standard metrics used). For $i,j=1,2,3$ compute
diagonal functions, derivatives and projections
\[
e_i=\tfrac12\ln\lambda_i,\quad d_i=\lambda_i^{-1},\quad f_i=-2\lambda_i^{-2}
\quad\text{and}\quad
\zeta_{ij}= \Sigb:(\boldsymbol{N}_i\otimes \boldsymbol{N}_j).
\] \\

\item \textbf{Coefficients $(i\neq j)$}

For three different eigenvalues $\lambda_{\bar{i}}\neq\lambda_{\bar{j}}\neq\lambda_{\bar{k}}$:
\[
\vartheta_{ij}:=\frac{e_i-e_j}{\lambda_i-\lambda_j},\qquad
\xi_{ij}:=\frac{\bigl(\vartheta_{ij}-\tfrac12 d_j\bigr)}{\lambda_i-\lambda_j},
\]
\[
\eta:=\sum_{i=1}^{3}\ \sum_{j\neq i}^{3}\ \sum_{\substack{k\neq i\\ k\neq j}}^{3}
\frac{e_i}{2(\lambda_i-\lambda_j)(\lambda_i-\lambda_k)}.
\]
For two equal eigenvalues $\lambda_{\bar{i}}=\lambda_{\bar{j}}\neq\lambda_{\bar{k}}$:
\[
\vartheta_{\bar{i}\bar{j}}=\vartheta_{\bar{j}\bar{i}}:=\tfrac12 d_{\bar{i}},\qquad
\xi_{\bar{i}\bar{j}}=\xi_{\bar{j}\bar{i}}:=\tfrac18 f_{\bar{i}},\qquad
\eta:=\xi_{\bar{k}\bar{i}},\qquad
\]
and, \corr{for all $i\neq j$ with $\bar{k}\in\{i,j\}$},
\[
\vartheta_{ij}:=\frac{e_i-e_j}{\lambda_i-\lambda_j}, \qquad
\xi_{ij}:=\frac{\bigl(\vartheta_{ij}-\tfrac12 d_j\bigr)}{\lambda_i-\lambda_j}.
\]
For three equal eigenvalues $\lambda_{\bar{i}}=\lambda_{\bar{j}}=\lambda_{\bar{k}}$:
\[
\vartheta_{ij}:=\tfrac12 d_{\bar{i}},\qquad
\xi_{ij}:=\tfrac18 f_{\bar{i}},\qquad
\eta:=\tfrac18 f_{\bar{i}}.
\] \\

\item \textbf{Nominal tensor operator}
\[
M^{ijaB}=n_{ia}N^{jB}+n_{ja}N^{iB},
\] \\

\item \textbf{Nominal projection tensors}
\[
\PP=\sum_{i=1}^{3}\tfrac12 d_i\,\boldsymbol{N}_i\otimes \boldsymbol{N}_i\otimes \boldsymbol{M}^{ii}
+\sum_{i=1}^{3}\ \sum_{j\neq i}^{3}\vartheta_{ij}\,\boldsymbol{N}_i\otimes \boldsymbol{N}_j\otimes \boldsymbol{M}^{ij}.
\]
\[
\begin{aligned}
\Sigb:\LL
&=\sum_{i=1}^{3}\tfrac14 f_i\,\zeta_{ii}\,\boldsymbol{M}^{ii}\otimes \boldsymbol{M}^{ii}
+\sum_{i=1}^{3}\ \sum_{j\neq i}^{3}\ \sum_{\substack{k\neq i\\ k\neq j}}^{3}
2\eta\,\zeta_{ij}\,\boldsymbol{M}^{ik}\otimes \boldsymbol{M}^{jk} \\
&\quad
+\sum_{i=1}^{3}\ \sum_{j\neq i}^{3}
2\xi_{ij}\Bigl[\zeta_{ij}\bigl(\boldsymbol{M}^{ij}\otimes \boldsymbol{M}^{jj}+\boldsymbol{M}^{jj}\otimes \boldsymbol{M}^{ij}\bigr)
+\zeta_{jj}\,\boldsymbol{M}^{ij}\otimes \boldsymbol{M}^{ij}\Bigr] \\
&\quad
+ (\boldsymbol{I}\F)^{-1}(\Sigb:\PP)\odot \boldsymbol{I}.
\end{aligned}
\]
\end{enumerate}
\end{mdframed}
\caption{Spectral decomposition and construction of projection tensors for the preprocessing and postprocessing phases used in our code. The constant standard (Cartesian) metrics are considered. 
Adapted from \cite{MiApLa2002}.}
\label{fig:Miehe-alg}
\end{figure}

The postprocessing part of the algorithmic workflow in Fig.~\ref{fig:alg-setup} uses the transformation tensors $\PP$ and $\LL$ from \eqref{eqsPPTL}. $\PP$ is used to compute both the first Piola--Kirchhoff stress, $\P$, from \eqref{eqsP} and the consistent tangent modulus, $\DF$, from \eqref{eqsDF}, while $\LL$ only appears in the assembly of $\DF$ in the tensor contraction $\Sigb : \LL$. Since some auxiliary variables appear in both assemblies of $\PP$ and $\Sigb : \LL$ (described in detail in Fig.~\ref{fig:Miehe-alg} and adopted from \cite{MiApLa2002}), we assemble them simultaneously in the \texttt{assembly\_PP\_SiLL} function, see also Code block~\ref{code77} in Appendix~\ref{app:PP-TL-code}. These assemblies are rather technical and include evaluation of (up to three-level) nested for-loops, where the indices of the inner loops depend on the indices of the outer ones. Therefore, vectorization of this part requires special treatment, where logical vectors \texttt{inds1}, \dots, \texttt{inds5} are used to distinguish between different cases of eigenvalue (in)equality (see part 2. in Fig.~\ref{fig:Miehe-alg}). These are consequently used in the final assemblies of $\PP$ and $\Sigb : \LL$, stored as three-dimensional arrays (for all integration points of $\calTh$), where the \texttt{amam}, \texttt{amtam}, \texttt{asam}, and \texttt{amt} functions listed in Subsection~\ref{subsec:vectorization} are utilized.

Since we are employing the first Piola stress in the virtual work principle \eqref{weak0FS} and working in the Cartesian metric, the left term of the Kronecker product
\begin{equation} \label{eq:KrDelta}
    (\boldsymbol{I}\F)^{-1}(\Sigb:\PP)\odot \boldsymbol{I},
\end{equation}
from the assembly of $\Sigb : \LL$ in the last relation in Fig.~\ref{fig:Miehe-alg} condenses merely to the second Piola--Kirchhoff stress, $\S$, cf. \eqref{eqsP} and \eqref{eqsSsi}. Taking the three-dimensional arrays of deformation gradients, $\F$, and the first Piola--Kirchhoff stress, $\P$, at all integration points and utilizing the \matlab{} \texttt{pagemldivide} function representing the page-wise left matrix divide, we can compute $\S$ at all integration points, which is subsequently used to compute the Kronecker delta products from \eqref{eq:KrDelta}. For better visual clarity, these products are not incorporated in the \texttt{assembly\_PP\_SiLL} function and are treated separately. Using Voigt notation and the symmetry of $\S$, the Kronecker delta product \eqref{eq:KrDelta} for any integration point of $\calTh$ is given by the mapping

\begingroup
\renewcommand{\arraystretch}{1.5} 
\setlength{\arraycolsep}{3pt}    
\begin{equation} \label{StoKrPr}
    \begin{pmatrix}
        \S_{11} & \S_{12} & \S_{13} \\
        \S_{12} & \S_{22} & \S_{23} \\
        \S_{13} & \S_{23} & \S_{33}
    \end{pmatrix}
    \quad \longrightarrow \quad
    \left(
    \begin{array}{ccc|ccc|ccc}
        \S_{11} & 0 & 0 & 0 & \S_{13} & \S_{12} & 0 & 0 & 0 \\
        0 & \S_{22} & 0 & \S_{23} & 0 & 0 & 0 & 0 & \S_{12} \\
        0 & 0 & \S_{33} & 0 & 0 & 0 & \S_{23} & \S_{13} & 0 \\
        \hline
        0 & \S_{23} & 0 & \S_{33} & 0 & 0 & 0 & 0 & \S_{13} \\
        \S_{13} & 0 & 0 & 0 & \S_{33} & \S_{23} & 0 & 0 & 0 \\
        \S_{12} & 0 & 0 & 0 & \S_{23} & \S_{22} & 0 & 0 & 0 \\
        \hline
        0 & 0 & \S_{23} & 0 & 0 & 0 & \S_{22} & \S_{12} & 0 \\
        0 & 0 & \S_{13} & 0 & 0 & 0 & \S_{12} & \S_{11} & 0 \\
        0 & \S_{12} & 0 & \S_{13} & 0 & 0 & 0 & 0 & \S_{11}
    \end{array}
    \right)
\end{equation}
\endgroup

Summarizing the crucial steps, using the output from the constitutive core -- current stress $\Sigb$, tangent modulus $\De$, and from the preprocessing phase -- the eigenvalues and eigenvectors of $\C$, we proceed with computing the following objects (see Fig.~\ref{fig:Miehe-alg}):
\begin{itemize}
    \item[$\bullet$] Transformation tensor $\PP$ (4th order) of size $6 \times 9 \times \nint$ and the product term $\Sigb : \LL$ (4th order) of size $9 \times 9 \times \nint$ using \eqref{eqsPPTL};
    \item[$\bullet$] The first Piola--Kirchhoff stress, $\P$, in the matrix and 3-D array representations, \texttt{P} and \texttt{P\_3D}, of size $9 \times \nint$ and $3 \times 3 \times \nint$, respectively, using \eqref{eqsP} and \texttt{amtam} function;
    \item[$\bullet$] Algebraic objects related to the transformed tangent modulus $\DF$ from \eqref{eqsDF}:
    \begin{itemize}
        \item[-] The first term, $\PP^{\top} : \De : \PP$, of size $9 \times 9 \times \nint$ using \texttt{amam}, \texttt{amtam} functions;
        \item[-] The second Piola--Kirchhoff stress, $\S$, of size $3 \times 3 \times \nint$ using \eqref{eqsSsi} and \texttt{pagemldivide} function;
        \item[-] The part given as the delta products of the second Piola--Kirchhoff tensor  of size $9 \times 9 \times \nint$, see \eqref{eq:KrDelta};
        \item[-] The product $\Sigb : \LL$ of size $9 \times 9 \times \nint$;
    \end{itemize}
    \item[$\bullet$] Additional objects required for computing new displacements:
    \begin{itemize}
        \item[-] Matrix $\mathbf{g}$ of internal forces at each node stored as $\mathbf{g} = B^{\top} P$ of size $3 \times n_n$;
    \end{itemize}
    \item[$\bullet$] Additional objects required for data analysis and visualization:
    \begin{itemize}
        \item[-] Cauchy stress, $\sigb$, of size $3 \times 3 \times \nint$ using \eqref{eqsCauchy} and \texttt{asam}, \texttt{amdet}, \texttt{amamt} functions.
    \end{itemize}
\end{itemize}

Code block~\ref{code7} provides details on the algorithmization.

\begin{mm}[caption={Vectorized assemblies in the postprocessing part.},label={code7}]
% logarithmic stress tensors (3-D array)
Si_3D = reshape([Si(1,:); Si(6,:); Si(5,:); Si(6,:); Si(2,:); Si(4,:); ...
                 Si(5,:); Si(4,:); Si(3,:)],3,3,n_int);  % 3 x 3 x nint
% assembly of 'PP' and part of 'Si:LL' (3-D array)
%         6 x 9 x nint      9 x 9 x nint
[PP_3D,SiLLpart_3D] = assembly_PP_SiLL(F_3D,Si_3D,eVal_3D(inds_diag),eVec_3D);
% first Piola-Kirchhoff stress tensors (matrix)
P = amtam(PP_3D,reshape(Si,6,1,n_int));  % 9 x 1 x nint
P = reshape(P,9,n_int);  % 9 x nint (squeeze)
% first Piola-Kirchhoff stress tensors (3-D array)
P_3D = reshape(P(inds_Voigt,:),3,3,n_int);  % 3 x 3 x nint
% second Piola-Kirchhoff stress tensors (3-D array)
S_3D = pagemldivide(F_3D,P_3D);  % 3 x 3 x nint
% complete contractions 'Si:LL' including the Kronecker delta products
SiLL_3D = SiLLpart_3D + Kron_delta_product(S_3D);  % 9 x 9 x nint
% consistent tangent moduli (3-D array)
DP_3D = amam(amtam(PP_3D,reshape(DSi,6,6,n_int)),PP_3D) + SiLL_3D;  % 9 x 9 x nint
% weighted tangent moduli (matrix)
DP_w = reshape(DP_3D,81,n_int).*repmat(WEIGHT,81,1);  % 81 x nint
% nodal internal forces (matrix)
g = B'*reshape(repmat(WEIGHT,9,1).*P, 9*n_int,1);
g = reshape(g,3,n_n);  % 3 x nn
\end{mm}

\subsection{Tangent matrix construction and Newton's method}

The final solution of the mechanical balance embodied in the system \eqref{eq:Newton-red} is found by semismooth Newton's method, see Algorithm 1 in \cite{CeSyVa2019}. The global vector of internal forces, denoted \texttt{g} in the algorithm, is obtained from the first Piola--Kirchhoff stress, \texttt{P}, and the sparse matrix representing the displacement-gradient operator at all integration points, \ttBd. The global tangent stiffness matrix, $\Kt$, crucial for proper convergence of Newton's method, is obtained as the product of the block-diagonal sparse matrix for the elastoplastic problem, \texttt{DP\_ws}, and the displacement-gradient matrix
\begin{equation}
 \texttt{K\_tan} = \ttBd^\top*\texttt{DP\_ws}*\ttBd, 
\end{equation}
hence, it is symmetric. Matrix \texttt{DP\_ws} is directly assembled from the constitutive tangents, $\DF$, at all integration points obtained in the postprocessing phase. In the present implementation, the global tangent stiffness matrix is assembled directly, without the decomposition into a precomputed elastic stiffness matrix and a ``plastic correction'' introduced in \cite{CeSyVa2019}. In the logarithmic strain formulation, the transformed tangent modulus depends on the current deformation even for a purely elastic response, so the correction generally remains nonzero at elastic integration points and the sparsity advantage of this decomposition is lost.

A further consideration is that the modular structure of the implementation allows the constitutive core to be easily substituted by another constitutive model, e.g., the shape memory alloy model presented in our previous work \cite{MHVKJSF2025}, where activation of inelastic processes occurs frequently, i.e., the purely elastic region is rather small. Any potential benefit of separating the elastic contribution would therefore be further limited by the small extent of this region, while requiring modifications of the constitutive core and not recovering the sparsity advantage discussed above.

In order to monitor the convergence of Newton's method at the $\ts$-th
time step and the $\l$-th iteration, we employ the relative increment
$\errN$ measured in the energy norm:
\begin{equation} \label{errN}
    \errN =
    \begin{cases}
        \dfrac{\| u_{\ts}^\l - u_{\ts}^{\l-1} \|_\textnormal{e}}
        {\| u_{\ts}^\l \|_\textnormal{e}
        + \| u_{\ts}^{\l-1} \|_\textnormal{e}},
        & \text{if }
        \| u_{\ts}^\l - u_{\ts}^{\l-1} \|_\textnormal{e}
        > \epsilon_{\rm min}, \\[3mm]
        0, & \text{otherwise}.
    \end{cases}
\end{equation}
Here, $\epsilon_{\rm min}>0$ is a prescribed threshold used to avoid
unnecessary evaluation of the relative increment when the absolute
increment is already sufficiently small.
The energy norm is defined by
\begin{equation} \label{energynorm}
    \| u \|_\textnormal{e}^2
    =
    \int_{\Omega}
    \operatorname{sym}(\nabla u)
    : \mathbb{C} :
    \operatorname{sym}(\nabla u)
    \,\dx,
\end{equation}
where $\operatorname{sym}(\nabla u)
= \tfrac12(\nabla u+\nabla u^\top)$ and $\mathbb{C}$ denotes the
elasticity tensor. This definition follows the approach in
\cite{CeSyVa2019}.
In the finite element implementation, the energy norm is evaluated
using the elastic stiffness matrix $\texttt{K\_el}$:
\begin{equation}
    \| u \|_\textnormal{e}^2
    \approx
    \texttt{U(:)}^\top*
    \texttt{K\_el}*
    \texttt{U(:)},
\end{equation}
where $\texttt{U(:)}$ is the vector of nodal displacement degrees of freedom.

Thanks to the considered independence of $\mathbb{C}$ from plastic strain (or other state variables), $\texttt{K\_el}$ can be precomputed before the main computation procedure starts (cf. the scheme in Fig.~\ref{fig:alg-setup}) based on elastic parameters of the material and the displacement-gradient operator, $\B$. Newton's method then iterates until the stopping criterion
\begin{equation} \label{eq:converg-crit}
\errN < \delta_\textnormal{tol} 
\end{equation}
is satisfied.

\section{Validation and vectorization performance assessment} 
\label{sec:simulations}

In this section, we illustrate improvements in computational efficiency and capabilities to handle geometric nonlinearities when the logarithmic strain approach is used. All computational times together with visualization outputs in this paper were obtained on 
ThinkPad T14 Gen 4 with AMD Ryzen 7 PRO processor and 32 GB RAM running \matlab{} R2024a. The package with codes can be downloaded at \cite{code-Zenodo} or at \\

\vspace{-2mm}
\url{https://www.mathworks.com/matlabcentral/fileexchange/185078} \\

Table~\ref{tab:mat-params} provides material parameters of the model with their values used for input in computations. Let us note that two common elastic parameters, the bulk modulus, $K = E/(3 - 6\nu)$, and the shear modulus, $G = E/(2 + 2\nu)$, used as the input to the constitutive model (see Section~\ref{subsec:const-mod}) are straightforwardly recovered thanks to presumed elastic isotropy. Although the specific values of the parameters were motivated by the behavior of metals and alloys, they do not represent a particular material. We apply a threshold $\epsilon_{\rm min} \sim 10^{-8}$, under which both displacements $u_{\ts}^{\l-1}$ and $u_{\ts}^{\l}$ are considered equal (see \eqref{errN}), and $\delta_\textnormal{tol} \sim 10^{-11}$ as the convergence criterion.

\begin{table}
\centering
\begin{tabular}{c c c c l}
\hline
Symbol & MATLAB variable & Value & Unit & Description \\
\hline
$E$   & \texttt{young} & 210    & GPa & Young's modulus \\
$\nu$ & \texttt{poisson} & 0.3125 & --  & Poisson's ratio \\
$Y$   & \texttt{sqrt(2/3)*Y} & 1      & GPa & Yield stress in uniaxial tension \\
$H^\textnormal{kin}$ & \texttt{a} & 1; 10 & GPa & Linear kinematic hardening coefficient \\
\hline
\end{tabular}
\caption{Material parameters used in numerical simulations. \corr{The yield stress enters the code through the expression \texttt{sqrt(2/3)*Y}.} See text for details.}
\label{tab:mat-params}
\end{table}

\subsection{Vectorization assessment} \label{subsec:vector-assess}

To demonstrate the computational benefits of vectorization, we consider a simple problem:
uniaxial stretching of the unit cube $\Omega=(0,1)^3$ to double its original length in the $x$-direction,  corresponding to a maximum prescribed displacement $u_{\rm max}=1$, i.e., 100\% nominal strain. We impose homogeneous Dirichlet boundary conditions $u_x=0$ on the wall $x=0$, $u_y=0$ on the wall $y=0$, and $u_z=0$ on the wall $z=0$. On the opposite wall $x=1$, the displacement component $u_x$ is prescribed uniformly and increased from zero to $u_{\rm max}$. All remaining displacement components on the walls are unconstrained.

\begin{figure}
\centering
\includegraphics[width=\textwidth]{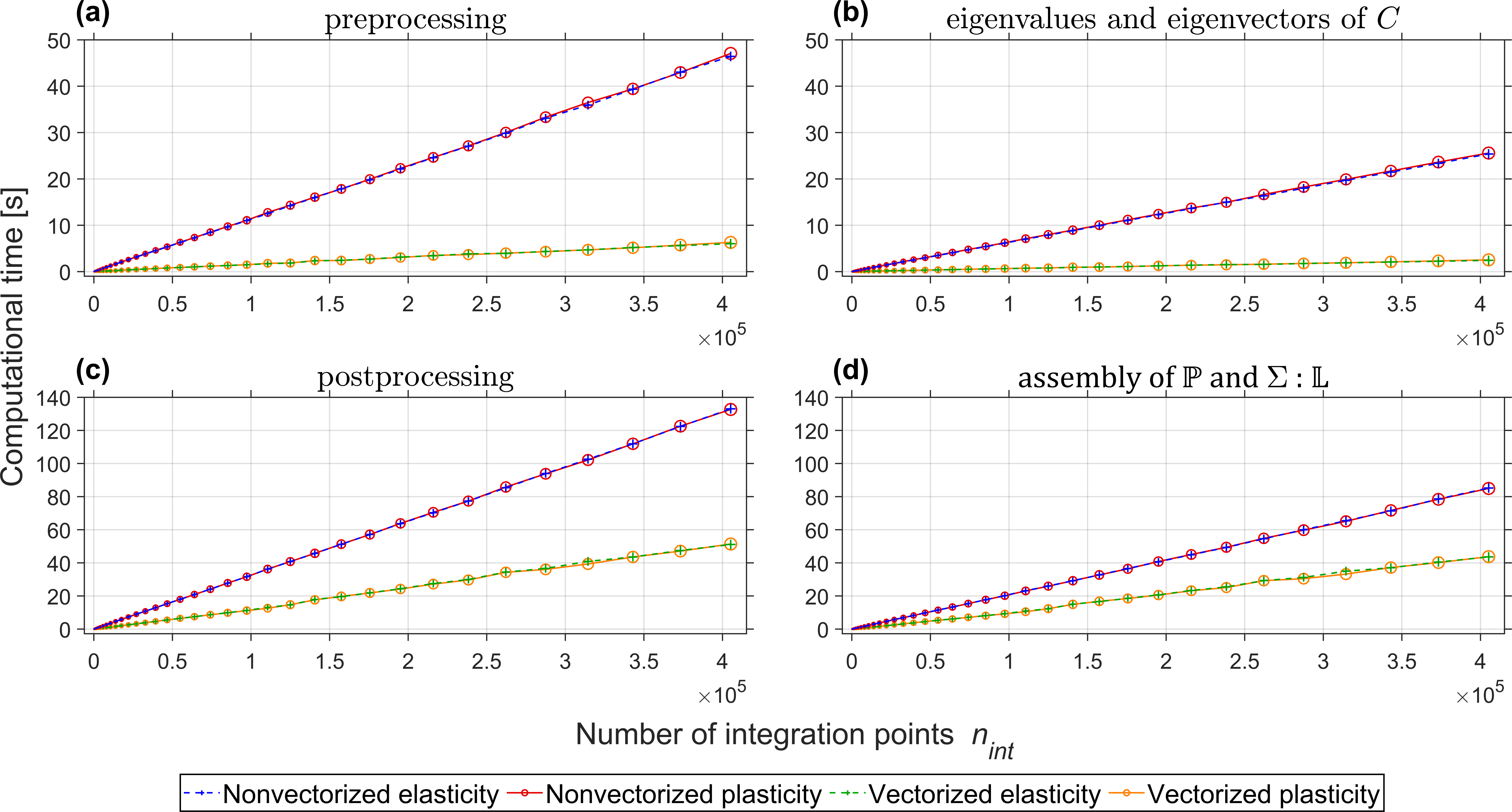}
\vspace{-3mm}
\caption{Computational times for preprocessing and postprocessing phases
in the logarithmic strain formulation for uniaxial stretching of a unit cube. Vectorized and nonvectorized implementations
are compared on \nmesh{} uniformly refined meshes for purely elastic
and fully plastic responses of an isotropic material. The panels show:
(a) preprocessing, including the evaluation of $\F$, $\C$, and $\E$;
(b) the spectral decomposition of $\C$, included in (a);
(c) postprocessing, including the evaluation of $\PP$,
$\Sigb:\LL$, the first Piola--Kirchhoff stress $\P$,
and the consistent tangent modulus $\DF$;
and (d) the evaluation of $\PP$ and $\Sigb:\LL$,
included in (c). Times are accumulated over all Newton iterations. See text for details.}
\label{fig:timesACDE}
\end{figure}

Depending on the choice of the yield stress in uniaxial tension, $Y$, we can switch between purely elastic, i.e., $Y$ set to a very high number, and fully plastic, i.e., $Y$ set close to zero, variants in this particular testing problem. In the purely elastic case, no plastic correction of the stress or tangent modulus is required throughout the simulation, whereas in the fully plastic case, both quantities are corrected at every integration point in every time step. Thanks to the character of the problem, the resulting strain and stress fields are homogeneous within the body. In addition, we compare the nonvectorized and vectorized \matlab{} implementations. Altogether, we consider four variants of the model: nonvectorized elasticity, nonvectorized plasticity, vectorized elasticity, and vectorized plasticity. For each variant, we perform \matlab{} computations on \nmesh{} uniformly refined meshes of the unit cube ($s^3$ elements with $s=1, \dots, \nmesh{}$; each element with 8 integration points). Each simulation is divided into 10 equidistant time steps. For a particular type of model, the total number of Newton iterations is independent of mesh refinement and reaches 39 for elasticity and 41 for plasticity. Newton's method converges in 3--4 iterations per time step, except for the first time step, which requires 5 iterations.\footnote{To compare the assembly costs, the results were normalized in Fig.~\ref{fig:timesACDE}.}

\begin{figure}
\centering
\includegraphics[width=\textwidth]{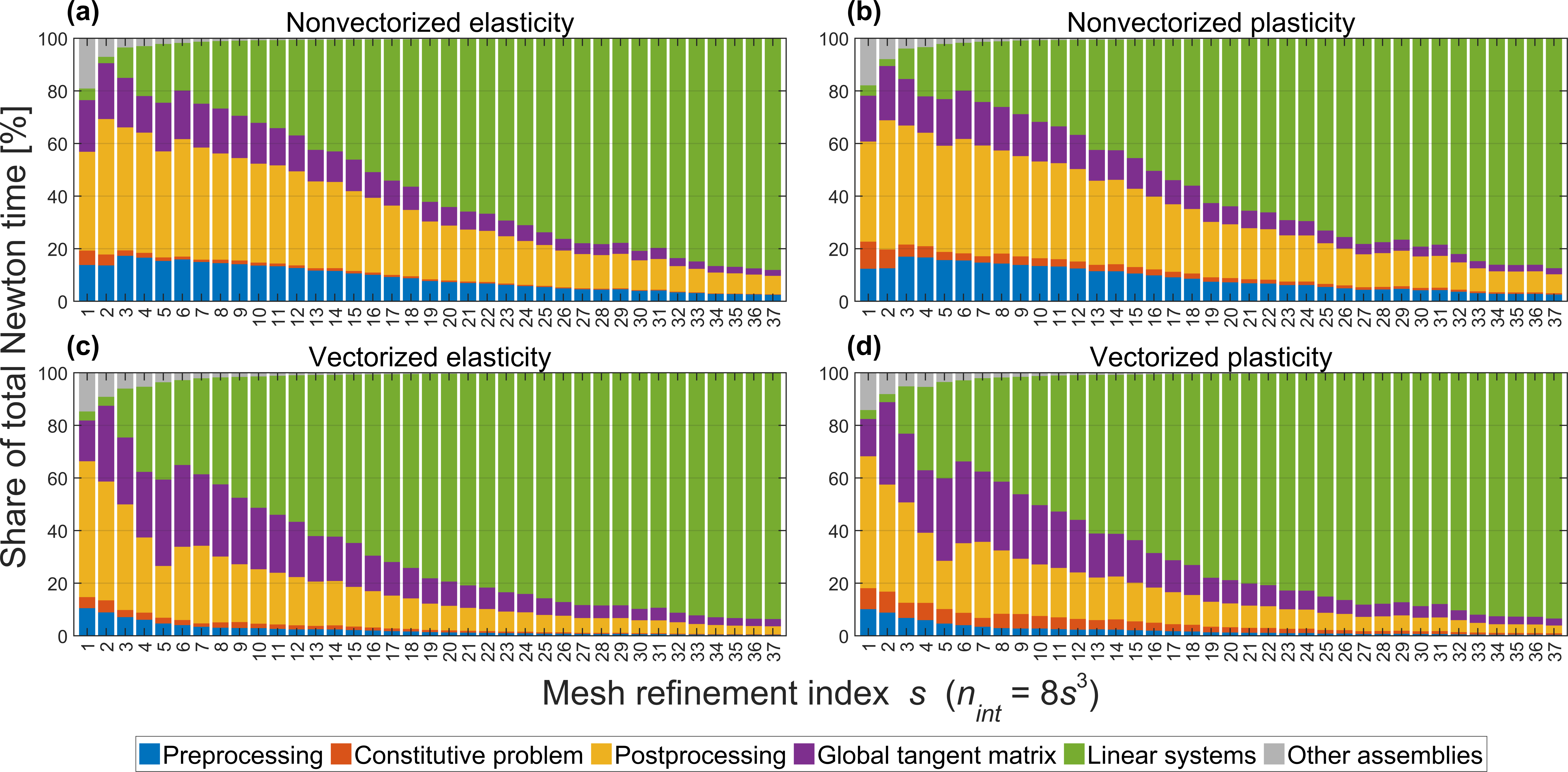}
\vspace{-3mm}
\caption{Percentage breakdown of the total time spent on Newton iterations for uniaxial stretching of a unit cube (see Subsec.~\ref{subsec:vector-assess} and Fig.~\ref{fig:alg-setup}). Results are shown for \nmesh{} uniformly refined meshes, indexed by $s$, with $n_{\mathrm{int}}=8s^3$ integration points, considering purely elastic and fully plastic responses of an
isotropic material. The panels correspond to:
(a) nonvectorized elasticity, (b) nonvectorized plasticity,
(c) vectorized elasticity, and (d) vectorized plasticity.
Each stacked bar shows the relative contributions of preprocessing, constitutive problem, postprocessing, global tangent matrix assembly, linear system solution, and the remaining operations within Newton's method. Percentages are calculated relative to the accumulated time of all Newton iterations for the corresponding mesh and case.}
\label{fig:timesperc}
\end{figure}

\begin{figure}
\centering
\includegraphics[width=\textwidth]{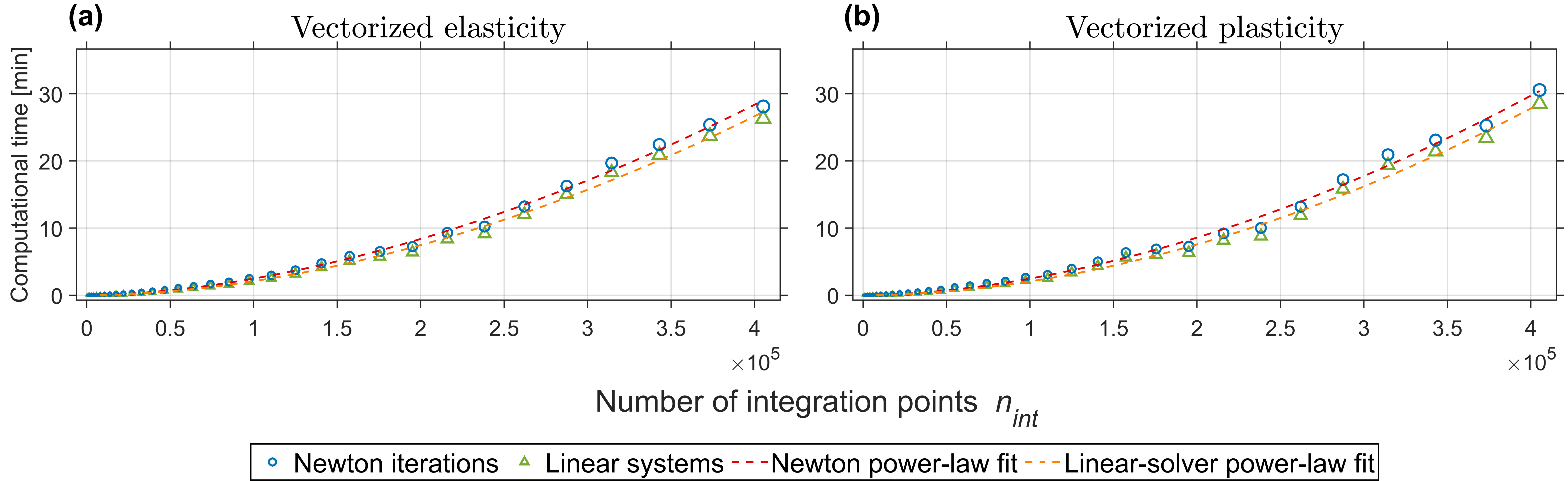}
\vspace{-3mm}
\caption{The accumulated computational time of all Newton iterations, step \textbf{(ii)} of stage \textbf{2.} in Fig.~\ref{fig:alg-setup}, and the total computational time of linear solvers, \textbf{(E)} in Fig.~\ref{fig:alg-setup}, for uniaxial stretching of a unit cube together with their power-law fits. Results are shown for vectorized code and \nmesh{} uniformly refined meshes with different numbers of integration points, considering purely elastic and fully plastic responses of an isotropic material.}
\label{fig:timesGN}
\end{figure}

Fig.~\ref{fig:timesACDE} contains plots comparing computational time with respect to the number of integration points for all four combinations in different stages of computation. \rt{The measured assembly times grow approximately linearly with the number of integration points, indicating nearly optimal scaling.} Moreover, as the data points for elasticity and plasticity are very close to each other in all cases, the computational times associated with the assembly operations in the preprocessing and postprocessing parts are largely insensitive to the constitutive response function in this benchmark.

Fig.~\ref{fig:timesACDE}(a) compares the total computational time of the preprocessing phase, denoted \textbf{(A)} in Fig.~\ref{fig:alg-setup}, where the logarithmic strain tensor $\E$ is computed at all integration points from the equation \eqref{eqsEC}. On the finest mesh, vectorization (see Code block~\ref{code1}) accelerates preprocessing by approximately a factor of 7.5. In the nonvectorized implementation, $\E$ is calculated by successively applying the \matlab{} function \texttt{eig}, which computes the eigenvalues and eigenvectors of a matrix, to the corresponding right Cauchy--Green tensors $\C$. The vectorized version operates instead with two 3-D arrays containing the eigenvalues and eigenvectors of $\C$ at all integration points via the \matlab{} function \texttt{pageeig}. Hence, all tensors $\E$ are computed using spectral decompositions of $\C$ by directly applying the logarithm to the eigenvalues and reconstructing $\E$ using the corresponding eigenvectors (common for diagonalizable matrices). This brings a computational speedup by a factor of 10. Fig.~\ref{fig:timesACDE}(b) compares computational time savings for the spectral decompositions alone in the preprocessing phase. Other time savings are due to applying functions from own-developed libraries presented in Section~\ref{subsec:vectorization} and avoiding related object conversions.

The computational time spent on the postprocessing phase, denoted \textbf{(C)} in Fig.~\ref{fig:alg-setup}, is shown in Fig.~\ref{fig:timesACDE}(c). Compared to Fig.~\ref{fig:timesACDE}(a), postprocessing is generally more time-consuming; the relative time reduction due to vectorization is lower, but the absolute reduction is higher. Fig.~\ref{fig:timesACDE}(d) provides computational times related to assemblies of the transformation tensor $\PP$ from \eqref{eqsPPTL} and the product $\Sigb : \LL$ from \eqref{eqsDF} in the first step of the postprocessing phase, which is denoted as \textbf{(a)} within \textbf{(C)} in Fig.~\ref{fig:alg-setup}. The nonvectorized version directly adopts the original algorithm of Miehe et al.~\cite{MiApLa2002} presented in Fig.~\ref{fig:Miehe-alg}. The vectorized version presented as a function in Code block~\ref{code77} avoids multiple nested for-loops and performs the operations at all integration points simultaneously. It operates with 3-D arrays of the eigenvalues and eigenvectors of $\C$ and of other auxiliary variables from Fig.~\ref{fig:Miehe-alg} and introduces the logical arrays \texttt{inds1}, $\hdots$, \texttt{inds5}. These determine which relations between the three eigenvalues $\lambda_i, \lambda_j, \lambda_k$ occur at each integration point, see the alternatives listed in step 2. in Fig.~\ref{fig:Miehe-alg}. Evaluation of $\PP$ and the contraction $\Sigb : \LL$ is approximately twice as fast, while complete postprocessing is accelerated by approximately a factor of 2.5. The computational time in Fig.~\ref{fig:timesACDE}(c) not covered in Fig.~\ref{fig:timesACDE}(d) is spent on determining variables $\P$ and $\DF$ from \eqref{eqsP} and \eqref{eqsDF}, respectively. That part mainly contains basic linear algebra operations, and the speedup there comes from the vectorization functions in Subsection~\ref{subsec:vectorization}. 

To put these results into the context of the complete Newton iteration procedure, i.e., step \textbf{(ii)} of stage \textbf{2.} in Fig.~\ref{fig:alg-setup}, we present the relative contributions of preprocessing \textbf{(A)}, constitutive problem \textbf{(B)}, postprocessing \textbf{(C)}, global tangent matrix assembly in \textbf{(D)}, linear system solution \textbf{(E)}, and the remaining operations in Fig.~\ref{fig:timesperc}; contributions are expressed as a percentage of the total Newton time for each corresponding mesh and variant.\footnote{Repeated benchmark runs show that local irregularities exhibited on some coarser meshes are related to run-to-run timing variability.} As the mesh is refined, solving the global linearized systems using MATLAB's backslash operator becomes the principal computational expense. Its relative contribution is larger in the vectorized implementation (93--94\% compared to 87--88\% on the finest mesh) because vectorization reduces the costs of preprocessing and postprocessing. The substantial speedups of these individual stages translate into a more modest reduction in the whole simulation time. Further substantial acceleration on larger meshes therefore requires improvements to the linear solution stage or modifications of the Newton iterative procedure (cf. \cite{Karatson2025, Sysala2025}), which are beyond the scope of the present work.

The von Mises elastoplastic routine adopted from \cite{CeSyVa2019} represents a relatively inexpensive component of the present benchmark. On the finest mesh, constitutive evaluation accounts for approximately 0.2\% and 0.6\% of the Newton time in the vectorized elastic and plastic calculations, respectively. This finding concerns the particular constitutive model and implementation considered here, which features a vectorized closed-form solution of the constitutive update. Thanks to the modular separation of constitutive evaluation from preprocessing and postprocessing, other constitutive models formulated within the internal variable framework can be incorporated into the same finite strain formulation (e.g., plasticity with damage \cite{HoJi2013, FrMeSt2023} or shape memory alloys \cite{MHVKJSF2025}), with additional field equations where required. For more complex models, constitutive evaluation may account for a considerably larger share of the total computation time. Its cost may be reduced through vectorization, parallelization, or tailored numerical approaches, e.g., \cite{Lewandowski2023}.

Finally, Fig.~\ref{fig:timesGN} provides the total computation time spent on Newton iterations and on solution of the linear system in our vectorized implementation, both for elastic (a) and plastic (b) models. It confirms that for both constitutive models the linear system provides the most severe computational burden and the absolute linear solution times remain similar between implementations. To quantify the scaling observed in the plots, we approximate the measured times by power laws $t\propto q^p$, where $q=8s^3$ is the number of integration points, over refinement indices $s=10,\ldots,\nmesh{}$. For preprocessing, postprocessing (cf. Fig.~\ref{fig:timesACDE}), constitutive evaluation, and global tangent matrix assembly, least-squares fits in logarithmic coordinates yield exponents ranging between 0.93 and 1.04, indicating approximately linear growth of time with problem size. For the accumulated Newton iteration times in elasticity and plasticity, fits yield exponents of 1.76 and 1.79, respectively, and for the corresponding linear-solver the exponents are approximately 1.84 and 1.88; the fitted power laws are plotted by dashed lines in Fig.~\ref{fig:timesGN}.\footnote{Let us note the resulting exponents describe the \emph{overall} growth over that range rather than a uniform increase between successive refinements.} \corr{Thus, the linear scaling of pre- and postprocessing confirms the robustness of their vectorization and contributes to the subquadratic growth of the total Newton-iteration time.} Let us note for completeness that the remaining processes of the computational workflow, i.e., the first and third stages in Fig.~\ref{fig:alg-setup} as well as parts \textbf{(i)} and \textbf{(iii)} of the second one, require negligible computational time ($<1\%$) for finer computational meshes (refinement index $s>5$).

\subsection{Computational simulations of structural elements}  \label{subsec:models}

We present simulations of the elastoplastic response of two geometric domains, an $\Omega$-shaped and a coil-shaped, which are motivated by structural elements in applications. The examples are chosen to demonstrate the overall performance of the implementation and, at the same time, the ability of the vectorized version of Miehe's approach to handle large rotations and large strains in computational analysis. In visualizations, we use the engineering scalars called \emph{von Mises stress} and \emph{equivalent plastic strain}\footnote{We note that $E^{\rm p}_{\textnormal{eq}}$ measures the magnitude of the current plastic-strain tensor; it should not be confused with the \emph{accumulated} equivalent plastic strain $\int_0^t\sqrt{2/3}\,\|\dot{\E}^{\rm p}(s)\|\,\md s$, which also records plastic flow during load reversals.}, which are defined as
\begin{equation}
    \sigma_\textnormal{vM}
    := \sqrt{\frac{3}{2}\operatorname{dev} \sigb :\operatorname{dev}\sigb},
    \qquad
    E^{\rm p}_\textnormal{eq}
    := \sqrt{\frac{2}{3}\operatorname{dev}\E^{\rm p}:\operatorname{dev}\E^{\rm p}},
    \quad \text{ with }
    \operatorname{dev}\boldsymbol{A}
    := \boldsymbol{A}-\frac{1}{3}\operatorname{tr}(\boldsymbol{A})\,\boldsymbol{I}.
\end{equation}
 Let us note that in the logarithmic strain approach, the yield condition is imposed on the relative \emph{logarithmic} stress (cf. \cite{HoBrReSi2022}), i.e., on $\boldsymbol{\Xi}=\Sigb-\boldsymbol X$ in our case. Hence, the Cauchy von Mises stress (obtained by transformation from $\Sigb$) cannot be used as a direct indicator of plastic activities. However, it can provide a fast clue for assessing the situation in monotonic loading paths starting in an unplasticized state.

\subsubsection{Benchmark 1: $\Omega$-shaped domain} \label{subsubsec:benchmark1}

\begin{figure}
\centering
\begin{minipage}[c]{0.37\textwidth}
\includegraphics[width=\textwidth]{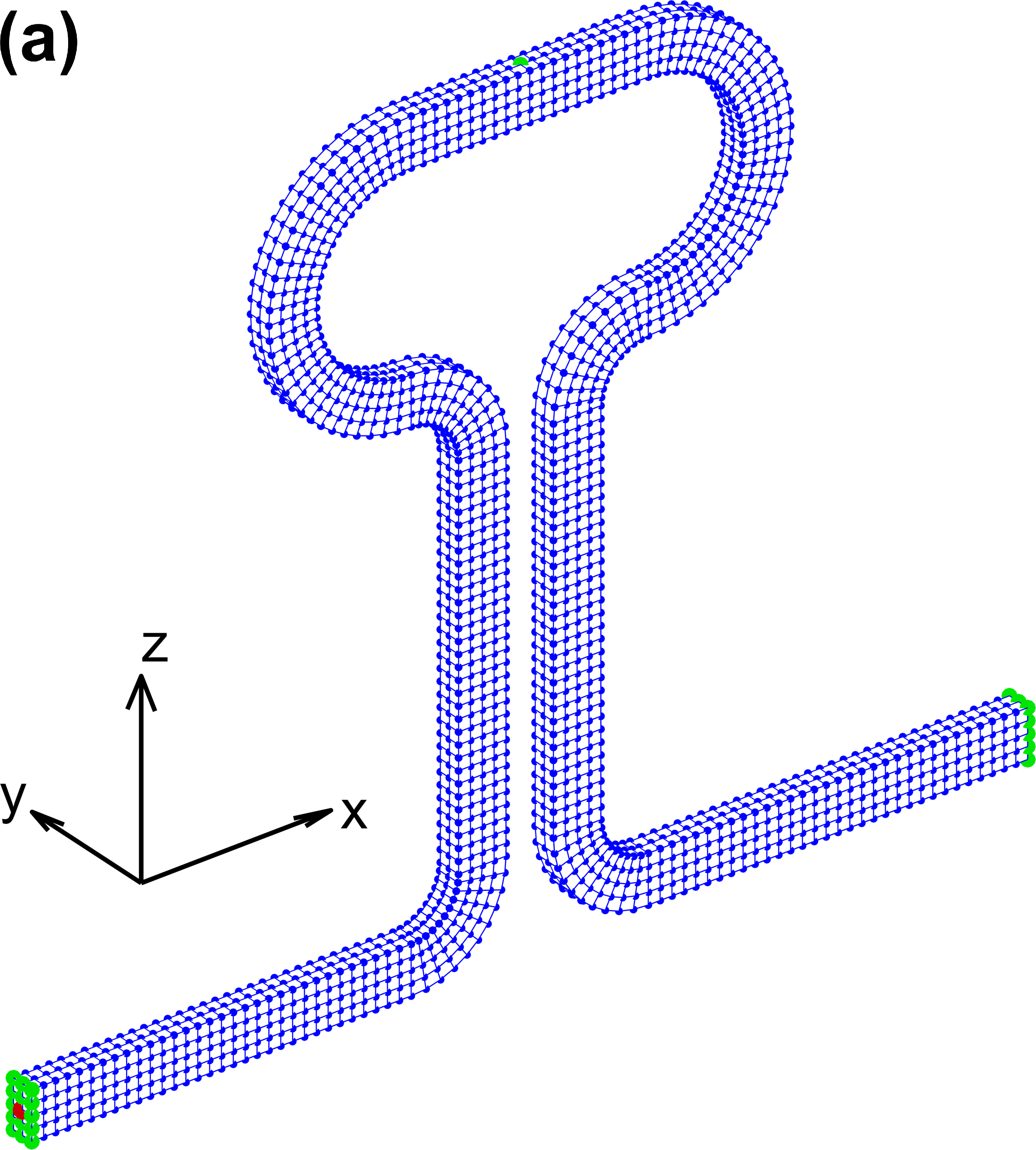}
\end{minipage}
\begin{minipage}[c]{0.58\textwidth}
\includegraphics[width=\textwidth]{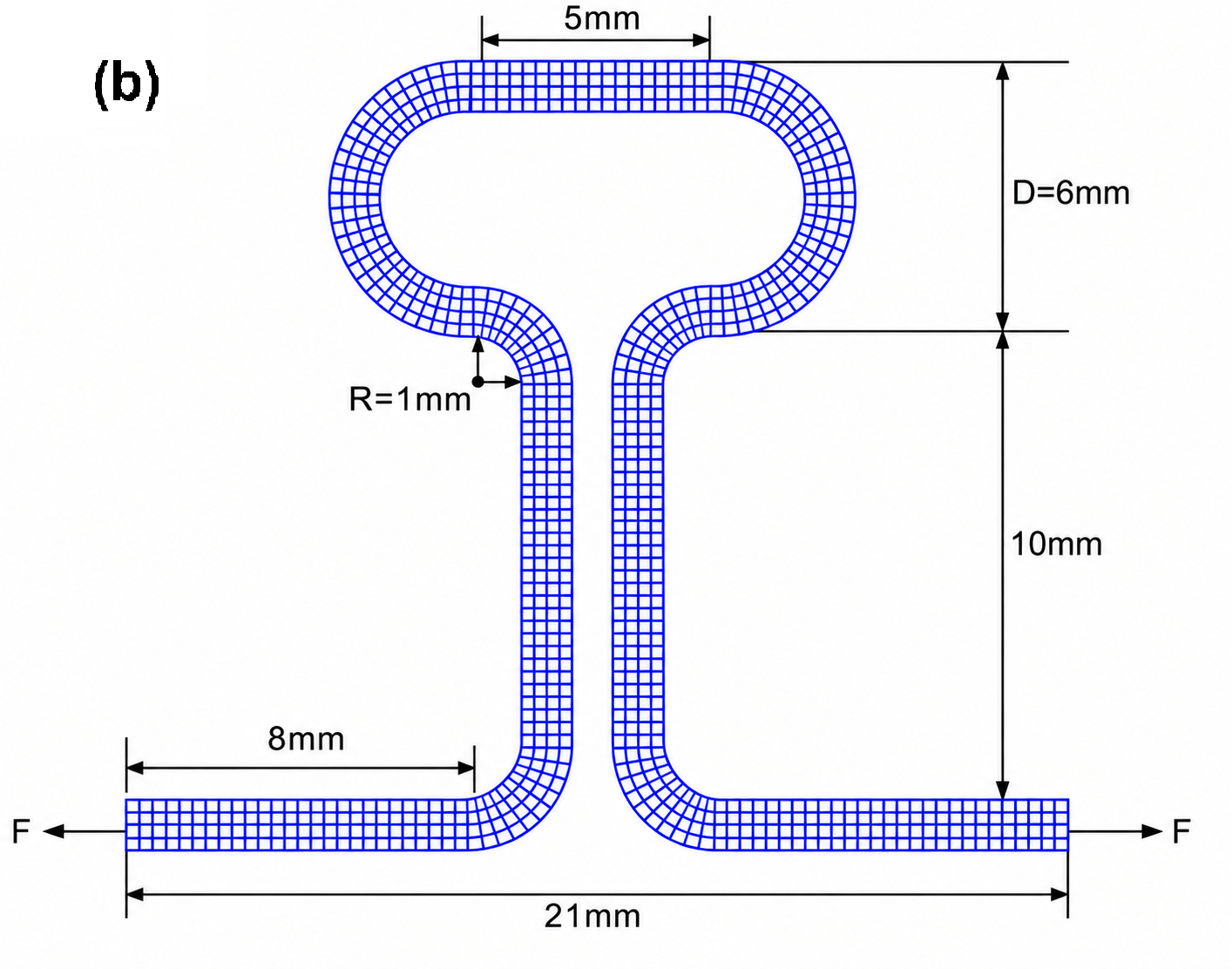}
\end{minipage}
\caption{$\Omega$-shaped domain from the first benchmark: its 3D view (a) and a 2D projection in $xz$-plane (b). A coarse computational mesh is used only for illustrative purposes; see text for details.}
\label{fig:OmegaMesh}
\end{figure}

\begin{figure}
    \centering
    \includegraphics[width=0.92\textwidth]{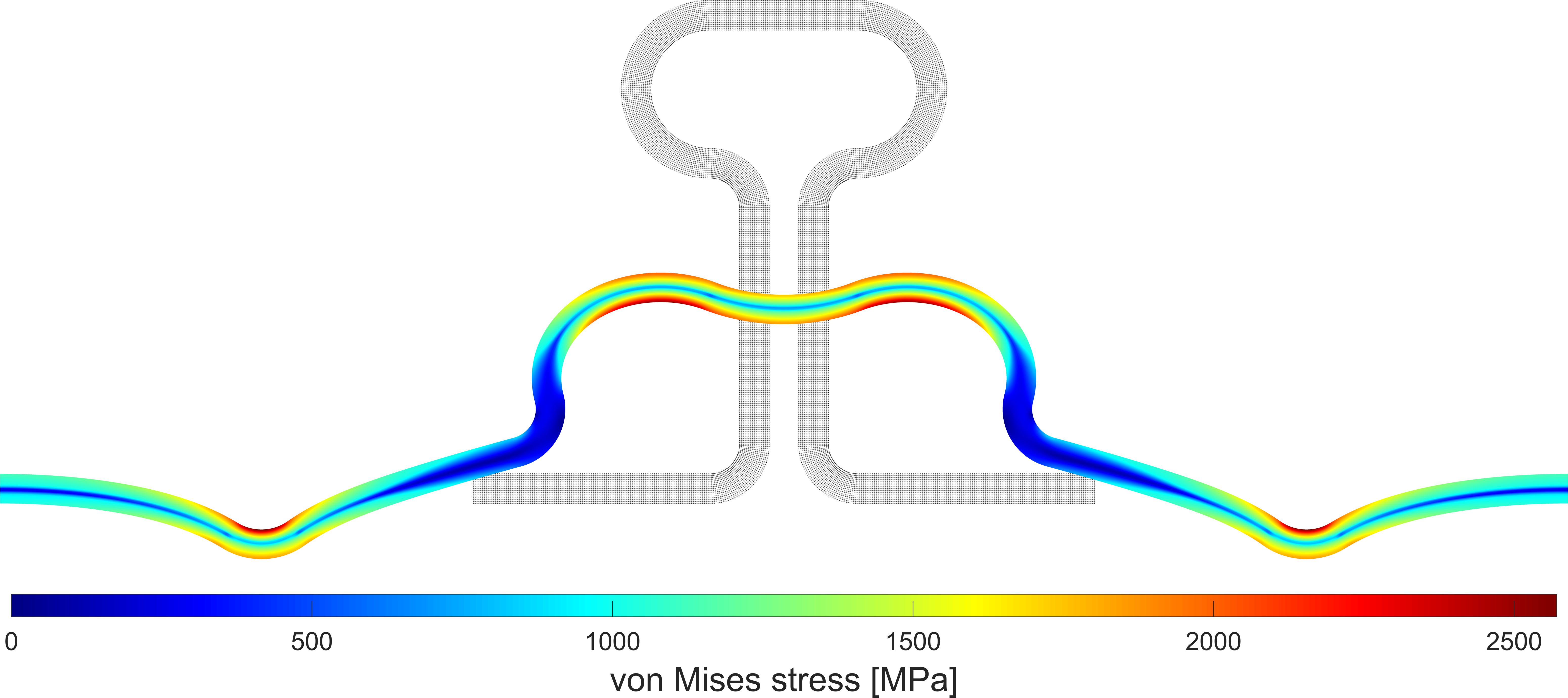}
    \caption{$\Omega$-shaped domain from Benchmark~\ref{subsubsec:benchmark1}: the initial configuration (grey mesh) and the final deformation with the underlying von Mises stress (in MPa).}
    \label{fig:OmegaSol_CAU}
\end{figure}

\begin{figure}
    \centering
    \includegraphics[width=0.92\textwidth]{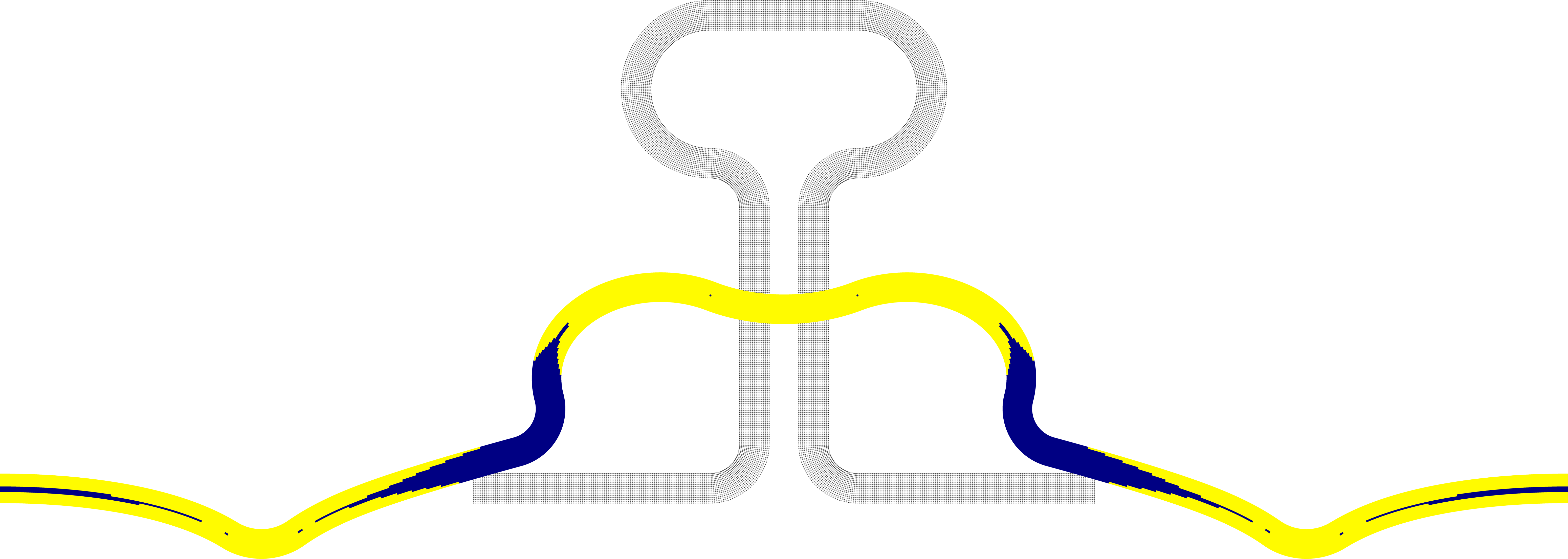}
    \caption{$\Omega$-shaped domain from Benchmark~\ref{subsubsec:benchmark1}: the initial configuration (grey mesh) and the final deformation, with elements where plasticity has been activated marked in yellow and elastically loaded ones marked in blue.}
    \label{fig:OmegaSol_EP}
\end{figure}

\begin{figure}
    \centering
    \includegraphics[width=0.92\textwidth]{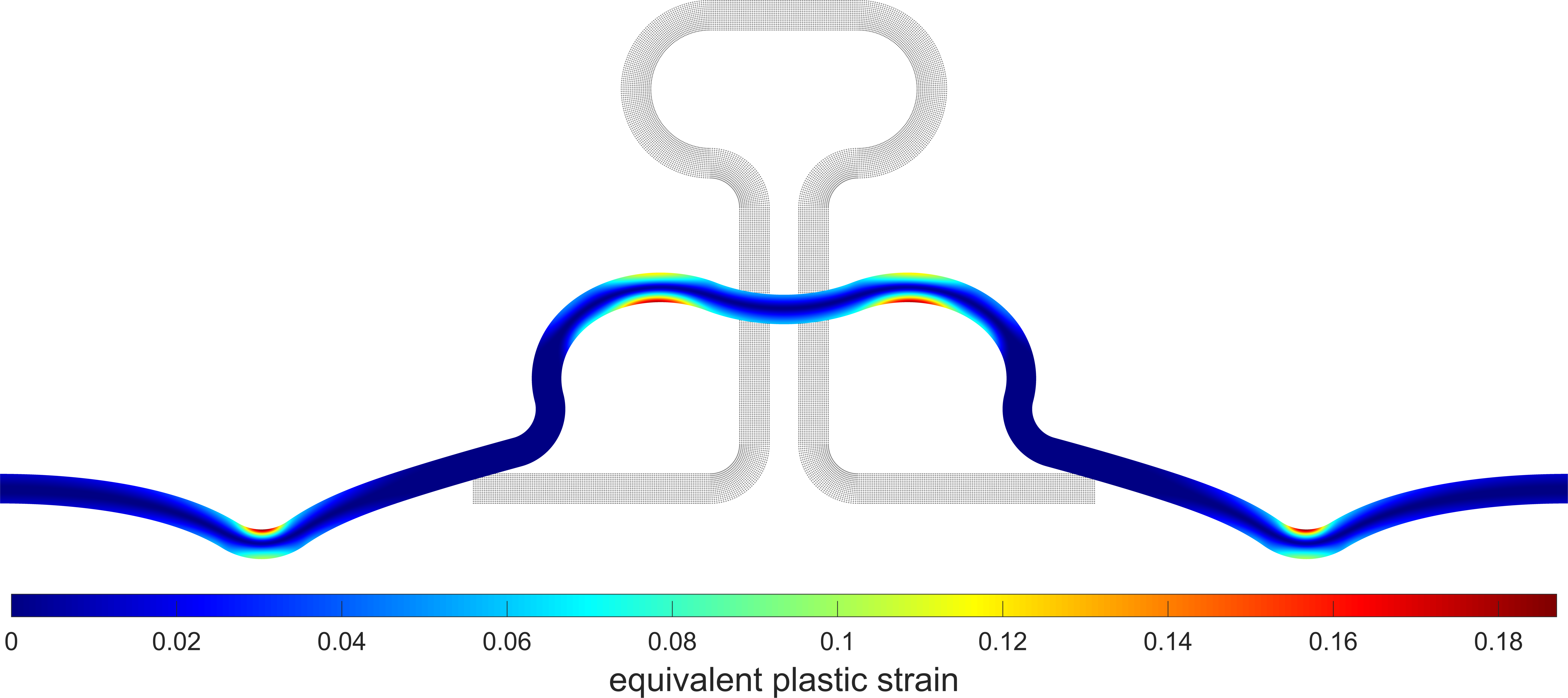}
    \caption{$\Omega$-shaped domain from Benchmark~\ref{subsubsec:benchmark1}: the initial configuration (grey mesh) and the final deformation with the underlying equivalent plastic strain.}
    \label{fig:OmegaSol_EPvM}
\end{figure}

We consider a three-dimensional $\Omega$-shaped domain; see Fig.~\ref{fig:OmegaMesh}(a), which can be considered as a building element of a long, recurring pattern in the $x$-direction. Its projection onto the $xz$-plane together with the specification of dimensions is depicted in Fig.~\ref{fig:OmegaMesh}(b). The width in the $y$-direction is set to $0.5$~mm. The outer left and right faces, whose nodes are marked by green and red points in the figure ($x=-10.5$~mm and $x=10.5$~mm, respectively), are moved apart. This is prescribed by nonhomogeneous Dirichlet boundary conditions in the $x$-direction with maximum displacements of $-u_{\rm max}$ and $+u_{\rm max}$ on the left and right sides, respectively, where $u_{\rm max}=16$~mm. Additionally, the central nodes of both faces are fixed in the $y$- and $z$-directions. To prevent rigid-body motion, the central node of the middle (symmetry) plane, $x=0$, is also fixed in the $y$-direction. The computational mesh consists of 64,512 hexahedral elements (4 cross-sectional layers with 16,128 elements each) and 516,096 corresponding (Gauss) integration points. For visual clarity, a coarser mesh with only 2016 elements (corresponding to \(1/32\) of the number of elements in the computational mesh) is superimposed in Fig.~\ref{fig:OmegaMesh}. The entire simulation consists of 240 uniform time steps with material parameters given in Table~\ref{tab:mat-params}, $a = 10$ GPa.

Figs.~\ref{fig:OmegaSol_CAU}, \ref{fig:OmegaSol_EP}, and~\ref{fig:OmegaSol_EPvM} show the deformed configuration at the final time step, along with the initial mesh configuration in grey. Thanks to the geometry and boundary conditions, the results do not vary with the thickness ($y$-direction) of the $\Omega$-shaped element. In Fig.~\ref{fig:OmegaSol_CAU}, the deformed configuration is overlaid with the distribution of von Mises stress. The most loaded parts are the arcs situated at the bottom and top of the structure, with equivalent stress reaching more than 2.5~GPa on the surfaces. From the mechanics perspective, the inner surfaces of the arcs are in tension, whereas the outer surfaces are compressed. The plotted Cauchy von Mises stress visualizes the intensity of the current stress state. It is not itself the decisive factor for the plastic flow -- the scalar function $\sqrt{3/2}\,\|\operatorname{dev}(\Sigb-\boldsymbol X)\|$ plays this role for the present model. However, for monotonic loading starting with no plastic strain, the comparison of von Mises stress with the yield stress provides a clue for identification of critical regions of the material. Fig.~\ref{fig:OmegaSol_EP} shows elements with nonzero current equivalent plastic strain, $E^{\rm p}_{\textnormal{eq}} > 0$, in yellow, and the remaining elements in blue; this identifies regions with plastic-strain history. Finally, as documented in Fig.~\ref{fig:OmegaSol_EPvM}, where the distribution of equivalent plastic strain is shown, the intensity of plastification varies significantly. The equivalent plastic strain grows rapidly near the four arcs, reaching maxima of about $0.19$ on the stretched inner surfaces. Such detailed information on the presence and intensity of the plastic process may be useful in designing real applications for optimizing their functionality and lifespan. 

\subsubsection{Benchmark 2: a coil} \label{subsubsec:benchmark2}

\begin{table}[htbp]
    \centering
    \renewcommand{\arraystretch}{1.25}
    \setlength{\tabcolsep}{10pt}
    \begin{tabular}{|c|c|c|}
        \hline
            {Snapshot}
            & $u_z^{\rm SN}$ [cm]
            & $u_z^{\rm SN}/u_z^{\rm SN, max}$ [1] \\
        \hline
        $S_1$ & 13.2 & $\sim 1/3$ \\
        \hline
        $S_2$ & 26.6 & $\sim 2/3$ \\
        \hline
        $S_3$ & 40.0 & $1$ \\
        \hline
        $S_4$ & 32.1 & $\sim 4/5$ \\
        \hline
    \end{tabular}
    \caption{Displacement values of the moved node of the deformed coil for the four snapshots discussed in the text.}
    \label{tab:snapshots}
\end{table}

\begin{figure}
\centering
\begin{minipage}[c]{0.55\textwidth}
\includegraphics[width=\textwidth]{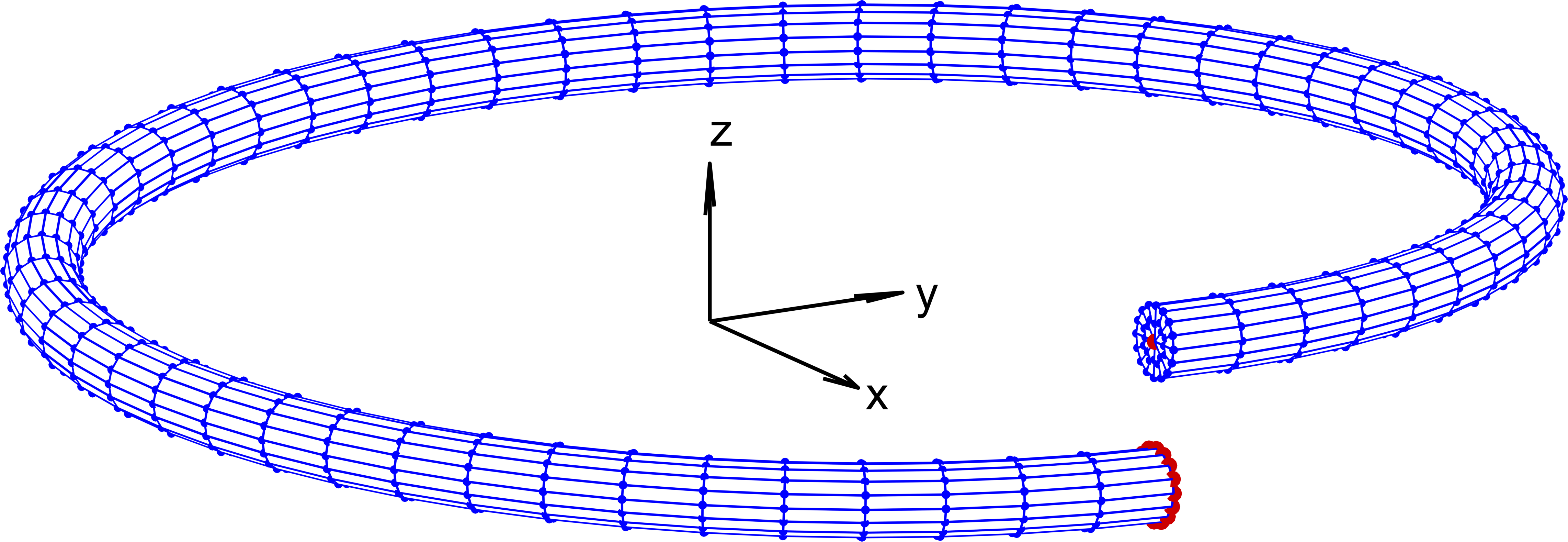}
\end{minipage} \qquad\qquad
\begin{minipage}[c]{0.22\textwidth}
\includegraphics[width=\textwidth]{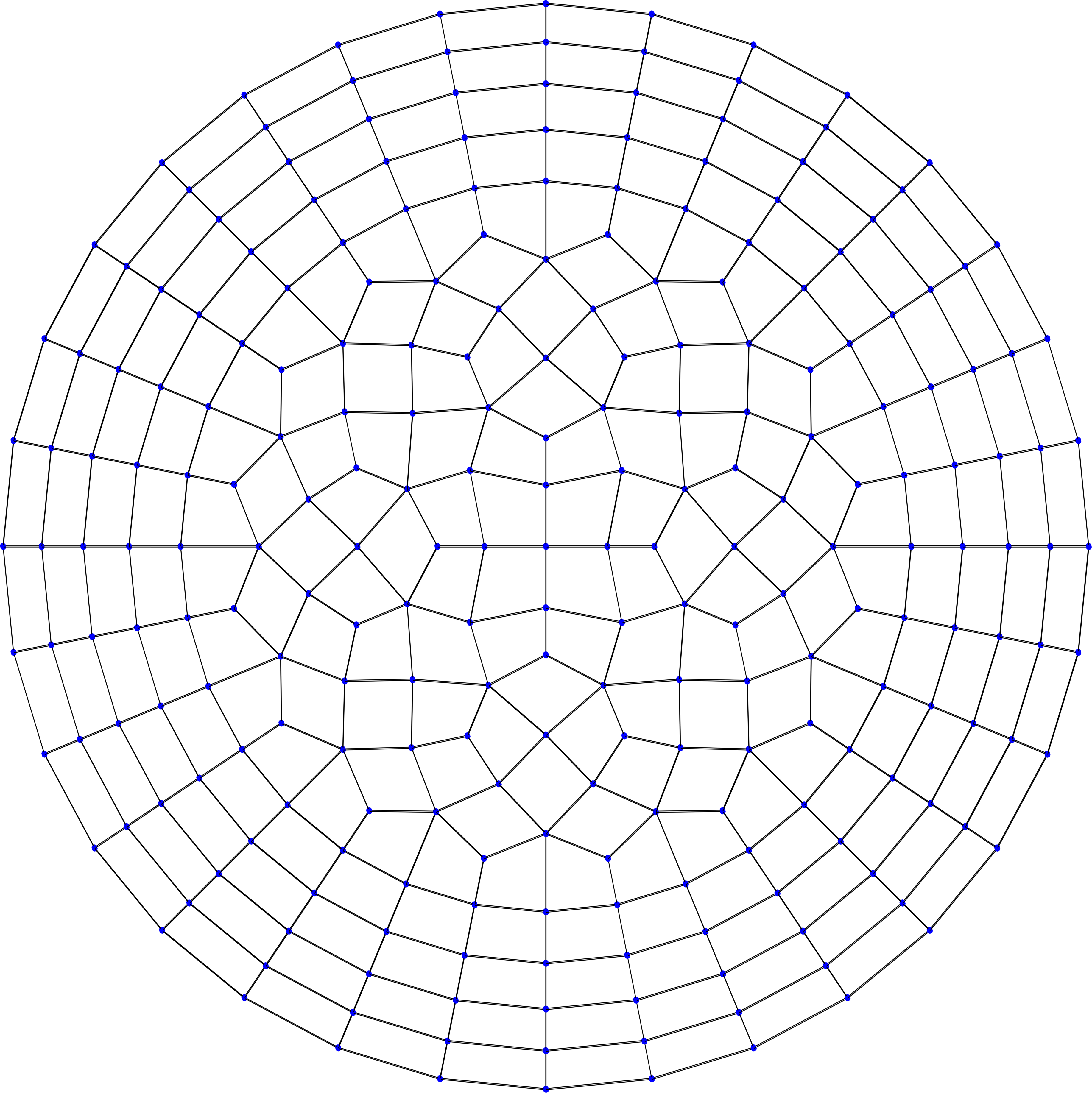}
\end{minipage}
    \caption{Coil domain from Benchmark~\ref{subsubsec:benchmark2} (left) and its circular cross-section (right) with the geometry of the FE mesh in black.}
\label{fig:springMesh2}
\end{figure}

\begin{figure}
    \centering
    \includegraphics[width=0.9\textwidth]{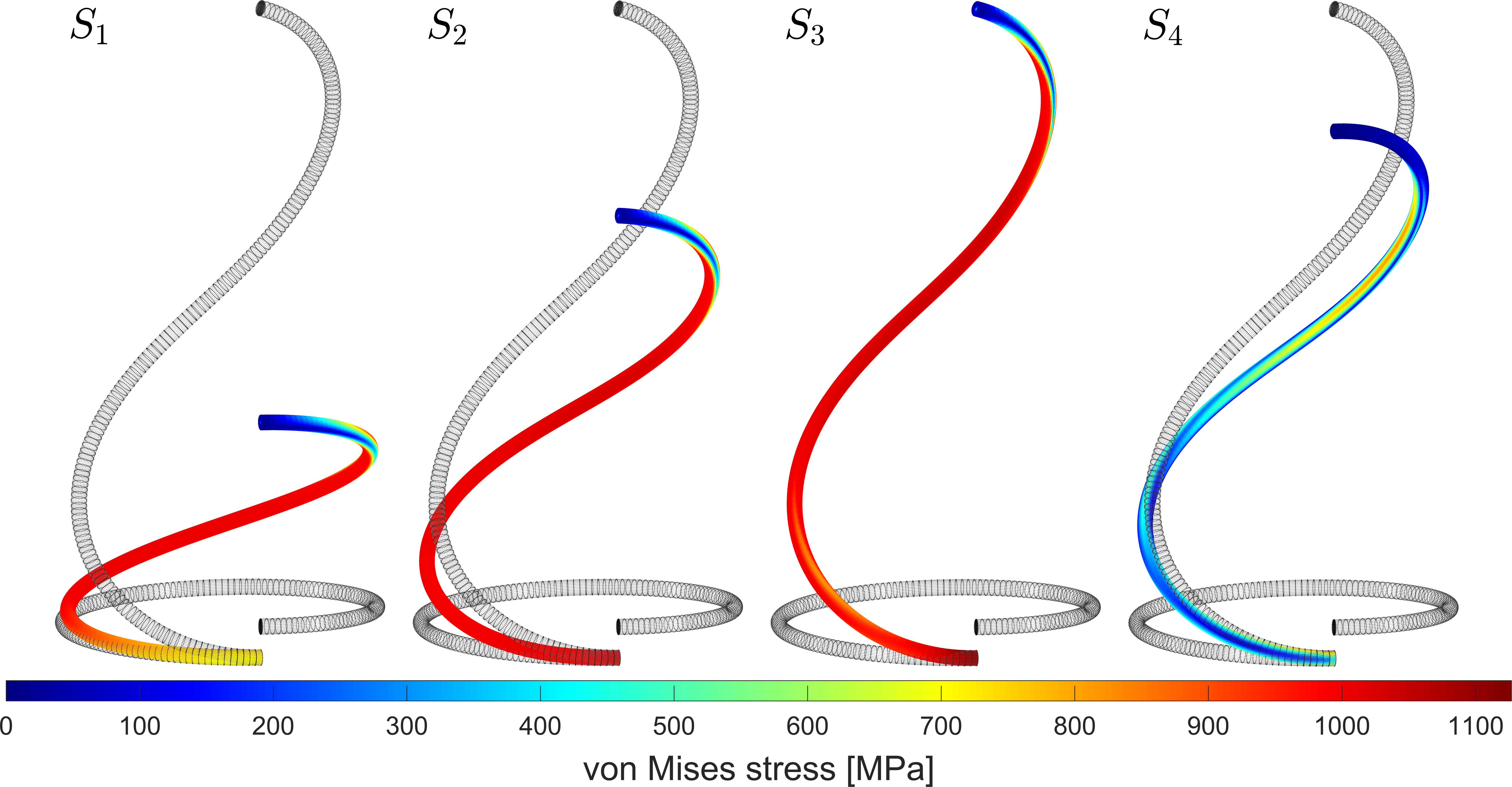}
    \caption{Benchmark~\ref{subsubsec:benchmark2}: von Mises stress (in MPa) at different stages of stretching $S_1-S_4$. See text for details.}
    \label{fig:coil-3D}
\end{figure}

\begin{figure}
    \centering
    \includegraphics[width=0.8\textwidth]{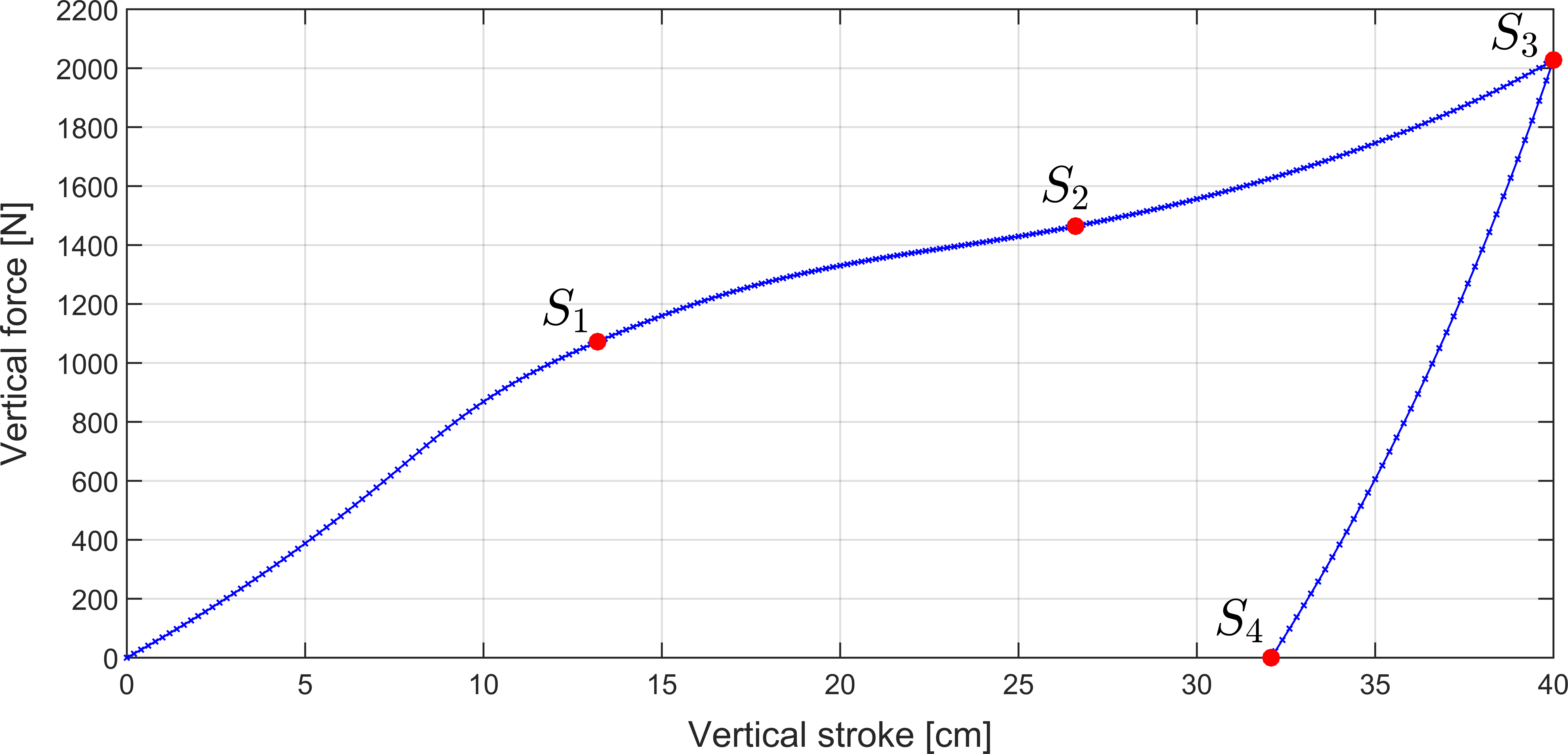}
    \caption{Benchmark~\ref{subsubsec:benchmark2}: the force-stroke relation with four particular stages of stretching $S_1-S_4$ marked by red points. See text for details.}
    \label{fig:springSol}
\end{figure}

\begin{figure}
    \centering
    \includegraphics[width=0.8\textwidth]{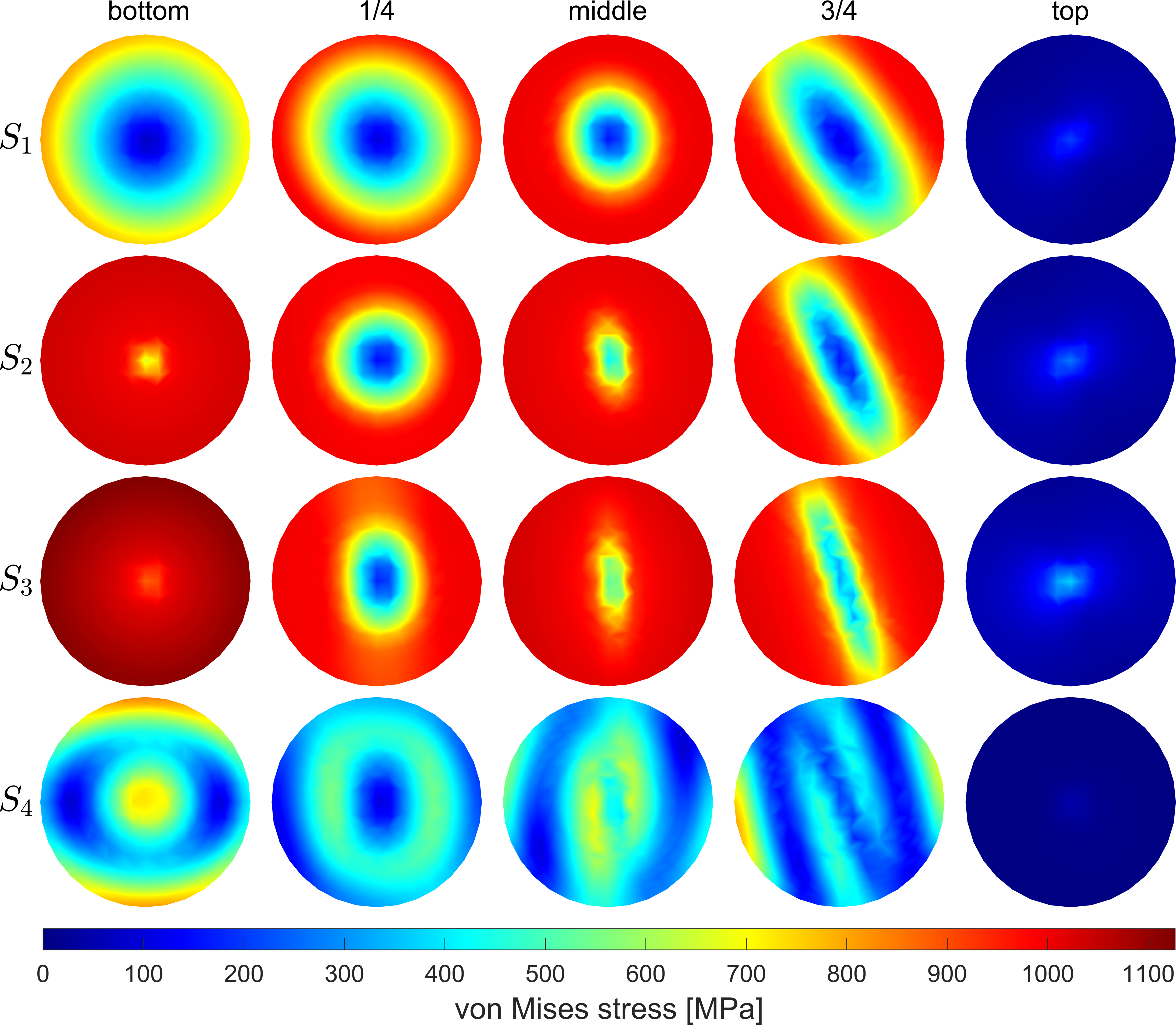}
    \caption{Benchmark~\ref{subsubsec:benchmark2}: von Mises stress (in MPa) on particular cross-sections (bottom, 1/4, middle, 3/4, and top with respect to the central fiber) for stages of stretching $S_1-S_4$. See text for details.}
    \label{fig:coil-cross-stress}
\end{figure}

In this benchmark, we consider a three-dimensional coil (ring) with a circular cross-section. The coil cross-section radius is 5~mm, the outer radius is 10~cm, and its height, i.e., the vertical distance between the central nodes of the two end cross-sections, is 2~cm. The computational domain is depicted in Fig.~\ref{fig:springMesh2} on the left. All nodes of the bottom end cross-section are fixed in all three spatial directions. The coil is loaded by prescribing 
displacement of the central node of the upper cross-section (marked in red) in the vertical $z$-direction, while its $x$- and $y$-displacements remain fixed. The coil is stretched up to $u_z = 40$~cm, then the stretch is released until the $z$-component of the force at the central node of the upper cross-section vanishes, which occurs at $u_z \approx 32.1$~cm. The computational mesh consists of 36,480 hexahedral elements, arranged in 160 layers along the main axis with 228 elements per layer (a coarser mesh is depicted in the left part of Fig.~\ref{fig:springMesh2}). The entire simulation comprises 240 uniform time steps (200 for loading, 40 for unloading). We consider material parameters from Table~\ref{tab:mat-params}; the kinematic hardening coefficient is now decreased to $a = 1$ GPa. 

Fig.~\ref{fig:coil-3D} shows four snapshots from the computation corresponding to the displacement of the elevated node ($u_z^{\rm SN}$) of 13.2~cm, 26.6~cm, 40.0~cm for loading and 32.1~cm for unloading, where the force in the vertical direction reaches zero; see Table~\ref{tab:snapshots} for an overview and Fig.~\ref{fig:springSol} for the evolution of the vertical component of force with stroke at the elevated node. In each snapshot, the initial configuration and the most stretched configuration are plotted as the meshed surface in grey, whereas the actual configuration is overlaid with the distribution of von Mises stress. 

Although the macroscopic force response up to $S_1$ in Fig.~\ref{fig:springSol} appears slightly nonlinear, the corresponding snapshot in Fig.~\ref{fig:coil-3D} shows substantial heterogeneity along the coil's central fiber, with very low stress near the moved end and higher stress below. A similar situation occurs during further loading (snapshots $S_1$--$S_3$): the moved end of the coil exhibits low von Mises stress; the surface of the central part is loaded close to $1000\,{\rm MPa}$, whereas the stress in the fixed end of the coil increases from values below $700\,{\rm MPa}$ for $S_1$ to values reaching the maximum of the whole simulation, over $1100\,{\rm MPa}$, for $S_3$. This shows a highly inhomogeneous stress state along the coil's central fiber, driven by the imposed boundary conditions; moreover, the large rotation of mesh elements farther from the fixed end is clearly visible. After unloading (snapshot $S_4$), the stress on the surface drops everywhere, but it now appears distributed inhomogeneously also on the outer edges of the coil cross-sections.

The last finding is further illustrated in Fig.~\ref{fig:coil-cross-stress}, where the distribution of von Mises stress within five cross-sections of the coil is presented for snapshots $S_1-S_4$ from Fig.~\ref{fig:coil-3D}. The five cross-sections were chosen at different relative positions along the central fiber of the coil. Defining the relative position of the fixed end and the moved end as 0 and 1, respectively, the cross-sections at relative positions 0 (=bottom), 1/4, 1/2 (=middle), 3/4 and 1 (=top) are shown. The plots confirm that the inhomogeneous stress distribution also occurs within the cross-section, caused by the combination of torsional and bending strain in this helical-spring-derived geometry (cf. \cite{Wahl1944, KaSu2021, SMS-DSCT}). During stretching, the dominant straining mode shifts from torsion to bending; e.g., see \cite{SMS-DSCT}. Hence, the circular symmetry of the stress distribution is broken, and mirror symmetry imposed by bending appears. This is clearly seen for the cross-section denoted $3/4$ and for the stress distribution at $S_4$. After unloading, the stress distribution is heavily influenced by the permanent plastic strain generated during stretching; see also Appendix~\ref{app:coil-cross-plast-strain}. Let us note that the inhomogeneous residual stresses can be investigated in more detail using dedicated FEM techniques \cite{Carstensen2016, Bammer2025}.

\section{Conclusions} \label{sec:conclusions}

More than two decades after the seminal work \cite{MiApLa2002}, the logarithmic strain framework has become a popular choice for upgrading engineering small-strain material models toward large deformations. The intrinsic multiplicative nature of finite-strain decomposition is mapped by a logarithmic function onto the additive structure of geometrically linear theory, allowing seamless integration of additive constitutive models without modification and shifting the effort to constitutive-core-independent, purely geometric pre- and post-processing phases. We have studied this shift both from the theoretical mathematical and practical computational perspectives, taking von Mises elastoplasticity with linear kinematic hardening as a clear and well-known prototype.

In the mathematical analysis, we have established the existence of energetic solutions for the gradient-regularized quasistatic model, under the assumptions of Theorem~\ref{thm:energetic}. The key challenge was addressing the elastic energy inheriting the additive nature from the small-strain ancestor; since this is typical for this class of models, the applied approach can serve as a methodological template for other class instances. In the computational part, we focused on the effective algorithmic treatment of tensor transformations in the pre- and post-processing phases, which are independent of the constitutive core and therefore vital for any practical implementation. We exploited the capabilities of the advanced vectorization tools and our own developed techniques in MATLAB code \cite{RaVa2013, MoVa2022, MoRaVaVa2026} and significantly accelerated the pre- and post-processing phases so that, for medium-sized FEM problems, the decisive computational limit is imposed by the solution of the linear system within the Newton method and potentially the constitutive update routine in case of more complex constitutive laws, not these vectorized phases. Let us note for completeness that we built on the implementation with optimized tangent matrix assembly \cite{CeSyVa2019}. \corr{In a broader context, this work also illustrates the versatility of our vectorization techniques: here they are applied to a geometrically nonlinear problem whose linearizations need not be elliptic, which extends the range of problems treated in our previous works \cite{RaVa2013, CeSyVa2019, MoVa2022}.}

The benchmark simulations demonstrated the applicability of the FEM implementation to resolving the mechanical response of geometries undergoing large displacements and rotations. They also illustrated that detailed information on physical quantities, e.g., stresses and internal variable(s) (plastic strains), can be provided. Let us note that advanced characterization techniques can provide experimental counterparts to such a type of data (e.g., \cite{SMS-DSCT, SMS-SNAKE, tong2004, chandra2024}), opening the possibilities of experimental validation of model predictions. 

The separation of geometric transformations from constitutive
evaluation facilitates extensions to other material models compatible with the logarithmic strain framework. Such extensions require appropriate constitutive equations and internal-variable updates, while the geometric preprocessing and postprocessing can be retained when the strain measure and stress conjugacy remain unchanged. Models involving additional fields and/or spatial gradients require corresponding extensions of the global formulation \cite{HoBrReSi2022, FrMeSt2023}.

\section*{Declaration of competing interest}
The authors declare that they have no known competing financial interests or personal relationships that could have appeared to influence the work reported in this paper.

\section*{Data availability}
Data are available at \cite{data-Zenodo}. The software is available at \cite{code-Zenodo}.

\section*{Acknowledgements}
This work has been financially supported by
the Czech Science Foundation [Project No. 24-10366S] and by the
Operational Programme Johannes Amos Comenius of the Ministry
of Education, Youth and Sport of the Czech Republic, within the
frame of project Ferroic Multifunctionalities (FerrMion) [Project No.
CZ.02.01.01/00/22\_008/0004591], co-funded by the European Union.

\printcredits

\appendix
\section{Appendix}

\subsection{Small-strain constitutive model of von Mises elastoplasticity with linear kinematic hardening} \label{app:von-mises}

Let us recall the constitutive model representing von Mises elastoplasticity with linear kinematic hardening in small strains. The state of the material is defined by the total strain $\boldsymbol{\varepsilon}$ (observable) and the plastic strain $\boldsymbol{\varepsilon}^\textnormal{p}$ (internal variable):
\begin{equation}
    \boldsymbol{\varepsilon} = \boldsymbol{\varepsilon}^{\rm el} + \boldsymbol{\varepsilon}^\textnormal{p}
\end{equation}
The free energy, $\psi$, is decomposed into an elastic part, $\psi^{\rm el}$, and a stored energy part due to the linear kinematic hardening, $\psi^{\rm kin}$:
\begin{equation}
    \psi(\boldsymbol{\varepsilon}, \boldsymbol{\varepsilon}^\textnormal{p}) = \psi^\textnormal{el}(\boldsymbol{\varepsilon} - \boldsymbol{\varepsilon}^\textnormal{p}) + \psi^\textnormal{kin}(\boldsymbol{\varepsilon}^\textnormal{p}) = \frac{1}{2} (\boldsymbol{\varepsilon} - \boldsymbol{\varepsilon}^\textnormal{p}) : \mathbb{C} : (\boldsymbol{\varepsilon} - \boldsymbol{\varepsilon}^\textnormal{p}) + \frac{1}{2} H^{\rm kin} \boldsymbol{\varepsilon}^\textnormal{p} : \boldsymbol{\varepsilon}^\textnormal{p}
\end{equation}
where $\mathbb{C}$ is the fourth-order elastic stiffness tensor and $H^{\rm kin}$ is the kinematic hardening modulus.

The thermodynamic forces are derived from the partial derivatives of the free energy:
\begin{align}
    \sigb &= \frac{\partial \psi}{\partial \boldsymbol{\varepsilon}} = \mathbb{C} : (\boldsymbol{\varepsilon} - \boldsymbol{\varepsilon}^\textnormal{p}) \\
    \boldsymbol{\beta} &= \frac{\partial \psi^{\rm kin}}{\partial \boldsymbol{\varepsilon}^\textnormal{p}} = H^{\rm kin} \boldsymbol{\varepsilon}^\textnormal{p} \\
    \boldsymbol{\chi} &= -\frac{\partial \psi}{\partial \boldsymbol{\varepsilon}^\textnormal{p}} = \sigb - \boldsymbol{\beta}
\end{align}
Here, $\sigb$ is the Cauchy stress, $\boldsymbol{\beta}$ is the backstress, and $\boldsymbol{\chi}$ is the relative stress (the force driving dissipation).

The dissipation (pseudo)potential takes the form
\begin{equation}
\phi(\dot{\boldsymbol{\varepsilon}}^{\rm p})  = 
\begin{cases}
\sqrt{2/3}\,Y\|\dot{\boldsymbol{\varepsilon}}^{\rm p}\| & \text{if }\;{\rm tr}(   \dot{\boldsymbol{\varepsilon}}^{\rm p}) = 0, \\
+\infty & \text{else}.
\end{cases}
\end{equation}
with $Y$ denoting the yield stress for uniaxial loading. Now, the solution of the constitutive problem (state update) can be formally written as a solution of an extremal problem \cite{HaFi2008}
\begin{equation}
   \operatorname*{min}_{\dot{\boldsymbol{\varepsilon}}^{\rm p}}
\left\{ \dot{\psi} + \phi \right\}.
\end{equation}

Although this is a less conventional formulation, it is fully equivalent to the ``classical'' alternative, e.g., \cite{SimHug98CI}, written as follows. First, the elastic domain $\mathcal{K}$ is defined by the von Mises yield criterion:
\begin{equation}
    f(\boldsymbol{\chi}) = \sqrt{\frac{3}{2} \operatorname{dev} \boldsymbol{\chi} : \operatorname{dev}\boldsymbol{\chi}} - Y \leq 0
\end{equation}
The evolution of the internal variable $\boldsymbol{\varepsilon}^\textnormal{p}$ follows the normality rule
\begin{equation}
    \dot{\boldsymbol{\varepsilon}}^\textnormal{p} = \dot{\gamma} \frac{\partial f}{\partial \boldsymbol{\chi}},
\end{equation}
which results in the associative flow rule:
\begin{equation}
    \dot{\boldsymbol{\varepsilon}}^\textnormal{p} = \dot{\gamma}\frac{3}{2} \frac{{\rm dev}(\boldsymbol{\chi})}{Y}
\end{equation}
subject to the Karush-Kuhn-Tucker conditions:
\begin{equation}
    \dot{\gamma} \geq 0, \quad f(\boldsymbol{\chi}) \leq 0, \quad \dot{\gamma} f(\boldsymbol{\chi}) = 0.
\end{equation}

\subsection{Existence of energetic solutions}
\label{app:en-solution}

We formulate the theorem on the existence of energetic solutions for the elastoplastic model introduced in Section~\ref{sec:maths} and briefly sketch its proof. We refer the reader to the classical works \cite{ortiz-repetto, mielke1, ball81} for further details of the technique.

\noindent
\begin{theorem}\label{thm:energetic}
Let \(\Omega\subset\R^3\) be bounded with a sufficiently regular boundary and let \(\Gamma_\textnormal{D}\subset\partial\Omega\) have positive surface measure. Assume
\begin{equation}
    r>3,
    \qquad
    q\ge2,
    \qquad
    s>\frac{3r}{r-3},
\end{equation}
and let \(\varepsilon,\kappa>0\). For the quadratic elastic energy with linear kinematic hardening used in the applications, the growth assumption below is satisfied with $q=2$.

Assume that
\begin{equation}
    \Psi^\textnormal{el}\in C^1(\Sym(3))
\end{equation}
and that there are constants \(c,C>0\) such that
\begin{equation}
    c|\boldsymbol Z|^2-C
    \le
    \Psi^\textnormal{el}(\boldsymbol Z)
    \le
    C(1+|\boldsymbol Z|^2)
\end{equation}
and
\begin{equation}
    |D\Psi^\textnormal{el}(\boldsymbol Z)|
    \le
    C(1+|\boldsymbol Z|)
    \qquad
    \text{for all }\, \boldsymbol Z\in\Sym(3).
\end{equation}

Assume that
\begin{equation}
    \Psi^{\rm kin}:\Symzero(3)\to[0,\infty]
\end{equation}
is convex, lower semicontinuous, and coercive in the sense that
\begin{equation}
    \Psi^\textnormal{kin}(\boldsymbol E^\textnormal{p})
    \ge
    c|\boldsymbol E^\textnormal{p}|^q-C.
\end{equation}

Assume that $\Phi$ is convex, positively one-homogeneous, and satisfies
\begin{equation}
    c_\Phi|\boldsymbol Z|
    \le
    \Phi(\boldsymbol Z)
    \le
    C_\Phi|\boldsymbol Z| \qquad
    \text{for all }\,\boldsymbol Z\in \Symzero(3)\, \text{ with }\, C_\Phi \geq c_\Phi > 0.
\end{equation}

Let
\begin{equation}
    \boldsymbol g(t,\boldsymbol x)=\boldsymbol x+\boldsymbol w(t,\boldsymbol x),\qquad
    \boldsymbol w\in C^1([0,T];C_b^2(\R^3;\R^3))
\end{equation}
be such that each \(\boldsymbol g_t\) is an orientation-preserving
diffeomorphism and
\begin{equation}
    \det\nabla\boldsymbol g_t
    \ge
    c_g>0
\end{equation}
uniformly in \(t\). Assume further that
\begin{equation}
    \sup_{t\in[0,T]}
    \left(
        \|\nabla\boldsymbol g_t\|_{L^\infty}
        +
        \|(\nabla\boldsymbol g_t)^{-1}\|_{L^\infty}
        +
        \|\partial_t\nabla\boldsymbol g_t\|_{L^\infty}
    \right)
    <\infty.
\end{equation}

Let the initial state
\begin{equation}
    (\boldsymbol\upsilon_0,\boldsymbol E^\textnormal{p}_0)\in\mathcal Q
\end{equation}
have finite energy and satisfy the initial global stability condition
\begin{equation}
    \mathcal E_\varepsilon
    (0,\boldsymbol\upsilon_0,\boldsymbol E^\textnormal{p}_0)
    \le
    \mathcal E_\varepsilon
    (0,\widehat{\boldsymbol\upsilon},\widehat{\boldsymbol E^\textnormal{p}})
    +
    \mathcal D(
        \boldsymbol E^\textnormal{p}_0,\widehat{\boldsymbol E^\textnormal{p}}
    )
\end{equation}
for all admissible competitors
\((\widehat{\boldsymbol\upsilon},\widehat{\boldsymbol E^\textnormal{p}})\in\mathcal Q\).

Then there exists an energetic solution
\begin{equation}
    (\boldsymbol\upsilon,\boldsymbol E^\textnormal{p}):[0,T]\to\mathcal Q
\end{equation}
such that
\begin{equation}
    \boldsymbol E^\textnormal{p}
    \in
    BV_{\mathcal D}
    \bigl(
        [0,T];
        L^1(\Omega;\Symzero(3))
    \bigr),
\end{equation}
and, for every \(t\in[0,T]\),
\begin{equation}
    \boldsymbol\upsilon(t)\in W^{2,r}(\Omega;\R^3),
    \qquad
    \boldsymbol E^\textnormal{p}(t)\in W^{1,q}(\Omega;\Symzero(3)).
\end{equation}
Moreover,
\begin{equation}
    \det\nabla\boldsymbol\upsilon(t,\boldsymbol x)>0
    \qquad
    \text{for every }\boldsymbol x\in\overline\Omega,
\end{equation}
and the physical deformation
\begin{equation}
    \boldsymbol y(t)
    :=
    \boldsymbol g_t\circ\boldsymbol\upsilon(t)
\end{equation}
satisfies
\begin{equation}
    \boldsymbol y(t,\boldsymbol x)
    =
    \boldsymbol g_t(\boldsymbol x)
    \qquad
    \text{on }\Gamma_\textnormal{D}.
\end{equation}
The energetic solution satisfies the global stability condition and the energy balance stated above.
\end{theorem}

\noindent \textbf{Brief sketch of the proof.} A central estimate used in the proof is
\begin{equation}
    \left|
        D\E(\boldsymbol F)
        [\boldsymbol L\boldsymbol F]
    \right|
    \le
    C|\boldsymbol L|,
    \qquad
    \boldsymbol F\in GL^+(3),
\end{equation}
with a constant independent of \(\boldsymbol F\). Consequently,
\begin{equation}
    \left|
        \partial_t
        \mathcal E_\varepsilon
        (t,\boldsymbol\upsilon,\boldsymbol E^\textnormal{p})
    \right|
    \le
    a(t)
    \left(
        1+
        \mathcal E_\varepsilon
        (t,\boldsymbol\upsilon,\boldsymbol E^\textnormal{p})
    \right)
\end{equation}
for some \(a\in L^1(0,T)\).

To justify the direct-method step, one first uses bounded energy to obtain a uniform $W^{2,r}$ bound on the deformation and a $W^{1,q}$ bound on the plastic variable. The embedding into $C^{1,\alpha}$, together with the inverse-determinant penalty, gives a uniform positive lower bound for the determinant on energy sublevels. Thus, the logarithmic strain is continuous along the strongly convergent deformation gradients. Compactness of the plastic variable then yields convergence of the elastic-energy term; convexity and weak lower semicontinuity of the gradient penalties and the dissipation handle the remaining terms.

The proof is based on the incremental problems
\begin{equation}
    (\boldsymbol\upsilon_\ts,\boldsymbol E^\textnormal{p}_\ts)
    \in
    \Argmin_{\boldsymbol\upsilon,\boldsymbol E^\textnormal{p}}
    \left\{
        \mathcal E_\varepsilon
        (t_\ts,\boldsymbol\upsilon,\boldsymbol E^\textnormal{p})
        +
        \mathcal D(
            \boldsymbol E^\textnormal{p}_{\ts-1},\boldsymbol E^\textnormal{p}
        )
    \right\},
\end{equation}
followed by the direct method, uniform energy estimates, a
\(\mathcal D\)-variation estimate for \(\boldsymbol E^\textnormal{p}\), Helly's selection
principle, and passage to the limit in the discrete stability and energy
inequalities.

The essential compactness mechanism is
\begin{equation}
    W^{2,r}(\Omega;\R^3)
    \Subset
    C^{1,\alpha}(\overline\Omega;\R^3),
    \qquad
    0<\alpha<1-\frac3r.
\end{equation}
Together with
\begin{equation}
    s>\frac{3r}{r-3},
\end{equation}
the inverse-determinant term yields a positive lower bound for
\(\det\nabla\boldsymbol\upsilon\) on bounded energy sublevels. Hence, the logarithmic strain is continuous along minimizing sequences. In particular, let us stress that no quasiconvexity or polyconvexity of
\begin{equation}
    \boldsymbol F
    \mapsto
    \Psi^\textnormal{el}\!\left(
        \frac12\log(\boldsymbol F^\top\boldsymbol F)
        -
        \boldsymbol E^\textnormal{p}
    \right)
\end{equation}
is required for the regularized problem.

While the regularizing terms in \eqref{eq:regul-energy} are necessary for the current analysis, the simulations in Section~\ref{sec:simulations} are performed for the non-regularized model, i.e., formally for $\varepsilon=0$ and $\kappa=0$, which is the formulation commonly accepted and used in computational practice \cite{CeSyVa2019, SimHug98CI}. Let us comment on the relation between the two settings from a theoretical point of view.

\emph{Time discretization.} Energetic solutions are obtained as limits of solutions of the time-incremental minimization problems, 
\begin{equation*}
    \min_{(\boldsymbol\upsilon,\boldsymbol E^\textnormal{p})\in\mathcal Q}
    \Big\{
        \mathcal E_\varepsilon(t_\ts,\boldsymbol\upsilon,\boldsymbol E^\textnormal{p})
        +\mathcal D\big(\boldsymbol E^\textnormal{p}_{\ts-1},\boldsymbol E^\textnormal{p}\big)
    \Big\},
\end{equation*}
i.e., from the incremental energy minimization principle, cf.~\eqref{eq:MPDP} and the discussion therein. For $\varepsilon=\kappa=0$, and up to the reformulation of the Dirichlet data through $\boldsymbol y(t)=\boldsymbol g_t\circ\boldsymbol\upsilon(t)$, the optimality conditions of these incremental problems are the time-discrete equilibrium equations \eqref{weak2} with vanishing volume and surface forces, together with the local return mapping at the quadrature points. Hence, the algorithm of Section~\ref{sec:numerics} computes stationary points of the same incremental functionals without the regularizing terms; since Newton's method is a local method, these stationary points need not be global minimizers.

\emph{Spatial discretization.} On a fixed finite element mesh, the incremental problems are finite-dimensional. At fixed mesh size, existence of an incremental minimizer can be established under appropriate discrete coercivity assumptions (including a quadrature rule that controls the discrete deformation gradients), provided the admissible set is nonempty, and orientation preservation is enforced at the quadrature points. The logarithmic-strain energy has a singularity as an eigenvalue of $\F^\top\F$ approaches zero, and its growth in the other singular values controls large discrete gradients. The precise discrete coercivity also depends on the boundary constraints and the hardening term; these conditions should not be inferred from continuity of the energy alone. The corresponding bounds are, however, not uniform with respect to the mesh size $h$, and without regularization one cannot expect compactness of the discrete solutions as $h\to0$. This is where the regularizing terms become indispensable.

\emph{The limit $\varepsilon\to0$.} For a fixed deformation with integrable second gradient and Jacobian bounded away from zero, the regularizing terms contribute only $O(\varepsilon)$ to the energy, so the regularized and the non-regularized incremental functionals are close on such configurations. The behaviour of the energetic solutions as $\varepsilon\to0$ is, however, not covered by our analysis. Since the elastic energy is not quasiconvex, minimizing sequences of the limit problem may develop fine-scale oscillations, and the limit is expected to involve a relaxed energy. This is consistent with the loss of ellipticity of the logarithmic-strain model reported in \cite{FrMeSt2022} at excessive strains. A rigorous analysis of the joint limit $\varepsilon\to0$, $h\to0$ is left for future work.

\subsection{Vectorized assembly of transformation operators in the postprocessing phase}\label{app:PP-TL-code}

We provide a code block with details on the vectorized implementation of part \textbf{(Ca)} of the postprocessing phase in Fig.~\ref{fig:alg-setup}.

\begin{mm}[caption={Vectorized assembly of the transformation tensor $\PP$ and contraction $\Sigb : \LL$ from \eqref{eqsPPTL} and \eqref{eqsDF}.},label={code77}]
function [PP, SiLL_part] = assembly_PP_SiLL(F, stress, Eigs, EigVecs)

tol = 1e-10;  % tolerance under which we can consider two eigenvalues to be the same

% 1 6 5
% 9 2 4  (Voigt notation)
% 8 7 3


N = size(F,3);  % number of integration points

n = amam(F,EigVecs);  % 3x3xN (non-symmetric)
d = 1./Eigs;  % 3xN
epsi = log(Eigs)/2;  % 3xN
dzeta = amam(amtam(EigVecs,stress),EigVecs);  % 3x3xN (symmetric)

Eigs12 = Eigs(1,:) - Eigs(2,:);  % 1xN
Eigs13 = Eigs(1,:) - Eigs(3,:);  % 1xN
Eigs23 = Eigs(2,:) - Eigs(3,:);  % 1xN

d1 = d(1,:);  d2 = d(2,:);  d3 = d(3,:);  % 1xN
epsi1 = epsi(1,:);  epsi2 = epsi(2,:);  epsi3 = epsi(3,:);  % 1xN

inds1 = abs(Eigs12)>tol & abs(Eigs23)>tol & abs(Eigs13)>tol;  % (1) all eigenvalues are different
inds2 = abs(Eigs12)>tol & abs(Eigs23)>tol & abs(Eigs13)<tol;  % (2) l1 == l3, l2 is different
inds3 = abs(Eigs12)>tol & abs(Eigs23)<tol;                    % (3) l2 == l3, l1 is different
inds4 = abs(Eigs12)<tol & abs(Eigs13)>tol;                    % (4) l1 == l2, l3 is different
inds5 = abs(Eigs12)<tol & abs(Eigs13)<tol;                    % (5) l1 == l2 == l3

theta12 = zeros(1,N);  theta13 = zeros(1,N);  theta23 = zeros(1,N);  % 1xN
if any(inds1)
    theta12(inds1) = (epsi1(inds1) - epsi2(inds1))./Eigs12(inds1);
    theta13(inds1) = (epsi1(inds1) - epsi3(inds1))./Eigs13(inds1);
    theta23(inds1) = (epsi2(inds1) - epsi3(inds1))./Eigs23(inds1);
end
if any(inds2)
    % never occurs as the eigenvalues are always sorted
end
if any(inds3)
    theta12(inds3) = (epsi1(inds3) - epsi2(inds3))./Eigs12(inds3);
    theta13(inds3) = (epsi1(inds3) - epsi3(inds3))./Eigs13(inds3);
    theta23(inds3) = 0.5*d2(inds3);
end
if any(inds4)
    theta12(inds4) = 0.5*d1(inds4);
    theta13(inds4) = (epsi1(inds4) - epsi3(inds4))./Eigs13(inds4);
    theta23(inds4) = (epsi2(inds4) - epsi3(inds4))./Eigs23(inds4);
end
if any(inds5)
    theta12(inds5) = 0.5*d1(inds5);
    theta13(inds5) = theta12(inds5);
    theta23(inds5) = theta12(inds5);
end

scalA = @(x) reshape(x,1,[]);

EV11 = scalA(EigVecs(1,1,:));  EV12 = scalA(EigVecs(1,2,:));  EV13 = scalA(EigVecs(1,3,:));   %  1xN
EV21 = scalA(EigVecs(2,1,:));  EV22 = scalA(EigVecs(2,2,:));  EV23 = scalA(EigVecs(2,3,:));   %  1xN
EV31 = scalA(EigVecs(3,1,:));  EV32 = scalA(EigVecs(3,2,:));  EV33 = scalA(EigVecs(3,3,:));   %  1xN
n11 = scalA(n(1,1,:));  n12 = scalA(n(1,2,:));  n13 = scalA(n(1,3,:));   %  1xN
n21 = scalA(n(2,1,:));  n22 = scalA(n(2,2,:));  n23 = scalA(n(2,3,:));   %  1xN
n31 = scalA(n(3,1,:));  n32 = scalA(n(3,2,:));  n33 = scalA(n(3,3,:));   %  1xN

M1 = reshape(2*[n11.*EV11; n21.*EV21; n31.*EV31; n21.*EV31; n11.*EV31; n11.*EV21; ...
                n31.*EV21; n31.*EV11; n21.*EV11],1,9,N);   %  1x9xN
M2 = reshape(2*[n12.*EV12; n22.*EV22; n32.*EV32; n22.*EV32; n12.*EV32; n12.*EV22; ...
                n32.*EV22; n32.*EV12; n22.*EV12],1,9,N);   %  1x9xN
M3 = reshape(2*[n13.*EV13; n23.*EV23; n33.*EV33; n23.*EV33; n13.*EV33; n13.*EV23; ...
                n33.*EV23; n33.*EV13; n23.*EV13],1,9,N);   %  1x9xN
M4 = reshape([n12.*EV13+n13.*EV12;  n22.*EV23+n23.*EV22;  n32.*EV33+n33.*EV32;  ...
              n22.*EV33+n23.*EV32;  n12.*EV33+n13.*EV32;  n12.*EV23+n13.*EV22;  ...
              n32.*EV23+n33.*EV22;  n32.*EV13+n33.*EV12;  n22.*EV13+n23.*EV12],1,9,N);   %  1x9xN
M5 = reshape([n11.*EV13+n13.*EV11;  n21.*EV23+n23.*EV21;  n31.*EV33+n33.*EV31;  ...
              n21.*EV33+n23.*EV31;  n11.*EV33+n13.*EV31;  n11.*EV23+n13.*EV21;  ...
              n31.*EV23+n33.*EV21;  n31.*EV13+n33.*EV11;  n21.*EV13+n23.*EV11],1,9,N);   %  1x9xN
M6 = reshape([n11.*EV12+n12.*EV11;  n21.*EV22+n22.*EV21;  n31.*EV32+n32.*EV31;  ...
              n21.*EV32+n22.*EV31;  n11.*EV32+n12.*EV31;  n11.*EV22+n12.*EV21;  ...
              n31.*EV22+n32.*EV21;  n31.*EV12+n32.*EV11;  n21.*EV12+n22.*EV11],1,9,N);   %  1x9xN

%=======================
% Assembly of PP matrix
%=======================

PP11 =  0.5 * asam(d1,amam(reshape([EV11.*EV11;   EV21.*EV21;   EV31.*EV31;  ...
                                  2*EV21.*EV31; 2*EV11.*EV31; 2*EV11.*EV21],6,1,N), M1));  %  6x9xN
PP22 =  0.5 * asam(d2,amam(reshape([EV12.*EV12;   EV22.*EV22;   EV32.*EV32;  ...
                                  2*EV22.*EV32; 2*EV12.*EV32; 2*EV12.*EV22],6,1,N), M2));  %  6x9xN
PP33 =  0.5 * asam(d3,amam(reshape([EV13.*EV13;   EV23.*EV23;   EV33.*EV33;  ...
                                  2*EV23.*EV33; 2*EV13.*EV33; 2*EV13.*EV23],6,1,N), M3));  %  6x9xN
PP12 =  2 * asam(theta12, amam(reshape([EV11.*EV12;  EV21.*EV22;  EV31.*EV32;  ...
EV21.*EV32 + EV22.*EV31;  EV11.*EV32 + EV12.*EV31;  EV11.*EV22 + EV12.*EV21],6,1,N), M6)); %  6x9xN
PP13 =  2 * asam(theta13, amam(reshape([EV11.*EV13;  EV21.*EV23;  EV31.*EV33;  ...
EV21.*EV33 + EV23.*EV31;  EV11.*EV33 + EV13.*EV31;  EV11.*EV23 + EV13.*EV21],6,1,N), M5)); %  6x9xN
PP23 =  2 * asam(theta23, amam(reshape([EV12.*EV13;  EV22.*EV23;  EV32.*EV33;
EV22.*EV33 + EV23.*EV32;  EV12.*EV33 + EV13.*EV32;  EV12.*EV23 + EV13.*EV22],6,1,N), M4)); %  6x9xN

PP =  PP11 + PP12 + PP13 + PP22 + PP23 + PP33;   %  6x9xN

%==========================
% Assembly of Si:LL matrix
%==========================

if nargout==2

f = -2./Eigs.^2;  % 3xN
f1 = f(1,:);  f2 = f(2,:);  f3 = f(3,:);  % 1xN

ksi12 = zeros(1,N);  ksi13 = zeros(1,N);  ksi23 = zeros(1,N);  % 1xN
ksi21 = zeros(1,N);  ksi31 = zeros(1,N);  ksi32 = zeros(1,N);  % 1xN
eta = zeros(1,N);  % 1xN
if any(inds1)
    ksi12(inds1) = +(theta12(inds1) - 0.5*d2(inds1))./Eigs12(inds1);  ...
    ksi21(inds1) = -(theta12(inds1) - 0.5*d1(inds1))./Eigs12(inds1);  ...
    ksi13(inds1) = +(theta13(inds1) - 0.5*d3(inds1))./Eigs13(inds1);  ...
    ksi31(inds1) = -(theta13(inds1) - 0.5*d1(inds1))./Eigs13(inds1);  ...
    ksi23(inds1) = +(theta23(inds1) - 0.5*d3(inds1))./Eigs23(inds1);  ...
    ksi32(inds1) = -(theta23(inds1) - 0.5*d2(inds1))./Eigs23(inds1);
    eta(inds1) = epsi1(inds1)./Eigs12(inds1)./Eigs13(inds1) - ...
                 epsi2(inds1)./Eigs12(inds1)./Eigs23(inds1) + ...
                 epsi3(inds1)./Eigs13(inds1)./Eigs23(inds1);
end
if any(inds2)
    1;
end
if any(inds3)
    ksi12(inds3) = +(theta12(inds3) - 0.5*d2(inds3))./Eigs12(inds3);  ...
    ksi21(inds3) = -(theta12(inds3) - 0.5*d1(inds3))./Eigs12(inds3);  ...
    ksi13(inds3) = +(theta13(inds3) - 0.5*d3(inds3))./Eigs13(inds3);  ...
    ksi31(inds3) = -(theta13(inds3) - 0.5*d1(inds3))./Eigs13(inds3);  ...
    ksi23(inds3) = 0.125*f2(inds3);   ksi32(inds3) = ksi23(inds3);
    eta(inds3) = ksi12(inds3);  % equal to 'ksi13'
end
if any(inds4)
    ksi12(inds4) = 0.125*f1(inds4);   ksi21(inds4) = ksi12(inds4);  ...
    ksi13(inds4) = +(theta13(inds4) - 0.5*d3(inds4))./Eigs13(inds4);  ...
    ksi31(inds4) = -(theta13(inds4) - 0.5*d1(inds4))./Eigs13(inds4);  ...
    ksi23(inds4) = +(theta23(inds4) - 0.5*d3(inds4))./Eigs23(inds4);  ...
    ksi32(inds4) = -(theta23(inds4) - 0.5*d2(inds4))./Eigs23(inds4);  ...
    eta(inds4) = ksi31(inds4);  % equal to 'ksi32'
end
if any(inds5)
    ksi12(inds5) = 0.125*f1(inds5);  ksi21(inds5) = ksi12(inds5);  ...
    ksi13(inds5) = ksi12(inds5);     ksi31(inds5) = ksi12(inds5);  ...
    ksi23(inds5) = ksi12(inds5);     ksi32(inds5) = ksi12(inds5);  ...
    eta(inds5) = 0.125*f1(inds5);
end

dzeta11 = scalA(dzeta(1,1,:));  dzeta22 = scalA(dzeta(2,2,:));  dzeta33 = scalA(dzeta(3,3,:)); % 1xN
dzeta12 = scalA(dzeta(1,2,:));  dzeta13 = scalA(dzeta(1,3,:));  dzeta23 = scalA(dzeta(2,3,:)); % 1xN


scalB = @(x) reshape(x,1,1,[]);
%% --- SiLL123 -----------------------------------------------------------
w1 = scalB(f1 .* dzeta11);
w2 = scalB(f2 .* dzeta22);
w3 = scalB(f3 .* dzeta33);
%% --- SiLL4 -------------------------------------------------------------
A4 = ...
    M4 .* scalB(ksi23.*dzeta33 + ksi32.*dzeta22) + ...
    2 * scalB(dzeta23) .* ( scalB(ksi23).*M3 + scalB(ksi32).*M2 ) + ...
    2 * scalB(eta.*dzeta13) .* M6;   %  1x9xN
%% --- SiLL5 -------------------------------------------------------------
A5 = ...
    M5 .* scalB(ksi13.*dzeta33 + ksi31.*dzeta11) + ...
    2 * scalB(dzeta13) .* ( scalB(ksi13).*M3 + scalB(ksi31).*M1 ) + ...
    2 * scalB(eta.*dzeta12) .* M4;   %  1x9xN
%% --- SiLL6 -------------------------------------------------------------
A6 = ...
    M6 .* scalB(ksi12.*dzeta22 + ksi21.*dzeta11) + ...
    2 * scalB(dzeta12) .* ( scalB(ksi12).*M2 + scalB(ksi21).*M1 ) + ...
    2 * scalB(eta.*dzeta23) .* M5;   %  1x9xN
%%
SiLL = 0.125*(amtam(M1.*w1,M1) + amtam(M2.*w2,M2) + amtam(M3.*w3,M3)) + ...
            amtam(A4, M4) + amtam(A5, M5) + amtam(A6, M6);
SiLL_part =  SiLL + amt(SiLL);  % 9x9xN

end

end
\end{mm}

\subsection{Distribution of von Mises plastic strain within cross-section of the coil from Benchmark~\ref{subsubsec:benchmark2}}
\label{app:coil-cross-plast-strain}

\begin{figure}
    \centering
    \includegraphics[width=0.82\textwidth]{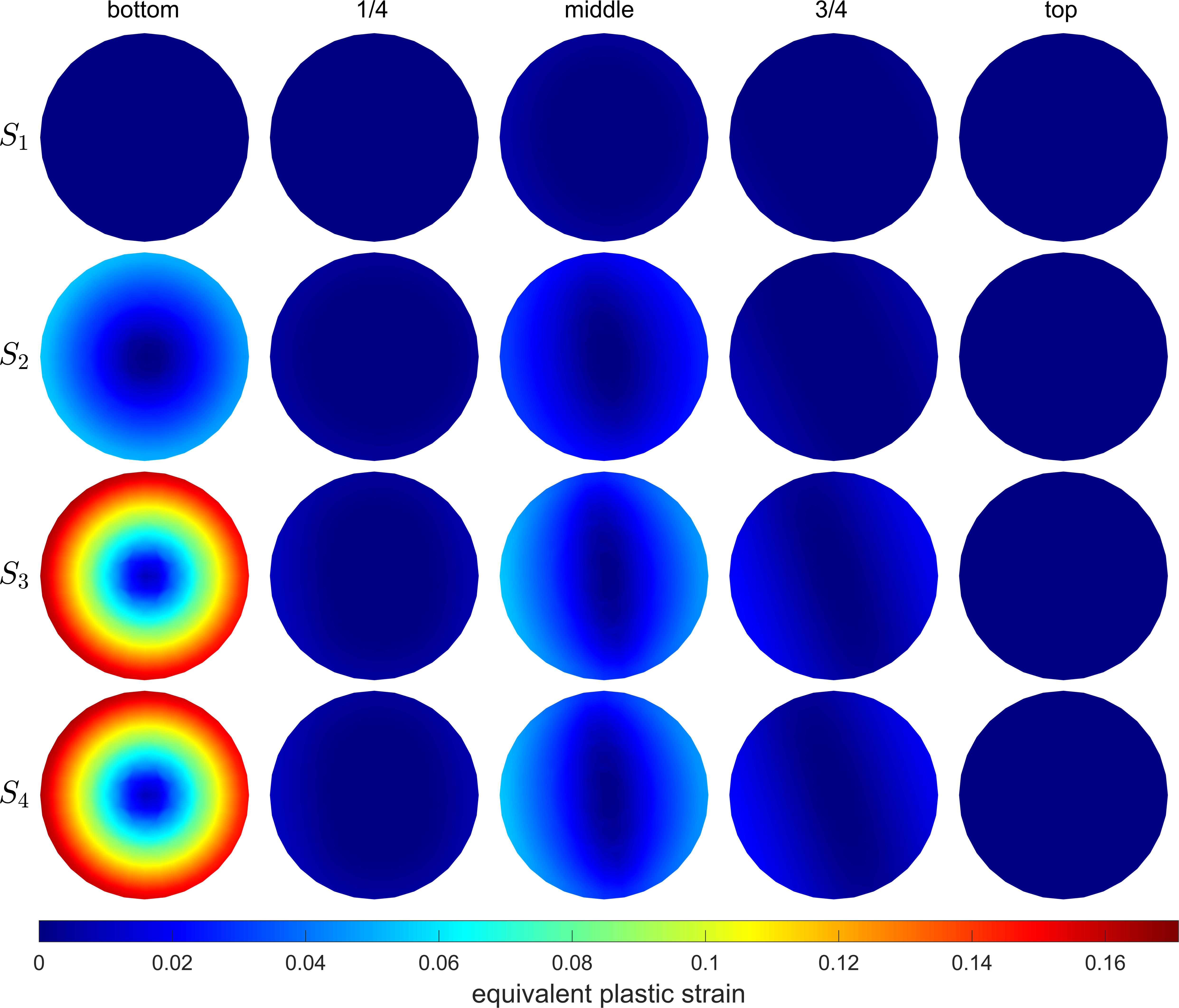}
    \caption{Benchmark~\ref{subsubsec:benchmark2}: equivalent plastic strain on particular coil cross-sections for stages of stretching $S_1-S_4$. See text for details.}
    \label{fig:coil-cross-plast-strain}
\end{figure}

Fig.~\ref{fig:coil-cross-plast-strain} provides the distribution of equivalent plastic strain through particular coil cross-sections for stages of stretching $S_1-S_4$, see Table~\ref{tab:snapshots}. The positions of cross-sections correspond to the relative distances 0 (=bottom), 1/4, 1/2 (=middle), 3/4, and 1 (=top) along the central fiber of the coil starting from the fixed end.

The plots confirm that up to stretching $S_1$, most of the material responds elastically; just the surface of the middle cross-section starts to plastify. Consistent with the increase in von Mises stress in Fig.~\ref{fig:coil-cross-stress}, the bottom cross-section shows the most pronounced increase in plasticity at $S_2$. At stretch $S_3$, the plastification becomes even more pronounced; the maximum equivalent plastic strain reaches $\approx 0.17$. Yet, the radial symmetry is preserved because of the dominant torsional loading, cf. Fig.~\ref{fig:coil-cross-stress}. However, such a symmetry is broken, e.g., for the middle cross-section, where the torsion is combined with bending. No plastic deformation has been initiated within the top cross-section. During unloading to $S_4$, the distribution of plastic strain is unaltered, which corresponds to purely elastic unloading.

\end{document}